\makeatletter
\@ifundefined{pdfpageheight}{\newdimen\pdfpageheight \pdfpageheight=\paperheight}{}
\@ifundefined{pdfpagewidth}{\newdimen\pdfpagewidth  \pdfpagewidth=\paperwidth}{}
\makeatother

\documentclass{article}
\usepackage{arxiv}
\usepackage{tikz}
\usepackage{subcaption}
\usepackage{graphicx,xcolor,xspace}
\usepackage{mathtools,amssymb,amsfonts,amsthm}
\usepackage[utf8]{inputenc} 
\usepackage[T1]{fontenc}    
\usepackage{booktabs}
\usepackage{makecell}
\usepackage{multirow}
\usepackage{comment}
\usepackage[colorlinks=true]{hyperref}       
\usepackage[nameinlink,capitalise,noabbrev]{cleveref}

\AddToHook{cmd/appendix/before}{%
    \crefalias{section}{appendix}%
    \crefalias{subsection}{appendix}
}

\newcommand{\linkref}[2]{\hyperref[#2]{#1~\ref{#2}}}

\newcommand{\propref}[1]{\linkref{Proposition}{#1}}

\newcommand{\corref}[1]{\linkref{Corollary}{#1}}

\usetikzlibrary{decorations.pathmorphing,patterns}
\newfont{\notapolice}{cmss8}
\newfont{\rempolice}{cmss9}
\newfont{\tablepolice}{cmtt10}

\newcommand{\Ham}[2]{{\mathcal H\left ({#1},{#2}\right )}}

\newcommand{\Kbody}{{\mathfrak B}}
\newcommand{\ba}{{\boldsymbol{a}}}
\newcommand{\bb}{{\boldsymbol{b}}}
\newcommand{\bx}{{\boldsymbol{x}}}
\newcommand{\bp}{{\boldsymbol{p}}}
\newcommand{\bv}{{\boldsymbol{v}}}
\newcommand{\trans}{{\textrm{t}}}
\newcommand{\blockZ}{{\boldsymbol \flat\hskip-0.2em \boldsymbol Z}}
\newcommand{\blockD}{{\boldsymbol \flat\hskip-0.15em\boldsymbol D}}
\newcommand{\blockS}{{\boldsymbol \flat\hskip-0.15em\boldsymbol S}}
\newcommand{\Zx}{{Z\hskip-0.1em\textrm{x}}}
\newcommand{\Dx}{{D\hskip-0.1em\textrm{x}}}
\newcommand{\Sx}{{S\hskip-0.1em\textrm{x}}}
\newcommand{\Zp}{{Z\hskip-0.1em\textrm{p}}}
\newcommand{\Dp}{{D\hskip-0.1em\textrm{p}}}
\newcommand{\Sp}{{S\hskip-0.1em\textrm{p}}}
\newcommand{\bZx}{{Z\hskip-0.1em\textrm{\bf x}  }}
\newcommand{\bDx}{{D\hskip-0.1em\textrm{\bf x}  }}
\newcommand{\bSx}{{S\hskip-0.1em\textrm{\bf x}  }}
\newcommand{\bZp}{{Z\hskip-0.1em\textrm{\bf p}  }}
\newcommand{\bDp}{{D\hskip-0.1em\textrm{\bf p}  }}
\newcommand{\bSp}{{S\hskip-0.1em\textrm{\bf p}  }}
\newcommand{\ZX}{{\boldsymbol{Z}\hskip-0.1em\textrm{\bf x}  }}
\newcommand{\DX}{{\boldsymbol{D}\hskip-0.1em\textrm{\bf x}  }}
\newcommand{\SX}{{\boldsymbol{S}\hskip-0.1em\textrm{\bf x}  }}
\newcommand{\ZP}{{\boldsymbol{Z}\hskip-0.1em\textrm{\bf p}  }}
\newcommand{\DP}{{\boldsymbol{D}\hskip-0.1em\textrm{\bf p}  }}
\newcommand{\SP}{{\boldsymbol{S}\hskip-0.1em\textrm{\bf p}  }}
\newcommand{\blockZx}{{\blockZ\hskip-0.1em\textrm{x}}}
\newcommand{\blockDx}{{\blockD\hskip-0.1em\textrm{x}}}
\newcommand{\blockSx}{{\blockS\hskip-0.1em\textrm{x}}}
\newcommand{\blockZp}{{\blockZ\hskip-0.1em\textrm{p}}}
\newcommand{\blockDp}{{\blockD\hskip-0.1em\textrm{p}}}
\newcommand{\blockSp}{{\blockS\hskip-0.1em\textrm{p}}}
\newcommand{\blockZX}{{\blockZ\hskip-0.1em\textrm{\bf x}}}
\newcommand{\blockDX}{{\blockD\hskip-0.1em\textrm{\bf x}}}
\newcommand{\blockSX}{{\blockS\hskip-0.1em\textrm{\bf x}}}
\newcommand{\blockZP}{{\blockZ\hskip-0.1em\textrm{\bf p}}}
\newcommand{\blockDP}{{\blockD\hskip-0.1em\textrm{\bf p}}}
\newcommand{\blockSP}{{\blockS\hskip-0.1em\textrm{\bf p}}}
\newcommand{\SE}{{\texttt{SE}\xspace }}
\newcommand{\PE}{{\texttt{PE}\xspace }}
\newcommand{\ZD}{{\texttt{ZD}\xspace }}
\newcommand{\ZDS}{{\texttt{ZDS}\xspace }}
\newcommand{\eH}{{\texttt{eH}\xspace }}
\newcommand{\ordH}{{\texttt{ordH}\xspace }}
\newcommand{\ex}{{\texttt{ex}\xspace }}
\newcommand{\ordx}{{\texttt{ordx}\xspace }}
\newcommand{\ea}{{\texttt{eA}\xspace }}

\newcommand{\eL}{{\texttt{eL}\xspace }}

\newcommand{\eq}{{\texttt{eq}\xspace }}
\newcommand{\ordq}{{\texttt{ordq}\xspace }}
\newcommand{\avrOne}{{ {\tt\small nb\_iter\_avg} }}
\newcommand{\avrTwo}{{ {\tt\small nb\_call\_avg} }}
\newcommand{\One}{{{\boldsymbol e}}}

\definecolor{review1}{RGB}{117,112,179}
\definecolor{review2}{RGB}{217,95,2}
\definecolor{review3}{RGB}{27,158,119}
\definecolor{allreviews}{RGB}{231,41,138}
\definecolor{todiscuss}{RGB}{127,150,127}

\newcommand{\ra}{\textcolor{black}}
\newcommand{\rb}{\textcolor{black}}

\newcommand{\rall}{\textcolor{black}}
\newcommand{\rallparagraph}{\color{black}}

\newcommand{\rabis}{\textcolor{review1}}
\newcommand{\rbbis}{\textcolor{review2}}

\newcommand{\rbparagraphbis}{\color{review2}}

\newtheorem{proposition}{Proposition}
\newtheorem{corollary}{Corollary}
\newtheoremstyle{myremark}
{3pt}
{3pt}
{\rempolice}
{}
{\itshape}
{:}
{.5em}
{}
\theoremstyle{myremark}
\newtheorem*{remark}{Remark}

\title{Structural schemes for Hamiltonian systems}

\author{
    \href{https://orcid.org/0000-0003-2295-5118}{\includegraphics[scale=0.06]{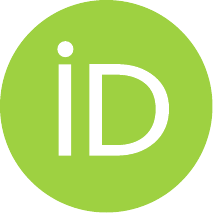}\hspace{1mm}Stéphane Clain} \\
	Centre of Mathematics, Coimbra University, Largo D. Dini, 3000-143 Coimbra, Portugal \\
	\texttt{clain@mat.uc.pt} \\
	\And
	\href{https://orcid.org/0000-0001-6180-9485}{\includegraphics[scale=0.06]{orcid.pdf}\hspace{1mm}Emmanuel Franck} \\
	Université de Strasbourg, CNRS, Inria, IRMA, F-67000, Strasbourg, France \\
	\texttt{emmanuel.franck@inria.fr} \\
	\And
	\href{https://orcid.org/0000-0002-3859-8517}{\includegraphics[scale=0.06]{orcid.pdf}\hspace{1mm}Victor Michel-Dansac} \\
	Université de Strasbourg, CNRS, Inria, IRMA, F-67000, Strasbourg, France \\
	\texttt{victor.michel-dansac@inria.fr} \\
}

\renewcommand{\headeright}{Article}
\renewcommand{\undertitle}{Article}
\renewcommand{\shorttitle}{Structural schemes for Hamiltonian system}

\begin{document}
\maketitle
\begin{abstract}
    We present an adaptation of the so-called structural method \cite{CMM23} for Hamiltonian systems, and redesign the method for this specific context, which involves two coupled differential systems. Structural schemes decompose the problem into two sets of equations: the physical equations, which describe the local dynamics of the system, and the structural equations, which only involve the discretization of a very compact stencil. \rall{They have desirable properties, such as unconditional stability and high-order accuracy \cite{CMMC25}}. We first give a general description of the scheme for the \rb{one-body} case (which corresponds to e.g. spring-mass interactions or pendulum motion), before extending the technique to the \rb{multi-body} case (treating e.g. the $n$-body system). The scheme is also written in the case of a non-separable system (e.g. a charged particle in an electromagnetic field). \rall{We prove that the method is exactly energy-preserving for quadratic Hamiltonians, and} we give numerical evidence of the method's efficiency, its capacity to preserve invariant quantities such as the total energy, and draw comparisons with the traditional symplectic methods.
\end{abstract}

\section{Introduction}
Hamiltonian systems, a class of ordinary differential equations (ODEs) that arise from Hamiltonian mechanics, play a fundamental role in the mathematical modelling of physical systems with conserved quantities, such as total energy for a closed system. The dynamics of such systems are described by the principle of stationary action, leading to the Euler-Lagrange equations and the associated Lagrangian function \cite{LeRe2004}.
Derived from the Hamiltonian function, Hamiltonian systems govern the time evolution of state variables such as positions and momenta.
Hamiltonian systems are ubiquitous in various fields, including celestial mechanics, where they describe planetary orbits, or plasma physics, capturing the behavior of charged particles. Their inherent structure, characterized by symplectic geometry and energy preservation, has to be considered in their numerical approximation.

Due to their rich structure, the numerical approximation of Hamiltonian systems presents several challenges. One of the main issues is the preservation of geometric features, such as symplecticity and conserved quantities like energy or angular momentum, which are integral to the physical behavior of the system.
Standard numerical methods, even of very high order, such as explicit or implicit Runge-Kutta schemes, often fail to maintain these properties during simulations, generating error accumulations for long run simulations that can become critical and totally spoil the quality of the solution \cite{Beust2003}. Additionally, the typically oscillatory nature of the solutions to such systems, particularly in applications like plasma physics, require specific methods to enable of handling multiple time scales without compromising accuracy or stability. On the other hand, the high-dimensional phase spaces involved in many problems further exacerbate computational costs, requiring efficient algorithms that balance precision with practicality.

These challenges have motivated the development of structure-preserving methods, such as symplectic integrators and geometric numerical schemes \cite{KaMe1987}, that ensure robust long-term stability by providing a good approximation of the total energy or other invariants. The design of schemes that preserve geometric structures for ODEs provide accurate and stable numerical methods. Some classes of partial differential equations, after well-adapted spatial discretization, are written as high-dimensional Hamiltonian ODEs. As examples, we mention the works of \cite{gempic1,gempic2} in plasma physics or \cite{shepherd1990symmetries} in  climatology where these approaches, using geometric structure-preserving time discretization, have shown their superiority.

Hamiltonian systems typically involve a system of two coupled ODEs, one describing the time evolution of a position in phase space, and the second one governing the time evolution of the momentum (or velocity). The simplest symplectic algorithm is a semi-implicit Euler integrator, where the first equation is integrated explicitly, and the second one implicitly. Therefore, it requires solving a (potentially non-linear) implicit equation at each time step.
However, this integrator is only first-order accurate. In fact, one of the most well-known symplectic algorithms for Hamiltonian systems is the \rb{St\"ormer}-Verlet method, see~\cite{Ver1967,HaiLubWan2003}, which can be seen as a second-order accurate version of the semi-implicit Euler scheme.
From this second-order algorithm, it is relatively easy to construct high-order methods thanks to the composition strategy, introduced in \cite{Suz1990,Yos1990,McL1995} and summarized in \cite{HaiLubWan2006}. Nevertheless, one drawback of this approach is that the number of function evaluations increases exponentially with the order of accuracy. Other examples of high-order schemes are provided in \cite{McLAte1992,SanCal1993}, and other composition schemes include the ones from \cite{KahLi1997}, which still \ra{suffer} from the same drawback in terms of number of function evaluations.

The structural method was recently introduced, in \cite{CMM23,CMMC25}, to solve ODE systems. It provides a systematic way of constructing implicit schemes with a highly compact stencil, which have an arbitrarily high order of accuracy together with unconditional stability. The technique relies on separating the physical part of the problem, i.e., the system of ODEs, from the discretization, which only involves the unknown quantities associated with the grid. In practice, the data is gathered into blocks, corresponding to several time steps to be solved simultaneously. This methodology was also adapted to one-dimensional boundary-valued problems in \cite{CPPL2023}. In this case, the structural method achieves a sixth-order accurate approximations using a tiny stencil of three points. \rall{Very recently, a complete study of the structural method for ODE systems \cite{CMMC25} has been carried out where linear stability, convergence, numerical errors and spectral resolution \rbbis{were analyzed} together with a wide spectrum of benchmarks to assess the accuracy and stability and the ability to handle stiff problems.}

In this paper, we develop and adapt the structural method to Hamiltonian systems, considering the specificities of such systems, namely the conservation of several invariants, and the two-variable nature of the problem. The paper is organized as follows. First, we present the structural method for ordinary differential equations in \cref{sec:structural_method_ODEs}, introducing two formulations, the first using only the ODE and the second adding information related to the time derivative of the ODE. Then, the structural method is applied to Hamiltonian \rb{one-body} problems in \cref{sec:Hamiltonian_scalar} and to general \rb{multi-body} Hamiltonian systems in \cref{sec:Hamiltonian_vector}, thus constructing high-order accurate and stable schemes for Hamiltonian systems. Finally, we present numerical results in \cref{sec:benchmarks} to demonstrate the efficiency and accuracy of the proposed schemes. Several Hamiltonian systems are considered, from simple separable examples to more complex non-separable ones, to illustrate the performance of the structural method in various contexts.

\section{The structural method for ODEs}
\label{sec:structural_method_ODEs}

We provide a short review of the structural method based on the so-called Physical and Structural Equations (denoted \PE\ and \SE\ respectively) detailed in \cite{CMM23,CMMC25}. To this end, let us consider the Ordinary Differential Equation (ODE) for $t\in [0,T]$
\begin{equation}\label{eq::EDO1}
    \dot x=f(x),\qquad x(0)=x_0.
\end{equation}
Traditional schemes blend the discretization with the physical equation, that is, the function that describe the dynamic of the physical system. For example, the popular Crank-\ra{Nicolson} scheme reads
$$
    \frac{x_{n+1}-x_n}{\Delta t}=\frac{f(x_{n})+f(x_{n+1})}{2}
$$
where the left-hand side is the time derivative discretization, while the right-hand side represents the physics (i.e., the function $f$ that characterizes the physical problem). Obviously, one can split the scheme into two equations, namely
$$
    D_{n+1}=f(Z_{n+1})
    \text{\qquad and \qquad}
    2(Z_{n+1}-Z_{n})-\Delta t(D_{n}+D_{n+1})=0,
$$
where $Z_n$ and $D_n$ are approximations of the Zeroth-order derivative $x(t_n)$ and the first-order Derivative  $\dot x(t_n)$, applying the same notations for $Z_{n+1}$ and $D_{n+1}$ at the time $t_{n+1}$.
The left relation is called the Physical Equation since it only involves the physics of the problem and not the discretization; the right one is called the Structural Equation since it only depends on the grid structure, and involves no physics.

An extension relies on taking the time derivative of the ODE, which reads $\ddot x=f'(x)\dot x$. Denoting by $S_n$ an approximation of the Second derivative $\ddot x(t_n)$, we now obtain two physical equations, respectively denoted by $\PE[1]$ and $\PE[2]$:
$$
    D_{n+1}=f(Z_{n+1})
    \text{\qquad and \qquad}
    S_{n+1}=f'(Z_{n+1})D_{n+1}.
$$
Since we have three unknowns $(Z_{n+1},D_{n+1},S_{n+1})$, we need one additional equation that represents the discretization.
We choose the following relation, exact for polynomial functions up to degree 4, and called the structural equation $\SE[1]$:
$$
    12(Z_{n+1}-Z_{n})-6\Delta t(D_{n+1}+D_n)+\Delta t^2(S_{n+1}-S_n)=0.
$$
Given the values $(Z_{n+1},D_{n+1},S_{n+1})$ at a time $t_n$, we seek solutions at $t_{n+1}$ of $\PE[1]$, $\PE[2]$, and $\SE[1]$. Note that the problem is fully implicit and requires solving a small nonlinear system as soon as $f$ is nonlinear.

Following the example, the idea of structural method consists in splitting the problem with, on the one hand, the  Physical Equations and, \ra{on} the other hand, the Structural Equations. The two sets of equations involve the unknowns over a block of size $R$ steps corresponding to the time step $t_{n+1}$ until $t_{n+R}$, given the initial configuration at the time $t_n$. We obtain a non-linear system combining the function approximations together with the derivatives for all the  intermediate steps.

We now detail the structural method for the scalar ODE \eqref{eq::EDO1} in the \ZD\ and \ZDS\ formulations, respectively involving the zeroth- and first-order derivatives in \cref{sec:ZD_generic}, and the zeroth-, first- and second-order derivatives in \cref{sec:ZDS_generic}.
The complete algorithm is then given, for each formulation, in \cref{sec:full_algo_scalar}.

\subsection{\ZD\ equations}
\label{sec:ZD_generic}

We first tackle the simple version called the \ZD\ scheme, where one only uses implicit combinations of the approximate function and first derivative as unknowns. The generic \ZD\ structural equation then reads
\begin{equation}\label{ZD_structural_equation}
    \sum_{r=0}^R a_{r,0} Z_{n+r}+a_{r,1} \Delta t\, D_{n+r}=0,
\end{equation}
where $(a_{r,s})_{r\in\{0,\dots,R\}, s\in\{0,1\}}$ are the coefficients of the \SE, independent of $n$ by assuming a uniform time discretization with time step $\Delta t$. It is important to remark that $r=0$ corresponds to the time $t_n$ where all the variables are given. In total, $2(R+1)$ coefficients characterize the  structural equation. They can be reshaped as a vector
$$
    \ba=\begin{bmatrix}
        a_{0,0}, & a_{1,0}, & \dots, & a_{R,0}, & a_{0,1}, & a_{1,1}, & \dots, & a_{R,1}
    \end{bmatrix}.
$$
On the other hand, the implicit problem involves $2R$ unknowns with $R$ physical equations $\PE[1], \dots, \PE[R]$ corresponding to the relations
$$
    D_{n+r}=f(Z_{n+r}), \text{\qquad for all } r \in \{1,\dots,R\}.
$$
Hence, the structural method requires $R$ linearly independent structural equations, i.e. $R$ vectors, $(\ba^m)_{m \in \{1,\dots,R\}}$ whose components are denoted by $(a^m_{r,s})_{r,s,m}$, to close the system of size $2R$.

We wish to provide this set of $R$ structural equations
by ensuring that the resulting scheme is high-order accurate.
To this end, for any vector $\ba$ and function $\phi$, we define the functional
\begin{equation}\label{functional_SE}
    E(\ba,\phi)=\sum_{r=0}^R a_{r,0} \phi(t_{r})+a_{r,1}\Delta t\, \phi'(t_{r}).
\end{equation}
Taking the particular case of polynomial functions $\displaystyle \pi_{\ell}(t)=\left (\frac{t}{\Delta t}\right )^{\ell-1}$ for $\ell \in \{1,\dots,2(R+1)\}$, we consider the linear system $E(\ba,\pi_\ell)=0$, which rewrites as the $2(R+1)$ equations
\begin{eqnarray*}
    \forall \ell \in \{1,\dots,2(R+1)\}, \quad
    E(a;\pi_{\ell})
    =
    \sum_{r=0}^R
    a_{r,0} \pi_{\ell}(r \Delta t) +
    a_{r,1}\Delta t\, \pi_{\ell}'(r\Delta t) = 0.
\end{eqnarray*}
We obtain a $2(R+1)\times 2(R+1)$ non-singular linear system given by $M\ba=0$, where the matrix $M$ gathers all the coefficients
$\pi_{\ell}(r \Delta t)$ and $\pi_{\ell}'(r \Delta t)$,
for all $\ell \in \{ 1,\dots,2(R+1)\}$ and $r\in\{0,\dots,R\}$.
\rall{
    \begin{remark}
        Notice that by construction the entries of matrix $M$ do not depend on the time step $\Delta t$, therefore the vectors of the kernel are independent of $\Delta t$.
    \end{remark}
}


To determine the $R$ structural equations, we simply withdraw the $R$ last lines of the matrix $M$, leading to a reduced matrix whose kernel has dimension $R$.
To build the structural equations,
we select an orthogonal basis of this kernel
and place it in the vector $(\ba^m)_{m \in \{1,\dots,R\}}$.
Thanks to this procedure, the obtained relations
\begin{equation}\label{R_ZD_structural_equation}
    \sum_{r=0}^R a^m_{r,0} Z_{n+r}+a^m_{r,1} \Delta t\,D_{n+r}=0,\quad m=1,\cdots,R,
\end{equation}
remain exact for polynomials up to degree $R+1$.
\rall{
    \begin{remark}
        The relation \eqref{R_ZD_structural_equation} does not correspond to a multistep method. Indeed, only  $Z_{n}$ and $D_n$ are given while $(Z_{n+r},D_{n+r})$, $r=1,\dots,R$ are the unknowns. To close the system, we impose the $R$ physical equations at the  time steps $t_{n+r}$, $r=1,\dots,R$.
    \end{remark}
}
\begin{remark}
    For small values of $R$ (say $R=2$), we derive an explicit analytic expression for the structural equation.
    However, this becomes intractable for larger $R>4$, and so we compute the kernel and determine an orthogonal basis to automatically provide the $R$ structural equations (see \cite{CMM23,CMMC25} for details).
\end{remark}
At the end of the day, we get a set of $R$ structural equations $\SE[1], \dots, \SE[R]$, linearly independent, exact for polynomials up to degree $R+1$.

\subsection{\texorpdfstring{$\ZDS$}{ZDS} equations}
\label{sec:ZDS_generic}

The \ZD\ scheme is quite effective and provides accurate solutions, but one can \ra{improve the accuracy of} the method by introducing the second derivatives.
Indeed, differentiating the original physical equation provides a second physical equation. Therefore, the scheme is more compact and accurate for the same block size.
More precisely, the \ZDS\ scheme involves the zeroth-, first- and second-order derivatives as unknowns, with the structural equation reading
\begin{equation}\label{ZDS_structural_equation}
    \sum_{r=0}^R a_{r,0} Z_{n+r}+a_{r,1}\Delta t\, D_{n+r}+a_{r,2} \Delta t^2\, S_{n+r}=0,
\end{equation}
with $\ba=(a_{r,s})_{r,s}$ its coefficients, independent of the time index $n$, the block number $r$ and the derivation order $s$.

The \ZDS\ method involves the $3R$ unknowns $(Z_{n+r},D_{n+r},S_{n+r})$ with $r \in \{1, \dots, R\}$.
The physical equations provide~$2R$ relations between these unknowns:
$$
    \forall r \in \{1, \dots, R\}, \qquad
    D_{n+r}=f(Z_{n+r})
    \text{\quad and \quad}
    S_{n+r}=f'(Z_{n+r})D_{n+r}.
$$
To close the system, we need $R$ linearly independent structural equations of type \eqref{ZDS_structural_equation} with respective coefficients~$(\ba^m)_{m \in \{1,\dots,R\}}$,
whose entries are $a^m_{r,d}$ with $r \in \{0,\dots, R\}$ and $d\in\{0,1,2\}$.

To this end, we introduce the functional
$$
    E(\ba,\phi)=\sum_{r=0}^R a_{r,0} \phi(r \Delta t)+a_{r,1}\Delta t\, \phi'(r\Delta t)+a_{r,2}\Delta t^2\, \phi''(r\Delta t)
$$
and recall the set of polynomial functions $\pi_{\ell}=t^{\ell-1}$. Then the $3(R+1)\time 3(R+1)$ linear system $E(\ba,\pi_{\ell})=0$ is non-singular and reads $M\ba=0$ where $M$ is an invertible matrix.
\rall{
    \begin{remark}
        Once again the entries of matrix $M$ do not depend on $\Delta t$, hence the vectors of the kernel are independent of $\Delta t$.
    \end{remark}
}
Eliminating the $R$ last rows of the matrix $M$, we get a kernel of dimension $R$ with an orthogonal basis $(\ba^m)_{m \in \{1,\dots,R\}}$ that provides the $R$ structural equation $\PE[1],\dots,\PE[R]$. Note that the relations
\begin{equation}\label{R_ZDS_structural_equation}
    \sum_{r=0}^R a^m_{r,0} Z_{n+r}+a^m_{r,1}\Delta t\, D_{n+r}+a^m_{r,2}\Delta t^2\, S_{n+r}=0
    \text{\qquad for } m \in \{1,\dots,R\}
\end{equation}
are exact for all polynomials of degree lower than $2(R+1)$. The solution of the $2R$ physical equations together with the $R$ structural equations provide an approximation of the $3R$ variables.
\rall{
    \begin{remark}
        The relations \eqref{R_ZDS_structural_equation} do not correspond to a two-derivative multistep method in the sense of the book of Butcher \cite[Chapter 4]{B08}, since only $(Z_n,D_n,S_n)$ are given whereas the other variables are unknowns of the problem. We use the expression $R$-block to highlight that we solve a non-linear system that involves all the unknowns from $t_{n+1}$ until $t_{n+R}$.
    \end{remark}
}
%
%
\subsection{Structural method for a scalar ODE}
\label{sec:full_algo_scalar}

Equipped with the physical and structural equations in the \ZD\ and \ZDS\ formulations, we now detail the full algorithm to compute a full $R$-block of solutions with each of these formulations.

\subsubsection{\ZD\ formulation}
\label{sec:full_algo_scalar_ZD}
Assume that we know the solution $(Z_n,D_n)$ at the time $t_n$. We seek the solution of the nonlinear problem with $2R$ unknowns
\begin{equation}
    \label{ZD_scalar_problem}
    \forall m \in \{1, \dots, R\}, \qquad
    D_{n+m}=f(Z_{n+m})
    \text{\quad and \quad}
    \sum_{r=0}^R a^m_{r,0} Z_{n+r}+a^m_{r,1}\Delta t\, D_{n+r}=0.
\end{equation}

Let $\displaystyle \blockZ_n=(Z_{n+r})_{r\in\{1,\dots,R\}}\in\mathbb R^R$ and $\displaystyle \blockD_n=(D_{n+r})_{r\in\{1,\dots,R\}}\in\mathbb R^R$ be the $R$-blocks for the approximations of $x$ and $\dot x$ at time $t_{n+r}$, for each $r \in \{1,\dots,R\}$. We rewrite the structural equations under a matrix form
\begin{equation}
    \label{eq:ZD_scheme_matrix_form}
    A_z\, \blockZ_n+\Delta t\, A_d\, \blockD_n+Z_n\, \ba_z+\Delta t\, D_n\, \ba_d=0,
\end{equation}
with
\begin{align*}
    A_z=(a^m_{r,0})_{m,r}, &
    \qquad \ba_z=\big (a^1_{0,0},a^2_{0,0},\cdots,a^R_{0,0}\big )^\trans, \\
    A_d=(a^m_{r,1})_{m,r}, &
    \qquad \ba_d=\big (a^1_{0,1},a^2_{0,1},\cdots,a^R_{0,1}\big )^\trans.
\end{align*}
Assuming the matrix $A_z\in\mathbb R^{R\times R}$ is non-singular, we  rewrite the system as
$$
    \blockZ_n+\Delta t\,B_d\,\blockD_n+Z_n\, \bb_z+\Delta t\,D_n\, \bb_d=0
$$
with $B_d=(A_z)^{-1}A_d$, $\bb_z=(A_z)^{-1}\ba_z$, $\bb_d=(A_z)^{-1}\ba_d$.

To solve problem \eqref{ZD_scalar_problem}, we proceed with a fixed-point method by producing a sequence $(\blockZ_n^{(k)},\blockD_n^{(k)})$ that converges to the solution. We sketch the algorithm hereafter.
\begin{itemize}
    \item {\bf Initialization.} To build $\blockZ_n^{(0)}$ and $\blockD_n^{(0)}$, we set for $r \in \{1,\dots,R\}$
          $$
              Z^{(0)}_{n+r}=Z^{(0)}_{n+r-1}+\Delta t D^{(0)}_{n+r-1},\quad D^{(0)}_{n+r}=f(Z^{(0)}_{n+r}),
          $$
          with $Z^{(0)}_{n}=Z_n$, $D^{(0)}_{n}=D_n$.
    \item {\bf Iteration.} We first compute a new approximation for the solution $\blockZ_n$ using the structural equations:
          $$
              \blockZ_n^{(k+1)}=-\Big (Z_n\,\bb_z+\Delta t\,D_n\,\bb_d+\Delta t\,B_d\,\blockD_n^{(k)}\Big )
          $$
          Then, we update the first derivative with the physical equations:
          $$
              \forall r \in \{1, \dots, R\}, \qquad
              D^{(k+1)}_{n+r}=f(Z^{(k+1)}_{n+r}).
          $$
    \item {\bf Stopping criterion.} We end the fixed-point when two successive solutions are close enough, according to the tolerance parameter \texttt{tol}, that is $\lVert \blockZ_n^{(k+1)}-\blockZ_n^{(k)}\rVert \leq \texttt{tol}$.
\end{itemize}
\begin{remark}
    Notice that $R=2$ provides a 4th-order unconditionally stable scheme while $R=4$ reaches sixth-order accuracy. Methods with $R=6$ and $R=8$ will also be experimented in the numerical section.\qed
\end{remark}
\begin{remark}
    The non-linear part of the problem does not require some local linearization, hence the computation of the derivative is straightforward. $\qed$
\end{remark}

\rall{When dealing with a uniform grid, analytic expressions of the structural equations are available, and given in \cref{table::SE_equation_ZD_analytic} (obtained from the coefficients in \cite[appendix A]{CMMC25}). The benchmarks we provide in the numerical section do not use these analytic expressions, but rather the coefficients that are automatically derived from the kernel for any $R$.}

\subsubsection{\ZDS\ formulation}
\label{sec:full_algo_scalar_ZDS}
Assume that we knew the solution $(Z_n,D_n,S_n)$ at the time $t_n$. We seek the solution of the following nonlinear problem with $3R$ unknowns:
for all $m \in \{1, \dots, R\}$,
\begin{equation}\label{ZDS_scalar_problem}
    D_{n+m}=f(Z_{n+m})
    \text{, }
    \quad S_{n+m}=f'(Z_{n+m})D_{n+m}
    \text{ and \quad}
    \sum_{r=0}^R a^m_{r,0} Z_{n+r}+a^m_{r,1}\Delta t\, D_{n+r}+a^m_{r,2}\Delta t^2 S_{n+r}=0.
\end{equation}

Denoting by $\blockS_n$ the $R$-block of the second derivative approximations, we rewrite the structural equations under a matrix form
\begin{equation}
    \label{eq:ZDS_scheme_matrix_form}
    A_z\, \blockZ_n+\Delta t\,A_d\, \blockD_n+\Delta t^2\,A_s\, \blockS_n+ Z_n\, \ba_z+\Delta t\,D_n\, \ba_d+\Delta t^2\, S_n\, \ba_s=0,
\end{equation}
with
\begin{align*}
    A_z =(a^m_{r,0})_{m,r}, & \qquad \ba_z =\big (a^1_{0,0},a^2_{0,0},\cdots,a^R_{0,0}\big )^\trans, \\
    A_d =(a^m_{r,1})_{m,r}, & \qquad \ba_d =\big (a^1_{0,1},a^2_{0,1},\cdots,a^R_{0,1}\big )^\trans, \\
    A_s =(a^m_{r,2})_{m,r}, & \qquad \ba_s =\big (a^1_{0,2},a^2_{0,2},\cdots,a^R_{0,2}\big )^\trans.
\end{align*}
\begin{remark}
    \rb{For the sake of simplicity, we still use the notation $A_z$, $A_d$, $\ba_z$ and $\ba_d$ introduced in the \ZD\, formulation}. Of course, the matrix entries are different from the ones corresponding to the $\ZD$ method.
\end{remark}
Assuming the matrix $A_z\in\mathbb R^{R\times R}$ is non-singular, we rewrite the system as
\begin{equation}\label{eq::solve_SE}
    \blockZ_n+\Delta t\,B_d\,\blockD_n+\Delta t^2\,B_s\,\blockS_n+Z_n\, \bb_z+\Delta t\, D_n\, \bb_d+\Delta t^2\,S_n\,\bb_s=0,
\end{equation}
with $B_d=(A_z)^{-1}A_d$, $B_s=(A_z)^{-1}A_s$, $\bb_z=(A_z)^{-1}\ba_z$, $\bb_d=(A_z)^{-1}\ba_d$,$\bb_s=(A_z)^{-1}\ba_s$.

A fixed-point procedure is applied to solve the system \eqref{ZDS_scalar_problem}, by producing a sequence $(\blockZ_n^{(k)},\blockD_n^{(k)},\blockS_n^{(k)})$ that converges to the solution. We give the algorithm for the \ZDS\ case.
\begin{itemize}
    \item {\bf Initialization.} To build $\blockZ_n^{(0)}$, $\blockD_n^{(0)}$, $\blockS_n^{(0)}$, we set for $r\in\{1,\dots,R\}$
          $$
              Z^{(0)}_{n+r}=Z^{(0)}_{n+r-1}+\Delta t D^{(0)}_{n+r-1}+ \frac{\Delta t^2}{2} S^{(0)}_{n+r-1},\quad D^{(0)}_{n+r}=f(Z^{(0)}_{n+r}),\quad S^{(0)}_{n+r}=f'(Z^{(0)}_{n+r})D^{(0)}_{n+r},
          $$
          with $Z^{(0)}_{n}=Z_n$, $D^{(0)}_{n}=D_n$, $S^{(0)}_{n}=S_n$.
    \item {\bf Iteration.} We first compute a new approximation for the solution  using the structural equations:
          $$
              \blockZ_n^{(k+1)}=-\Big (Z_n\,b_z+\Delta t\, D_n\,b_d+\Delta t^2\,S_n\, b_s
              +\Delta t\,B_d\,\blockD_n^{(k)}+\Delta t^2\, B_s\, \blockS_n^{(k)}\Big ).
          $$
          Then, we update the first and second derivatives with the physical equations, by computing
          $$
              \forall r \in \{1, \dots, R\}, \qquad
              D^{(k+1)}_{n+r}=f(Z^{(k+1)}_{n+r})
              \text{\quad and \quad}
              S^{(k+1)}_{n+r}=f'(Z^{(k+1)}_{n+r})D^{(k+1)}_{n+r}.
          $$
    \item {\bf Stopping criterion.} We end the fixed-point when two successive solutions are close enough according to the tolerance parameter \texttt{tol}, that is $\Vert \blockZ_n^{(k+1)}-\blockZ_n^{(k)}\Vert \leq \texttt{tol}$.
\end{itemize}
\begin{remark}
    It has been shown \rbbis{in~\cite{CMMC25}} that $R=1$ provides a 4th-order unconditionally stable scheme while $R=2$ reaches sixth-order accuracy. Methods with $R=3$ and $R=4$ are also investigated in the numerical section.\qed
\end{remark}
\begin{remark}
    Given a grid, the structural equations' coefficients are computed independently of the problem, and their evaluations may be obtained in a pre-processing process or be stored with the grid points. \qed
\end{remark}

\rall{Once again, analytic expressions of the structural equations are available for uniform grids. We reproduce the results of~\cite{CMMC25} in \cref{table::SE_equation_ZDS_analytic}. As mentioned above, we only use the kernel method to compute the coefficients, making it possible to carry out the simulation for any $R$.}
%
%
\rb{
    \begin{remark}
        The authors showed in \cite{CMMC25} that the structural method is equivalent to a multi-derivative Runge-Kutta (MDRK) method as long as we use all the available physical equations. Note that the resulting MDKT is very specific since we guarantee that all the intermediate steps enjoy the same order of accuracy than the last step. Such a property is not true in general for the MDRK.
    \end{remark}
}
\section{Structural method for Hamiltonian problems: one-body systems}
\label{sec:Hamiltonian_scalar}

We reach the main novelty of the present work by adapting the structural method to Hamiltonian systems. The main difference is the introduction of a second variable $p$, that requires to handle both the approximations of $t\to x(t)\in\mathbb R$ and $t\to p(t)\in\mathbb R$, coupled through Hamilton's equations.
This section is dedicated to \rb{one-body} Hamiltonian equations, while \rb{multi-body} systems will be considered in the next section.
We first introduce the physical equations for the \rb{one-body} Hamiltonian problem in \cref{sec:PE_scalar_Hamiltonian}, and then detail the structural equations and the algorithm to solve the coupled physical and structural equations in \cref{sec:SE_scalar_Hamiltonian}.

%
%
\subsection{Physical equations}
\label{sec:PE_scalar_Hamiltonian}

Consider a smooth function $\mathcal H: \mathbb{R}^2 \to \mathbb{R}$ that takes as input \rb{the position $x$ and the momentum $p$ of some body subject to some exterior forces}.
We seek solutions $x=x(t)$ and $p=p(t)$ such that $\Ham{x(t)}{p(t)}$ is constant. Differentiating in time gives $\partial_x \Ham{x}{p}\dot x+\partial_p \Ham{x}{p}\dot p=0$, and we define the trajectories as the solution of the ODE
$$
    \dot x=\partial_p \Ham{x}{p}
    \text{\quad and \quad}
    \dot p=-\partial_x \Ham{x}{p},
$$
with the initial condition  $x(0)=x_0$, $p(0)=p_0$.
This ODE makes up the first physical equations $\PE[1]$.
We reformulate the problem within the $Z$, $D$ framework by denoting by $\Zx$ and $\Dx$ the approximations of $x$ and $\dot x$, and adopt similar notations for $\Zp$ and $\Dp$. The physical equations then read
\begin{align}
    \Dx & = \partial_p\Ham{\Zx}{\Zp}, \label{PE1x} \\
    \Dp & =-\partial_x\Ham{\Zx}{\Zp}. \label{PE1p}
\end{align}
For each time step, we have $2$ physical equations, with $4$ unknowns in total, and so we need $2$ structural equations to close the system.

In the \ZDS\ case, to provide the second physical system $\PE[2]$, we differentiate in time the first physical equations, and we obtain
$$
    \ddot x=\partial_{xp}\Ham{x}{p}\dot x+\partial_{pp}\Ham{x}{p}\dot p
    \text{\qquad and \qquad}
    \ddot p=-\partial_{xx}\Ham{x}{p}\dot x-\partial_{xp}\Ham{x}{p}\dot p.
$$
Reformulating the problem within the $Z$, $D$, $S$ framework yields the second physical equations $\PE[2]$:
\begin{align}
    \Sx & = \partial_{xp}\Ham{\Zx}{\Zp}\Dx+\partial_{pp}\Ham{\Zx}{\Zp}\Dp, \label{PE2x} \\
    \Sp & =-\partial_{xx}\Ham{\Zx}{\Zp}\Dx-\partial_{xp}\Ham{\Zx}{\Zp}\Dp. \label{PE2p}
\end{align}
Note that we have $4$ physical equations, with $6$ unknowns in total, and so $2$ structural equations are required to close the system for each time step.
%
%
\subsection{Structural equations and algorithms}
\label{sec:SE_scalar_Hamiltonian}

Equipped with the physical and structural equations in the \ZD\ and \ZDS\ formulations, coupling $x$ and $p$, we now give the algorithm to solve them in both cases.

\subsubsection{\ZD\ schemes}
Approximations of the function $x$ and the derivative $\dot x$ are connected via the structural equations \eqref{R_ZD_structural_equation}, and similarly for the function $p$ and its derivative.
Hence, the structural equations read
\begin{equation*}
    \forall m\in\{1,\dots,R\}, \qquad
    \begin{dcases}
        \sum_{r=0}^R a^m_{r,0} \Zx_{n+r}+a^m_{r,1}\Delta t\, \Dx_{n+r} = 0, \\
        \sum_{r=0}^R a^m_{r,0} \Zp_{n+r}+a^m_{r,1}\Delta t\, \Dp_{n+r} = 0.
    \end{dcases}
\end{equation*}
It is important to note that the $x$ and $p$ use the \textbf{same} structural equations (the same coefficients $a^m_{r,s}$)  and \textbf{only differ} through the physical ones. Let denote by
$$
    \blockZx_n=\big (\Zx_{n+1},\Zx_{n+2},\dots,\Zx_{n+R}\big )^\trans
    \text{\qquad and \qquad}
    \blockZp_n=\big (\Zp_{n+1},\Zp_{n+2},\dots,\Zp_{n+R}\big )^\trans,
$$
the respective $R$-block vectors for the functions $x$ and $p$. Similarly, the $R$-block vectors for the first derivatives are denoted by $\blockDx_n$ and $\blockDp_n$. The structural equations for the $\ZD$ scheme of size $R$ then read
\begin{align}
    0 & =\blockZx_n+B_d\,\Delta t\,\blockDx_n+\Zx_n\, \bb_z+\Dx_n\,\Delta t\, \bb_d,\label{ZD_SEx} \\
    0 & =\blockZp_n+B_d\,\Delta t\,\blockDp_n+\Zp_n\, \bb_z+\Dp_n\,\Delta t\, \bb_d.\label{ZD_SEp}
\end{align}
\begin{remark}
    These two linear systems are independent and should be treated in parallel.
\end{remark}
To solve the ODEs \ra{stemming} from the Hamiltonian, we proceed in a very similar way as in \cref{sec:full_algo_scalar_ZD}. We produce two sequences $\Big (\blockZx_n^{(k)},\blockDx_n^{(k)}\Big )$ and $\Big (\blockZp_n^{(k)},\blockDp_n^{(k)}\Big )$ that converge to the solution. The fixed-point algorithm is then given by the following iterative procedure.
\begin{itemize}
    \item {\bf Initialization.} To build $\blockZx_n^{(0)}$, $\blockDx_n^{(0)}$, $\blockZp_n^{(0)}$, $\blockDp_n^{(0)}$, we set for $r \in \{1,\dots,R\}$
          \begin{eqnarray*}
              &&\Zx^{(0)}_{n+r}=\Zx^{(0)}_{n+r-1}+\Delta t \Dx^{(0)}_{n+r-1},\\
              &&\Zp^{(0)}_{n+r}=\Zp^{(0)}_{n+r-1}+\Delta t \Dp^{(0)}_{n+r-1},\\
              &&\Dx^{(0)}_{n+r}= \partial_p\Ham{\Zx^{(0)}_{n+r}}{\Zp^{(0)}_{n+r}},\\
              &&\Dp^{(0)}_{n+r}=-\partial_x\Ham{\Zx^{(0)}_{n+r}}{\Zp^{(0)}_{n+r}},
          \end{eqnarray*}
          with $\Zx^{(0)}_{n}=\Zx_n$, $\Dx^{(0)}_{n}=\Dx_n$ and $\Zp^{(0)}_{n}=\Zp_n$, $\Dp^{(0)}_{n}=\Dp_n$.

    \item {\bf Iteration.} We first compute a new approximation for the solution  using the structural equations
          \begin{eqnarray*}
              \blockZx_n^{(k+1)}&=&-\Big (\Zx_n\,\bb_z+\Dx_n\,\Delta t\,\bb_d+B_d\,\Delta t\,\blockDx_n^{(k)}\Big ),\\
              \blockZp_n^{(k+1)}&=&-\Big (\Zp_n\,\bb_z+\Dp_n\,\Delta t\,\bb_d+B_d\,\Delta t\,\blockDp_n^{(k)}\Big ),
          \end{eqnarray*}
          and then update the first derivatives with the physical equations \ra{derived} from the Hamiltonian
          \begin{eqnarray*}
              &&\Dx^{(k+1)}_{n+r}=\partial_p\Ham{\Zx^{(k+1)}_{n+r}}{\Zp^{(k+1)}_{n+r}},\\
              &&\Dp^{(k+1)}_{n+r}=-\partial_x\Ham{\Zx^{(k+1)}_{n+r}}{\Zp^{(k+1)}_{n+r}}.
          \end{eqnarray*}
    \item {\bf Stopping criterion.} We end the fixed-point when two successive solutions are close enough according to the tolerance parameter \texttt{tol}, that is $\Vert \blockZx^{(k+1)}-\blockZx^{(k)}\Vert \leq \texttt{tol}$.
\end{itemize}
\begin{remark}
    As in the case of the scalar ODE problems, no local linearization is required, leading to elementary formulations. Nevertheless, we show in~\cite{CMMC25} that the fixed-point method converges under a CFL-like condition.$\qed$
\end{remark}

\subsubsection{\ZDS\ schemes}
Approximations of the function $x$ and its two derivatives $\dot x$ and $\ddot x$ are connected via the structural equations \eqref{R_ZDS_structural_equation}, as are $p$ and its derivatives.
Therefore, the structural equations read
\begin{equation*}
    \forall m\in\{1,\dots,R\}, \qquad
    \begin{dcases}
        \sum_{r=0}^R a^m_{r,0} \Zx_{n+r}+a^m_{r,1}\Delta t\, \Dx_{n+r}+a^m_{r,2}\Delta t^2\, \Sx_{n+r} = 0, \\
        \sum_{r=0}^R a^m_{r,0} \Zp_{n+r}+a^m_{r,1}\Delta t\, \Dp_{n+r}+a^m_{r,2}\Delta t^2\, \Sp_{n+r} = 0.
    \end{dcases}
\end{equation*}
It is important to note that $x$ and $p$ once again use exactly the same structural equations and only differ through the physical ones. We complete the notation with
$$
    \blockSx_n=\big (\Sx_{n+1},\Sx_{n+2},\cdots,\Sx_{n+R}\big )^\trans, \qquad
    \blockSp_n=\big (\Sp_{n+1},\Sp_{n+2},\cdots,\Sp_{n+R}\big )^\trans,
$$
the R-block vectors for the second derivatives. The two structural equations of size $R$ then read
\begin{align}
    0=\blockZx_n+B_d\,\Delta t\,\blockDx_n+B_s\,\Delta t^2\,\blockSx_n+\Zx_n\, \bb_z+\Delta t\,\Dx_n\, \bb_d+\Sx_n\,\Delta t^2\,\bb_s,\label{ZDS_SEx} \\
    0=\blockZp_n+B_d\,\Delta t\,\blockDp_n+B_s\,\Delta t^2\,\blockSp_n+\Zp_n\, \bb_z+\Delta t\,\Dp_n\, \bb_d+\Sp_n\,\Delta t^2\,\bb_s.\label{ZDS_SEp}
\end{align}

To solve the ODE \ra{stemming} from the Hamiltonian, we proceed similarly as in the previous section. We produce two sequences $\Big (\blockZx_n^{(k)},\blockDx_n^{(k)},\blockSx_n^{(k)}\Big )$ and $\Big (\blockZp_n^{(k)},\blockDp_n^{(k)},\blockSp_n^{(k)}\Big )$ that converge to the solution. The fixed-point algorithm is detailed below.
\begin{itemize}
    \item {\bf Initialization.} To build $\blockZx_n^{(0)}$, $\blockDx_n^{(0)}$, $\blockSx_n^{(0)}$, $\blockZp_n^{(0)}$, $\blockDp_n^{(0)}$, $\blockSp_n^{(0)}$, we set for $r \in \{1,\dots,R\}$
          \begin{align*}
              \Zx^{(0)}_{n+r} & =\Zx^{(0)}_{n+r-1}+\Delta t \Dx^{(0)}_{n+r-1}+ \frac{\Delta t^2}{2} \Sx^{(0)}_{n+r-1},                                                   \\
              \Zp^{(0)}_{n+r} & =\Zp^{(0)}_{n+r-1}+\Delta t \Dp^{(0)}_{n+r-1}+ \frac{\Delta t^2}{2} \Sp^{(0)}_{n+r-1},                                                   \\
              \Dx^{(0)}_{n+r} & = \partial_p\Ham{\Zx^{(0)}_{n+r}}{\Zp^{(0)}_{n+r}},                                                                                      \\
              \Dp^{(0)}_{n+r} & =-\partial_x\Ham{\Zx^{(0)}_{n+r}}{\Zp^{(0)}_{n+r}},                                                                                      \\
              \Sx^{(0)}_{n+r} & = \partial_{xp}\Ham{\Zx^{(0)}_{n+r}}{\Zp^{(0)}_{n+r}}\Dx^{(0)}_{n+r}+\partial_{pp}\Ham{\Zx^{(0)}_{n+r}}{\Zp^{(0)}_{n+r}}\Dp^{(0)}_{n+r}, \\
              \Sp^{(0)}_{n+r} & =-\partial_{xx}\Ham{\Zx^{(0)}_{n+r}}{\Zp^{(0)}_{n+r}}\Dx^{(0)}_{n+r}-\partial_{xp}\Ham{\Zx^{(0)}_{n+r}}{\Zp^{(0)}_{n+r}}\Dp^{(0)}_{n+r}.
          \end{align*}
          with $\Zx^{(0)}_{n}=\Zx_n$, $\Dx^{(0)}_{n}=\Dx_n$, $\Sx^{(0)}_{n}=\Sx_n$ and $\Zp^{(0)}_{n}=\Zp_n$, $\Dp^{(0)}_{n}=\Dp_n$, $\Sp^{(0)}_{n}=\Sp_n$.

    \item {\bf Iteration.} We first compute a new approximation for the solution  using the structural equations
          \begin{eqnarray*}
              \blockZx_n^{(k+1)}&=&-\Big (\Zx_n\,\bb_z+\Dx_n\,\Delta t\,\bb_d+\Sx_n\,\Delta t^2\, \bb_s+
              B_d\,\Delta t\,\blockDx_n^{(k)}+B_s\,\Delta t^2\, \blockSx_n^{(k)}\Big )\\
              \blockZp_n^{(k+1)}&=&-\Big (\Zp_n\,\bb_z+\Dp_n\,\Delta t\,\bb_d+\Sp_n\,\Delta t^2\, \bb_s+
              B_d\,\Delta t\,\blockDp_n^{(k)}+B_s\, \Delta t^2\,\blockSp_n^{(k)}\Big )
          \end{eqnarray*}
          and then update the first and second derivative with the physical equations \ra{derived} from the Hamiltonian
          \begin{eqnarray*}
              &&\Dx^{(k+1)}_{n+r}=\partial_p\Ham{\Zx^{(k+1)}_{n+r}}{\Zp^{(k+1)}_{n+r}},\\
              &&\Dp^{(k+1)}_{n+r}=-\partial_x\Ham{\Zx^{(k+1)}_{n+r}}{\Zp^{(k+1)}_{n+r}},\\
              &&\Sx^{(k+1)}_{n+r}=\partial_{xp}\Ham{\Zx^{(k+1)}_{n+r}}{\Zp^{(k+1)}_{n+r}}\Dx^{(k+1)}_{n+r}+\partial_{pp}\Ham{\Zx^{(k+1)}_{n+r}}{\Zp^{(k+1)}_{n+r}}\Dp^{(k+1)}_{n+r},\\
              &&\Sp^{(k+1)}_{n+r}=-\partial_{xx}\Ham{\Zx^{(k+1)}_{n+r}}{\Zp^{(k+1)}_{n+r}}\Dx^{(k+1)}_{n+r}-\partial_{xp}\Ham{\Zx^{(k+1)}_{n+r}}{\Zp^{(k+1)}_{n+r}}\Dp^{(k+1)}_{n+r}.
          \end{eqnarray*}
    \item {\bf Stopping criterion.} We end the fixed-point algorithm when two successive solutions are close enough, according to the tolerance parameter \texttt{tol}, that is $\Vert \blockZx^{(k+1)}-\blockZx^{(k)}\Vert \leq \texttt{tol}$.
\end{itemize}

\begin{remark}
    The coupling between the primal $x$ and dual variables $p$ takes place in the physical equations, hence the computation of $\blockZx_n^{(k+1)}$ and $\blockZp_n^{(k+1)}$ are independent and can be carried out in parallel. \qed
\end{remark}
\begin{remark}
    We have just considered a very basic fixed-point method for the sake of simplicity, but any accelerating techniques could be exploited. In practice, the predictor is good enough and very few iterations (3 to 12 in practice, depending on the tolerance and order of the method) are required to provide the accurate approximation. Of course, for stiff problems, a more sophisticated fixed-point procedure may be implemented.\qed
\end{remark}

{\rallparagraph
    \section{Conservation properties of the structural schemes}
    \label{sec:conservation_properties}

    Before moving on to \rb{multi-body systems}, we first check a few conservation properties of the \ZD\ and \ZDS\ schemes applied to \rb{one-body} Hamiltonian problems.
    Namely, we discuss their (the lack of) symplecticity in \cref{sec:non_symplecticity}, and their exact preservation of quadratic Hamiltonians in \cref{sec:exact_quadratic}.

    \subsection{Structural scheme and symplecticity}
    \label{sec:non_symplecticity}

    The symplecticity is the ability to exactly preserve geometric properties of the flow. In \cite{L88,S88}, a set of conditions for the Butcher tableau coefficients is provided, corresponding to $R(R+1)/2$ additional constraints on the coefficients given by
    $$
        \alpha_{R+1,m}\alpha_{m,s}+\alpha_{R+1,s}\alpha_{s,m}=\alpha_{R+1,m}\alpha_{R+1,s},\quad 1\leq m,s\leq R+1,
    $$
    taking into account that $\beta_m=\alpha_{R+1,m}$ in the context of the structural equations.

    In conclusion, we have $R(R+1)/2$ additional constraints on RK matrix $A$ and vector $b$ that have already been defined by the $R+2$ accuracy conditions (to be exact for polynomials up to degree $R+1$) together with the mandatory condition that the kernel is at least of dimension $R$ to provide $R$ linearly independent structural equations. Consequently, except one very particular case ($R=1$), the $R(R+1)/2$ conditions cannot be satisfied, hence the symplecticity is not achieved.

    Moreover, Hairer {\it et al.} \cite{HMS94} prove that there is no two-derivative symplectic RK scheme. Consequently, the $\ZDS$ scheme can not be symplectic.
}

{\rallparagraph
    \subsection{Exact conservation of quadratic Hamiltonians}
    \label{sec:exact_quadratic}

    In 2006, Ernst Hairer \cite{H06} shows that exact energy conservation and symplecticity for Hamiltonian systems cannot be obtained simultaneously. Moreover, symplectic integrators are not necessarily the main objective and researchers promote the energy or momentum conservation strategy for the long-time simulations \cite{BDZ23}. A first quantity that we expect to preserve at the numerical level is the Hamiltonian. We prove in this section that the structural schemes $\ZD$ or $\ZDS$ with $R$ blocks exactly preserve quadratic Hamiltonians every $R$ steps.


    \subsubsection{The quadratic Hamiltonian}
    Let us consider the quadratic Hamiltonian of the form
    \begin{equation}
        \label{eq:general_quadratic_Hamiltonian}
        \Ham{x}{p} = \frac 1 2 p^2 + \gamma x p + \frac 1 2 \omega^2 x^2,
    \end{equation}
    with $\gamma, \omega \in \mathbb{R}$. Then, the associated ODEs read
    \begin{equation*}
        \dot x = \gamma x + p,
        \text{\qquad and \qquad}
        \dot p = -\omega^2 x - \gamma p,
    \end{equation*}
    or in the compact form
    $$
        \blockDx=\gamma \blockZx+\blockZp, \text{\qquad and \qquad} \blockDp=-\omega^2 \blockZx-\gamma\blockZp.
    $$
    We also consider the ODEs satisfied by the second derivatives
    \begin{equation*}
        \ddot x = (\gamma^2 -\omega^2) x,
        \text{\qquad and \qquad}
        \ddot p = (\gamma^2 -\omega^2) p,
    \end{equation*}
    which we rewrite under the form
    $$
        \blockSx=(\gamma^2 -\omega^2)\blockZx,\text{\qquad and \qquad} \blockSp = (\gamma^2 -\omega^2)\blockZp,
    $$
    to design the $\ZDS$ scheme.

    For the sake of simplicity, we deal with the separable case, taking $\gamma = 0$, and  assume $\omega = 1$. Then, the equations reduce to
    \begin{equation}
        \label{eq:compact_quadratic_physical_equation}
        \blockDx=\blockZp,\quad \blockDp=-\blockZx,\quad \blockSx=-\blockZx,\quad  \blockSp =-\blockZp.
    \end{equation}
    \subsubsection{\texorpdfstring{Hamiltonian preservation for the $\ZD$ scheme}{Hamiltonian preservation for the ZD scheme}}
    We begin the study by considering the $\ZD$ scheme, and prove the following proposition.
    \begin{proposition}
        \label{prop:energy_conservation_ZD}
        Let $\mathcal{H}$ be the quadratic separable Hamiltonian given by \eqref{eq:general_quadratic_Hamiltonian} with $\gamma = 0$ and $\omega = 1$.
        We define the vector
        \begin{equation*}
            \bb_H=\Big( (I_R + \Delta t^2 B_d^2)^{-1} (\bb_z + \Delta t^2 B_d \bb_d) \Big)^{\otimes2}
            +
            \Delta t^2 \Big( (I_R + \Delta t^2 B_d^2)^{-1} (\bb_d - B_d \bb_z ) \Big)^{\otimes2}
        \end{equation*}
        where the expression $\bb^{\otimes2}$ is the Hadamard square (the square component-wise) of some vector $\bb$, $I_R$ is the identity matrix of size $R$, while  $B_d$, $\bb_z$, and $\bb_d$ are the coefficients of the structural equations \eqref{ZD_SEx}--\eqref{ZD_SEp}.
        Then, the $R$-block $\ZD$ scheme exactly preserves the Hamiltonian at the $r$\textsuperscript{th} substep, that is
        \begin{equation*}
            \Ham{\Zx_{n+r}}{\Zp_{n+r}} = \Ham{\Zx_n}{\Zp_n},
        \end{equation*}
        if and only if $\bb_{H,r}=1$.

    \end{proposition}
    \begin{proof}
        As a first step, we plug the physical equations \eqref{eq:compact_quadratic_physical_equation}
        into the structural equations \eqref{ZD_SEx}--\eqref{ZD_SEp},
        and we obtain
        \begin{align}
            \label{eq:ZD_SEx_proof}
            \blockZx_n + \Delta t \, B_d\,\blockZp_n & = - \Zx_n\, \bb_z - \Delta t \, \Zp_n\, \bb_d, \\
            \label{eq:ZD_SEp_proof}
            \blockZp_n - \Delta t \, B_d\,\blockZx_n & = - \Zp_n\, \bb_z + \Delta t \, \Zx_n\, \bb_d.
        \end{align}
        With
        \begin{equation*}
            \blockZx_n =
            \begin{pmatrix}
                \Zx_{n+1} \\ \Zx_{n+2} \\ \vdots \\ \Zx_{n+R}
            \end{pmatrix}
            \text{\quad and \quad}
            \blockZp_n =
            \begin{pmatrix}
                \Zp_{n+1} \\ \Zp_{n+2} \\ \vdots \\ \Zp_{n+R}
            \end{pmatrix},
        \end{equation*}
        the goal is to exhibit conditions on the coefficients of the structural equations such that
        \begin{equation}
            \label{eq:goal_energy_conservation}
            \frac 1 2 \Zp_{n+r}^2 + \frac 1 2 \Zx_{n+r}^2 = \frac 1 2 \Zp_n^2 + \frac 1 2 \Zx_n^2
        \end{equation}
        for some $1 \leq r \leq R$.
        To that end, we plug $\blockZp_n$ from \eqref{eq:ZD_SEp_proof} into \eqref{eq:ZD_SEx_proof} and $\blockZx_n$ from \eqref{eq:ZD_SEx_proof} into \eqref{eq:ZD_SEp_proof}, to obtain
        \begin{align*}
            (I_R + \Delta t^2 B_d^2) \, \blockZx_n
             & =
            - \Zx_n \, (\bb_z + \Delta t^2 B_d \bb_d)
            - \Delta t \, \Zp_n \, (\bb_d - B_d \bb_z), \\
            (I_R + \Delta t^2 B_d^2) \, \blockZp_n
             & =
            - \Zp_n \, (\bb_z + \Delta t^2 B_d \bb_d)
            + \Delta t \, \Zx_n \, (\bb_d - B_d \bb_z).
        \end{align*}
        Squaring the vectors $\blockZx_n$ and $\blockZp_n$ component-wise and summing them up, we get
        $$
            \blockZx_n^{\otimes2} + \blockZp_n^{\otimes2} = \big( \Zx_n^2 + \Zp_n^2 \big)\bb_H.
        $$
        Energy conservation \eqref{eq:goal_energy_conservation} is thus satisfied for step $r$ if $\bb_{H,r}=1$.
    \end{proof}

    Plugging the matrices and vectors of the \ZD\ schemes (present in \cref{table::SE_equation_ZD_analytic}) into \propref{prop:energy_conservation_ZD}, we obtain the following result, showing that $\bb_{H,r}=1$ for $R \in \{1, 2, 3, 4\}$.

    \begin{corollary}
        \label{cor:energy_conservation_ZD}
        The \ZD\ scheme with $R \in \{1, 2, 3, 4\}$ blocks satisfies the conditions of \propref{prop:energy_conservation_ZD} for $r = R$, {\it i.e.}, $\Ham{\Zx_{n+R}}{\Zp_{n+R}} = \Ham{\Zx_n}{\Zp_n}$.
    \end{corollary}
    \begin{proof}
        We use symbolic computation, providing a rational function depending on the powers of $\Delta t$, and we highlight the leading term.
        \begin{itemize}
            \item $R=1$: we get the vector $\bb_H=1$
            \item $R=2$: we get the vector
                  $\bb_{H} = \begin{pmatrix} 1-\Delta t^4/c+O(\Delta t^6) \\ 1 \end{pmatrix}$, with $c=18$.
            \item $R=3$: we get the vector
                  $\bb_{H} = \begin{pmatrix} 1+8\Delta t^6/c+O(\Delta t^8)\\1+\Delta t^6/c+O(\Delta t^8) \\ 1 \end{pmatrix}$, with $c=18$.
            \item $R=4$: we get the vector $\bb_{H} =
                      \begin{pmatrix}
                          1+135\Delta t^6/c+O(\Delta t^8) \\
                          1+80\Delta t^6/c+O(\Delta t^8)  \\
                          1+135\Delta t^6/c+O(\Delta t^8) \\
                          1
                      \end{pmatrix}$, with $c=3600$.
        \end{itemize}
        In each case, we have $H(\Zx_n,\Dx_n)=H(\Zx_{n+R},\Dx_{n+R})$, while the intermediate solutions slightly deviate from $H(\Zx_n,\Dx_n)$. Moreover, they do not increase or decrease with respect to the final time~$T$, and they satisfy
        $H(\Zx_{n+r},\Dx_{n+r})=H(\Zx_{n+r+R},\Dx_{n+r+R})$.
    \end{proof}
    \begin{remark}
        The general case of a quadratic Hamiltonian \eqref{eq:general_quadratic_Hamiltonian} with $\gamma, \omega \in \mathbb{R}$ is treated in \cref{sec:non_separable_quadratic_Hamiltonians}.
        Basically, the results of \corref{cor:energy_conservation_ZD} carry over to this more general, non-separable setting.
        Hence, while it is not symplectic,
        the $R$-block \ZD\ scheme exactly preserves quadratic Hamiltonians every $R$ steps,
        which is a desirable property for long-time simulations.
    \end{remark}

    \subsubsection{\texorpdfstring{Hamiltonian preservation for the $\ZDS$ scheme}{Hamiltonian preservation for the ZDS scheme}}
    We obtain a similar property for the $\ZDS$ scheme, which we state hereafter.
    \begin{proposition}
        \label{prop:energy_conservation_ZDS}
        Let $\mathcal{H}$ be the quadratic separable Hamiltonian given by \eqref{eq:general_quadratic_Hamiltonian} with $\gamma = 0$ and $\omega = 1$.
        We define the vector
        \begin{equation}\label{eq:ZDS_stability_H}
            \bb_H = \Big( (I_R + \Delta t^2 \bar B^2 B_d^2)^{-1}
            (\bar B \bb_d + \Delta t^2 \bar B B_d \bb_d) \Big)^{\otimes2}
            +
            \Delta t^2 \Big( (I_R + \Delta t^2 \bar B^2 B_d^2)^{-1}
            (\bar B \bb_d - \bar B B_d \bar B \bar \bb) \Big)^{\otimes2},
        \end{equation}
        with $\bar B = (I_R - \Delta t^2 B_s)^{-1}$ and $\bar \bb = \bb_z - \Delta t^2 \bb_s$,
        where the expression $\bb^{\otimes2}$ is the Hadamard square (component-wise square) of some vector $\bb$, $I_R$ is the identity matrix of size $R$, while $B_d$, $B_s$, $\bb_z$, $\bb_d$, and $\bb_s$ are the coefficients of the structural equations \eqref{ZDS_SEx}--\eqref{ZDS_SEp}.
        The $R$-block $\ZDS$ scheme exactly preserves the Hamiltonian at the $r$\textsuperscript{th} substep, that is
        \begin{equation*}
            \Ham{\Zx_{n+r}}{\Zp_{n+r}} = \Ham{\Zx_n}{\Zp_n},
        \end{equation*}
        if and only if $\bb_{H,r}=1$.
    \end{proposition}

    \begin{proof}
        The proof follows the same lines as the one of \propref{prop:energy_conservation_ZD},
        starting with plugging the physical equations    \eqref{eq:compact_quadratic_physical_equation}
        into the structural equations \eqref{ZDS_SEx}--\eqref{ZDS_SEp}, yielding
        \begin{align*}
            \blockZx_n + \Delta t \, B_d\,\blockZp_n - \Delta t^2 B_s\,\blockZx_n
             & =
            - \Zx_n\, \bb_z - \Delta t \, \Zp_n\, \bb_d + \Delta t^2 \bb_s\,Zx_n, \\
            \blockZp_n - \Delta t \, B_d\,\blockZx_n - \Delta t^2 B_s\,\blockZp_n
             & =
            - \Zp_n\, \bb_z + \Delta t \, \Zx_n\, \bb_d + \Delta t^2 \bb_s\,Zp_n.
        \end{align*}
        Rearranging the terms, we get
        \begin{align}
            \label{eq:ZDS_SEx_proof}
            \blockZx_n + \Delta t \, \bar B \, B_d\,\blockZp_n
             & =
            - \Zx_n \, \bar B \bar \bb - \Delta t \, \Zp_n\, \bar B \, \bb_d, \\
            \label{eq:ZDS_SEp_proof}
            \blockZp_n - \Delta t \, \bar B \, B_d\,\blockZx_n
             & =
            - \Zp_n \, \bar B \bar \bb + \Delta t \, \Zx_n\, \bar B \, \bb_d,
        \end{align}
        where we have set
        \begin{equation*}
            \bar B = (I_R - \Delta t^2 B_s)^{-1}
            \text{\quad and \quad}
            \bar \bb = \bb_z - \Delta t^2 \bb_s.
        \end{equation*}
        After substitution, we get
        $$
            [I_R+(\Delta t \bar B B_d)^2]\blockZx=-\Zx\Big \{ \bar B\bar \bb+\Delta t^2 \bar B B_d \bar B \bb_d \Big \}
            -\Zp\Big \{\Delta t \bar B \bb_d- \Delta t \bar B B_d \bar B \bar \bb \Big \},
        $$
        $$
            [I_R+(\Delta t \bar B B_d)^2]\blockZp=\Zx\Big \{ \Delta t \bar B \bb_d-\Delta t \bar B B_d \bar B \bar \bb \Big \}
            -\Zp\Big \{\bar B \bar \bb+\Delta t^2 \bar B B_d \bar B \bb_d \Big \}.
        $$
        Let define
        $$
            \bb_1=[I_R+(\Delta t \bar B B_d)^2]^{-1}\Big \{ \bar B\bar \bb+\Delta t^2 \bar B B_d \bar B \bb_d \Big \},\quad
            \bb_2=[I_R+(\Delta t \bar B B_d)^2]^{-1}\Big \{ \Delta t \bar B \bb_d-\Delta t \bar B B_d \bar B \bar \bb\Big \}.
        $$
        Then, we have $\blockZx^{\otimes 2}+\blockZp^{\otimes 2}=\bb_H(\Zx^2+\Zp^2)$ with $\bb_H=\bb_1^{\otimes 2}+\bb_2^{\otimes 2}$, which is nothing but the expected relation \eqref{eq:ZDS_stability_H}.
    \end{proof}
    We now analyze the energy conservation of the \ZDS\ scheme using the coefficients from \cref{table::SE_equation_ZDS_analytic} to compute the vector $\bb_{H}$ for $R\in\{1,2,3,4\}$.
    \begin{corollary}
        \label{cor:energy_conservation_ZDS}
        The \ZDS\ scheme with $R \in \{1, 2, 3, 4\}$ blocks satisfies the conditions of \propref{prop:energy_conservation_ZDS} for $r = R$, {\it i.e.}, $\Ham{\Zx_{n+R}}{\Zp_{n+R}} = \Ham{\Zx_n}{\Zp_n}$.
    \end{corollary}
    \begin{proof}
        We use a symbolic computation software that provides the entries of the vector $\bb_H$ as rational functions, given hereafter.
        \begin{itemize}
            \item $R=1$: we get the vector $\bb_H=1$
            \item $R=2$: we get the vector $\bb_{H} = \begin{pmatrix} 1-120\Delta t^6/c+O(\Delta t^8) \\ 1 \end{pmatrix}$,
                  with $\displaystyle c=129600$
            \item $R=3$: we get the vector $\bb_{H} =
                      \begin{pmatrix} 1+672\Delta t^{10}/c+O(\Delta t^{12})\\1+672\Delta t^{10}/c+O(\Delta t^{12}) \\ 1 \end{pmatrix}$,
                  with $\displaystyle c= 25401600$
            \item $R=4$: we get the vector $\bb_{H} =
                      \begin{pmatrix}
                          1-127575\Delta t^{12}/c_1+O(\Delta t^{14}) \\
                          1-302400\Delta t^{12}/c_2+O(\Delta t^{14}) \\
                          1-300875\Delta t^{12}/c_1+O(\Delta t^{14}) \\
                          1
                      \end{pmatrix}$,
                  with $\left \{ \begin{matrix}
                          c_1=22861440000 \\
                          c_2=51438240000
                      \end{matrix} \right .$
        \end{itemize}
        Once again, we have $H(\Zx_n,\Dx_n)=H(\Zx_{n+R},\Dx_{n+R})$ while the intermediate solutions slightly deviate from the $H(\Zx_n,\Dx_n)$. Furthermore, the non-dependence on the Hamiltonian on the final time~$T$ comes from the property
        $H(\Zx_{n+r},\Dx_{n+r})=H(\Zx_{n+r+R},\Dx_{n+r+R})$.
    \end{proof}
}

%
%
%
%
\section{Structural method for Hamiltonian problems: multi-body systems}
\label{sec:Hamiltonian_vector}

We now generalize the method to the situation where the Hamiltonian system involves $K$ vector-valued bodies. Such problems arise when modelling planetary systems, charged particles in electromagnetic fields, multiple coupled mass-damped-spring systems, among many other cases. The main issue of this extension lies in handling the technical difficulties arising when discretizing such complex systems. Overcoming these difficulties requires new notation and specific operators, which we introduce in \cref{sec:notation_vector_Ham}. The structural algorithm is then recast using this notation, in \cref{sec:algo_vector_Ham}.
%
%
\subsection{Notation}
\label{sec:notation_vector_Ham}
The Hamiltonian system involves $K$ distinct bodies characterized by their positions $\bx^k\in\mathbb R^I$ and momenta $\bp^k \in\mathbb R^I$. We gather the $\bx$ and $\bp$ in the ``body space'' $\Kbody=\mathbb R^{I\times K}$, and we introduce the matrix notation
$$
    X=\begin{bmatrix}
        x_1^1  & x_1^2 & \cdots & x_1^K  \\
        x_2^1  & x_2^2 & \cdots & x_2^K  \\
        \vdots &       &        & \vdots \\
        x_I^1  & x_I^2 &        & x_I^K  \\
    \end{bmatrix}
    \text{\qquad and \qquad}
    P=\begin{bmatrix}
        p_1^1  & p_1^2 & \cdots & p_1^K  \\
        p_2^1  & p_2^2 & \cdots & p_2^K  \\
        \vdots &       &        & \vdots \\
        p_I^1  & p_I^2 &        & p_I^K  \\
    \end{bmatrix}.
$$
Let $X,P \in \Kbody$ and $\Ham{X}{P}\in \mathbb R$ be the Hamiltonian. We adopt the following notation.

\begin{itemize}
    \item The gradients $\nabla_X \Ham{X}{P}\in \Kbody$ and $\nabla_P \Ham{X}{P}\in \Kbody$ are given by
          $$
              \nabla_X \Ham{X}{P}=\left [\partial_{x_i^k}\Ham{X}{P} \right ]_{k\in\{1,\dots,K\}, i\in\{1,\dots,I\}}.
          $$
          Note that the gradient does not output a vector but is reshaped to provide an element in $\Kbody$.
    \item The second derivatives $\nabla_{XX} \Ham{X}{P}$, $\nabla_{PP} \Ham{X}{P}$, $\nabla_{XP} \Ham{X}{P}=\nabla_X\left (\nabla_P \Ham{X}{P} \right )$, and $\nabla_{PX} \Ham{X}{P}=\nabla_P\left (\nabla_X \Ham{X}{P} \right )$ belong to $\Kbody^2$. Note that the order $\nabla_{XP}$ and $\nabla_{PX}$ are, in general, different.
          We introduce the following computation rules, using Einstein summation:
          \begin{eqnarray*}
              \nabla_{PX} \Ham{X}{P}\cdot \dot p&=&\left [\partial_{p^m_j}\partial_{x_i^k}\Ham{X}{P}\,\dot p^m_j \right ]_{k\in\{1,\dots,K\}, i\in\{1,\dots,I\}},\\
              \nabla_{XP} \Ham{X}{P}\cdot \dot x&=&\left [\partial_{x^m_j}\partial_{p_i^k}\Ham{X}{P}\,\dot x^m_j \right ]_{k\in\{1,\dots,K\}, i\in\{1,\dots,I\}}.
          \end{eqnarray*}
\end{itemize}
We also define the tensors that gather all the data we shall handle in the structural method.
\begin{itemize}
    \item We denote by $\ZX_n\approx \bx(t_n)\in \Kbody$ an approximation of the positions of the $K$ bodies at the time $t_n$. Similar notation is used for the first and second derivatives, $\DX_n\approx \dot \bx(t_n)\in \Kbody$ and $\SX_n\approx \ddot \bx(t_n)\in \Kbody$.
          We adopt the same notations $\ZP_n$, $\DP_n$ and $\SP_n$ for the momentum $\bp$. Moreover, for anybody $k$, the row vector $\ZX_n[k]$ represents its space coordinates in $\mathbb R^I$.
    \item Given a block of size $R$, we introduce the $\Kbody$-valued vector $\blockZX_n\in\Kbody^R$ that gathers the coordinates of all the $K$ bodies from time $t_{n+1}$ until time $t_{n+R}$ with the convention  $\blockZX_n[r]=\ZX_{n+r}\in \Kbody$ and $\blockZX_n[r][k]=\ZX_{n+r}[k]\in \mathbb R^I$.
    \item We introduce similar notation for $\blockZP_n$ and the derivatives, namely $\blockDX_n$, $\blockDP_n$, $\blockSX_n$, $\blockSP_n$, elements of $\Kbody^R$. The main difference compared to the previous section is that $\ZX_{n+r}$ now belongs to the body space $\Kbody$ instead of being real values.
\end{itemize}

A last point concerns the linear algebra with elements of $\Kbody$. To this end, we introduce the following notation
\begin{itemize}
    \item Let $\ZX\in\Kbody$ and $\alpha\in\mathbb R$. Then $\alpha \ZX\in\Kbody$ is the usual product of a matrix with a real number.
    \item Let $\ZX\in\Kbody$ and $a=(a_r)_r\in\mathbb R^R$, then $a\otimes \ZX$ is the vector in $\Kbody^R$ given by
          $$
              a\otimes \ZX=\left [\begin{array}{c}
                      a_1 \ZX \\
                      a_2 \ZX \\
                      \vdots  \\
                      a_R \ZX \\
                  \end{array} \right ].
          $$
    \item Let  $\blockZX\in\Kbody^R$ and $A=(a_{r,m})_{r,m}$ a square matrix of dimension $R$, then $A\times \blockZX$ is the vector in $\Kbody^R$ given by
          $$
              A\times \blockZX=\left [\begin{array}{c}
                      y_1    \\
                      y_2    \\
                      \vdots \\
                      y_r    \\
                  \end{array} \right ], \qquad
              \text{ with } y_r=\sum_{m=1}^R a_{r,m}\blockZX[m]\in\Kbody
              \text{\quad for all } r \in \{1,\dots,R\}.
          $$
\end{itemize}
The \ZD\ structural equations \eqref{ZD_SEx}--\eqref{ZD_SEp} adapted to vectors of $\Kbody^R$ read
\begin{eqnarray}
    0=\blockZX_n+B_d\times \blockDX_n+b_z\otimes\ZX_n+b_d\otimes\DX_n,\label{ZD_vSEx}\\
    0=\blockZP_n+B_d\times \blockDP_n+b_z\otimes\ZP_n+b_d\otimes\DP_n,\label{ZD_vSEp}
\end{eqnarray}
while the \ZDS\ structural equations \eqref{ZDS_SEx}--\eqref{ZDS_SEp} read
\begin{eqnarray}
    0=\blockZX_n+B_d\times \blockDX_n+B_s\times \blockSX_n+b_z\otimes\ZX_n+b_d\otimes\DX_n+b_s\otimes\SX_n,\label{ZDS_vSEx}\\
    0=\blockZP_n+B_d\times \blockDP_n+B_s\times \blockSP_n+b_z\otimes\ZP_n+b_d\otimes\DP_n+b_s\otimes\SP_n.\label{ZDS_vSEp}
\end{eqnarray}
\begin{remark}
    We highlight again that we use the same notation for matrices $B_d$ and vectors $b_z$, $b_d$ for the sake of simplicity, but they are different for the $\ZD$ or the $\ZDS$ methods. Additionally, the computation of the products is largely parallelizable, which may be taken advantage of in an HPC context. $\qed$
\end{remark}

\subsection{Algorithm}
\label{sec:algo_vector_Ham}

We rephrase the fixed-point method with the notation introduced in \cref{sec:notation_vector_Ham}.
Once again, we need an initialization and an iteration routine.
We just present the $\ZDS$ version.
Indeed, the simpler $\ZD$ version only consists in extracting from the $\ZDS$ version the two first physical equations \eqref{body_PE1_x}--\eqref{body_PE1_p} and structural equations \eqref{ZD_SEx}--\eqref{ZD_SEp}.

Given the initial state $(\ZX_n,\DX_n,\SX_n)$ and $(\ZP_n,\DP_n,\SP_n)$, we define two sequences $(\blockZX_n^{(k)},\blockDX_n^{(k)},\blockSX_n^{(k)})$ and $(\blockZP_n^{(k)},\blockDP_n^{(k)},\blockSP_n^{(k)})$ of elements of $\Kbody^R$ that converge to the $R$-block solution. We recall the convention
$\blockZX_n^{(k)}[r]=\ZX_{n+r}\in \Kbody$ for all $r \in \{1,\dots, R\}$.
\begin{itemize}
    \item {\bf Initialization.} To initialize the $R$-size blocks $\blockZX_n^{(0)}$, $\blockDX_n^{(0)}$, $\blockSX_n^{(0)}$, $\blockZP_n^{(0)}$, $\blockDP_n^{(0)}$, $\blockSP_n^{(0)}$ for all $r \in \{1,\dots, R\}$, we use the second order Taylor expansion in space $\Kbody$:
          \begin{eqnarray*}
              &&\ZX^{(0)}_{n+r}=\ZX^{(0)}_{n+r-1}+\Delta t\, \DX^{(0)}_{n+r-1}+ \frac{\Delta t^2}{2}\, \SX^{(0)}_{n+r-1},\\
              &&\ZP^{(0)}_{n+r}=\ZP^{(0)}_{n+r-1}+\Delta t\, \DP^{(0)}_{n+r-1}+ \frac{\Delta t^2}{2}\, \SP^{(0)}_{n+r-1},
          \end{eqnarray*}
          where we compute the first and second derivatives through the physical equations
          \begin{eqnarray*}
              &&\DX^{(0)}_{n+r}= \nabla_P\Ham{\ZX^{(0)}_{n+r}}{\ZP^{(0)}_{n+r}},\\
              && \DP^{(0)}_{n+r}=-\nabla_X\Ham{\ZX^{(0)}_{n+r}}{\ZP^{(0)}_{n+r}},\\
              &&\SX^{(0)}_{n+r}= \nabla_{XP}\Ham{\ZX^{(0)}_{n+r}}{\ZP^{(0)}_{n+r}}\DX^{(0)}_{n+r}+
              \nabla_{PP}\Ham{\ZX^{(0)}_{n+r}}{\ZP^{(0)}_{n+r}}\Dp^{(0)}_{n+r},\\
              &&\SP^{(0)}_{n+r}=-\nabla_{XX}\Ham{\ZX^{(0)}_{n+r}}{\ZP^{(0)}_{n+r}}\DX^{(0)}_{n+r}
              -\nabla_{PX}\Ham{\ZX^{(0)}_{n+r}}{\ZP^{(0)}_{n+r}}\DP^{(0)}_{n+r}
          \end{eqnarray*}
          with $\ZX^{(0)}_{n}=\ZX_n$, $\DX^{(0)}_{n}=\DX_n$, $\SX^{(0)}_{n}=\SX_n$ and
          $\ZP^{(0)}_{n}=\ZP_n$, $\DP^{(0)}_{n}=\DP_n$, $\SP^{(0)}_{n}=\SP_n$.

    \item {\bf Iteration.} We first compute new $R$-size blocks $\blockZX_n^{(k+1)}$ and $\blockZP_n^{(k+1)}$, approximating the positions and momenta, using the system of structural equations:
          \begin{eqnarray*}
              \blockZX_n^{(k+1)}&=&-\Big (b_z\otimes\ZX_n+b_d\otimes\DX_n+b_s\otimes\SX_n+B_d\times\blockDX_n^{(k)}+B_s\times\blockSX_n^{(k)}\Big ),\\
              \blockZP_n^{(k+1)}&=&-\Big (b_z\otimes\ZP_n+b_d\otimes\DP_n+b_s\otimes\SP_n+B_d\times\blockDP_n^{(k)}+B_s\times\blockSP_n^{(k)}\Big ).
          \end{eqnarray*}
          Then, for each $r \in \{1,\dots, R\}$, we update the first and second derivatives at each time step $t_{n+r}$ by setting
          \begin{eqnarray*}
              &&\DX^{(k+1)}_{n+r}= \nabla_P \Ham{\ZX^{(k+1)}_{n+r}}{\ZP^{(k+1)}_{n+r}}, \\
              &&\DP^{(k+1)}_{n+r}=-\nabla_X \Ham{\ZX^{(k+1)}_{n+r}}{\ZP^{(k+1)}_{n+r}},\\
              &&\SX^{(k+1)}_{n+r}= \nabla_{XP} \Ham{\ZX^{(k+1)}_{n+r}}{\ZP^{(k+1)}_{n+r}}\DX^{(k+1)}_{n+r}+
              \nabla_{PP} \Ham{\ZX^{(k+1)}_{n+r}}{\ZP^{(k+1)}_{n+r}}\DP^{(k+1)}_{n+r},\\
              &&\SP^{(k+1)}_{n+r}=-\nabla_{XX} \Ham{\ZX^{(k+1)}_{n+r}}{\ZP^{(k+1)}_{n+r}}\DX^{(k+1)}_{n+r}
              -\nabla_{PX} \Ham{\ZX^{(k+1)}_{n+r}}{\ZP^{(k+1)}_{n+r}}\DP^{(k+1)}_{n+r}.
          \end{eqnarray*}
    \item {\bf Stopping criterion.} We end the fixed-point scheme once two successive solutions are close enough, according to the tolerance parameter \texttt{tol}, that is $\Vert \ZX_{n+R}^{(k+1)}-\ZX_{n+R}^{(k)}\Vert \leq \texttt{tol}$.
\end{itemize}
\begin{remark}
    Once again, we highlight that the procedure is highly parallelizable since the update of the derivative for each time steps are independent.
\end{remark}
\begin{remark}
    \rallparagraph
    Since the same structural equation is applied to each component of the system of ODEs, the results of \cref{sec:conservation_properties} are directly applicable.
\end{remark}
%
%
%
%

\section{Benchmarks}
\label{sec:benchmarks}

We propose and analyze a series of benchmarks to assess the properties of the scheme.
Namely, we check the accuracy when an analytical solution is available, as well as the numerical preservation of several invariants, like the Hamiltonian, for long-time simulations.
In some cases, we also check the preservation of other invariant quantities of interest.
The method has been implemented in {\tt C++} in IEEE 754 {\tt quadruple} and {\tt octuple} precision using the library {\tt qd} developed by  \cite{Bailey2000}. The new structural method will be compared to standard symplectic schemes already implemented in {\tt julia 1.11}, see~\cite{bezanson2017julia}.
These schemes are provided by the \texttt{DifferentialEquations.jl} library, see \cite{RacNie2017}, and we use the second-order accurate \texttt{McAte2} scheme from \cite{McLAte1992}, the fourth-order accurate \texttt{CalvoSanz4} scheme from \cite{SanCal1993}, and the sixth- and eighth-order accurate \texttt{KahanLi6} and \texttt{KahanLi8} schemes from \cite{KahLi1997}.

The benchmarks and their settings are summarized in \cref{tab:test_cases}.
Moreover, \cref{sec:complexity} contains a study of the complexity of the \ZD\ and \ZDS\ methods, compared to the classical ones, in separable and non-separable settings, as well as a study of the computational cost of the fixed-point method.
\rb{Finally, \cref{sec:computation_time} provides a comparison of the computation time of the \ZD\ and \ZDS\ methods with respect to the reference symplectic schemes.}

\begin{table}[!ht]\centering
    \caption{\rempolice Summary of the benchmarks.}
    \label{tab:test_cases}
    \begin{tabular}{lcccc}
        \toprule
        \textbf{Benchmark}                     & \makecell{\textbf{Has an exact} \\ \textbf{solution}} & \makecell{\textbf{Additional} \\ \textbf{Invariants}} & \makecell{\textbf{Separable} \\ \textbf{Hamiltonian}} & \textbf{Section}        \\
        \midrule
        Two springs, two masses system         & Yes                                                   & No                                                    & Yes                                                   & \cref{sec:2S2M}         \\
        One-dimensional pendulum problem       & Yes                                                   & No                                                    & Yes                                                   & \cref{sec:pendulum}     \\
        Two-dimensional Kepler problem         & No                                                    & Yes                                                   & Yes                                                   & \cref{sec:2D_Kepler}    \\
        Outer Solar system                     & No                                                    & Yes                                                   & Yes                                                   & \cref{sec:solar_system} \\
        Particle in a 3D electromagnetic field & No                                                    & No                                                    & No                                                    & \cref{sec:3D_EM_field}  \\
        \bottomrule
    \end{tabular}
\end{table}

In each benchmark, the time interval $[0,T]$ is uniformly divided into $N$ steps of constant size $\Delta t$, and set $t_n=n\Delta t$, for all $n \in \{0,\dots,N\}$ to represent the subdivision of this interval. Denoting by $\bar q(t_n)$ the exact quantity at time $t_n$ (position, momentum, or invariant) and $q_n$ its approximation, we evaluate the error by computing $\eq_n(\rall{\Delta t})=|q_n-\bar x(t_n)|$, and we define the maximum error, together with the order, by
$$
    \eq(\Delta t)=\max_{n}\eq_n(\Delta t),  \qquad
    \ordq(\Delta t_1,\Delta t_2)=\frac{\log\Big (\eq(\Delta t_1)\big /\eq(\Delta t_2)\Big )}{\log(\Delta t_1\big /\Delta t_2)}.
$$
For instance, $\ex$ will represent the position error and $\eH$ the deviation of the Hamiltonian. All simulations have been carried out with {\tt quadruple} precision since the {\tt double} precision is not enough to correctly compute the tiny errors obtained when using very high order accurate methods (methods of eighth or tenth order, for instance).

\ra{We mention a general property of the structural schemes that the variables connected by Physical Equations enjoy the same order of accuracy (see the section on consistency errors in~\cite{CMMC25}). Consequently, to avoid an abundance of tables, we only present the errors and convergence order for the variable $\Zx$.}

%
%
%
%
\subsection{Two springs, two masses system}
\label{sec:2S2M}

We consider a system composed of two bodies linked by two springs, following \cref{Two_spring_fig}. The system is described by the Hamiltonian functional, whose inputs are two vectors $x=(x_1,x_2)$ and $p=(p_1,p_2)$ in $\mathbb R^2$, and whose expression is
\begin{equation}
    \label{eq:hamiltonian_2M2S}
    \mathcal{H}(x, p) =  \frac{p_1^2}{2 m_1} + \frac{p_2^2}{2 m_2} +  \frac 1 2 k_1 x_1^2 + \frac 1 2 k_2 (x_2 - x_1)^2.
\end{equation}
\begin{figure}[!ht]
    \centering
    \begin{tikzpicture}
        \fill [pattern = north east lines] (-1,1) rectangle (-1.2,0);
        \draw[thick] (-1,1) -- (-1,0);
        \node[circle, draw, inner sep=2.5mm] (a) at (2,0.5) {$m_1$};
        \draw[decoration={aspect=0.3, segment length=1.5mm, amplitude=3mm,coil},decorate] (-1,0.5) -- (a) node[midway, below, yshift=-3mm] {$k_1$};
        \node[circle, draw, inner sep=2.5mm] (b) at (7,0.5) {$m_2$};
        \draw[decoration={aspect=0.3, segment length=2mm, amplitude=3mm,coil},decorate] (a) -- (b) node[midway, below, yshift=-3mm] {$k_2$};
    \end{tikzpicture}
    \caption{\rempolice Schematic representation of the two masses-spring system.}
    \label{Two_spring_fig}
\end{figure}

The analytical solution $(x,p)$ is given by
\begin{equation*}
    \begin{aligned}
        x_1(t) & =A \cos{\left( {{\omega }_1} t+{{\alpha }_1}\right) }+B \cos{\left( {{\omega }_2} t+{{\alpha }_2}\right) },                                                                                                                                                         \\
        x_2(t) & =\frac{A \left( {k_1}+{k_2}-{m_1} {{\omega }_1^{2}}\right) }{{k_2}} \cos{\left( {{\omega }_1} t+{{\alpha }_1}\right) }+\frac{B {k_2}}{{k_2}-{m_2} {{\omega }_2^{2}}} \cos{\left( {{\omega }_2} t+{{\alpha }_2}\right) },                                            \\
        p_1(t) & =-{m_1} \left( A {{\omega }_1} \sin{\left( t {{\omega }_1}+{{\alpha }_1}\right) } + B {{\omega }_2} \sin{\left( t {{\omega }_2}+{{\alpha }_2}\right) }\right),                                                                                                      \\
        p_2(t) & =-m_2 \left(\frac{A{\omega }_1 \left( {k_1}+{k_2}-{m_1} {{\omega }_1^{2}}\right) }{{k_2}} \sin{\left( {{\omega }_1} t+{{\alpha }_1}\right) }+\frac{B{{\omega }_2} {k_2}}{{k_2}-{m_2} {{\omega }_2^{2}}} \sin{\left( {{\omega }_2} t+{{\alpha }_2}\right) } \right),
    \end{aligned}
\end{equation*}
where $\omega_1$ and $\omega_2$ are positive solutions to
\begin{equation*}
    m_1 m_2 \omega^4 - (m_1 k_2 + m_2 k_1 + m_2 k_2) \omega^2 + k_1 k_2 = 0.
\end{equation*}

\subsubsection{Physical equations}

Hamilton's equation $\dot x=\nabla_p \mathcal H(x,p)$ and $\dot p=-\nabla_x \mathcal H(x,p)$, applied to the Hamiltonian \eqref{eq:hamiltonian_2M2S}, yields the following dynamical system:
$$
    \dot x=\left ( \begin{array}{c}p_1/m_1\\p_2/m_2 \end{array} \right ),
    \qquad \dot p=-\left ( \begin{array}{l}k_1 x_1+k_2(x_1-x_2)\\ k_2(x_2-x_1) \end{array} \right ).
$$
Differentiating with respect to time leads to the system
$$
    \ddot x=\left ( \begin{array}{c}\dot p_1/m_1\\\dot p_2/m_2 \end{array} \right ), \qquad
    \ddot p=-\left ( \begin{array}{l}k_1 \dot x_1+k_2(\dot x_1-\dot x_2)\\ k_2(\dot x_2-\dot x_1) \end{array} \right ).
$$
We deduce the first group of structural equation $\PE[1]$, given by
$$
    \Dx=\left [\begin{array}{cc} 1/m_1& 0\\0 &1/m_2 \end{array}\right ] \Zp,\qquad
    \Dp=-\left [\begin{array}{cr} k_1+k_2& -k_2\\-k_2 &k_2 \end{array}\right ] \Zx,
$$
while the second group of physical equation $\PE[2]$ reads
$$
    \Sx=\left [\begin{array}{cc} 1/m_1& 0\\0 &1/m_2 \end{array}\right ] \Dp,\qquad
    \Sp=-\left [\begin{array}{cr} k_1+k_2& -k_2\\-k_2 &k_2 \end{array}\right ] \Dx.
$$

\subsubsection{Numerical tests}
The analytical solution is computed for the physical parameters $k_1= 1$, $k_2= 5$, $m_1= 2$, $m_2= 1$ that provides the two frequencies
\begin{equation*}
    \omega_1= \sqrt{\frac{8-3\sqrt6}{2}}\approx 0.57
    \quad \text{and} \quad
    \omega_2= \sqrt{\frac{8+3\sqrt6}{2}}\approx 2.77.
\end{equation*}
Moreover, the initial condition is made of the analytical solution taken at time $t=0$, with $A=1$, $B=2$, $\alpha_1=\pi/2$ and $\alpha_2=-\pi/4$. Numerical simulations have been carried out with the $\ZD$, $\ZDS$ and classical schemes. The error on the position and the conservation of the Hamiltonian are the two main indicators we use to assess the quality of the simulation.

\begin{figure}[!ht]
    \includegraphics{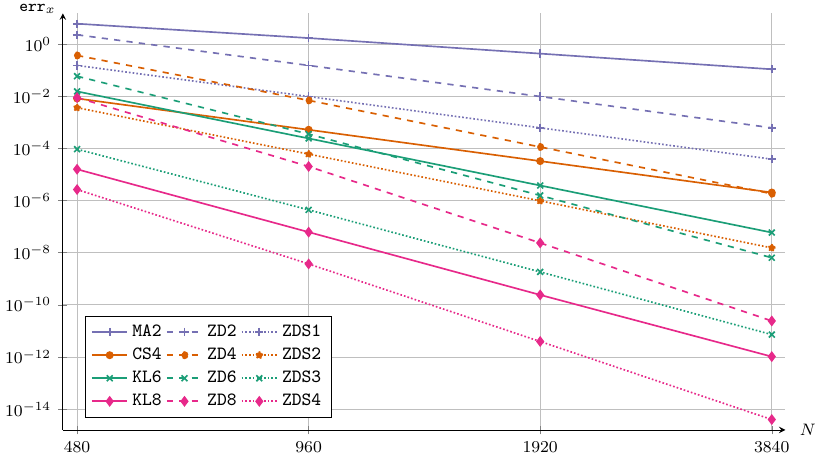}
    \caption{\rempolice Two springs, two masses test case from \cref{sec:2S2M}: error on the position at $T=100$s.}
    \label{fig:two_spring-mass_accuracy}
\end{figure}

\paragraph{Accuracy and order of convergence.} The final time is set to $T=100$s, and we take $N=480$, $960$, $1920$ and $3840$. We recall that the smallest period is \rb{$T_2 \approx 2.27$}, corresponding to $\omega_2$. Thus, the coarse grid $N=480$ corresponds to about $11$ points to cover a full period.
We display in \cref{fig:two_spring-mass_accuracy} the maximal error and order of convergence for the $\ZD$, $\ZDS$ and classical schemes, respectively. The optimal order is clearly achieved in each situation. \rbbis{For a given order of accuracy, the classical schemes generally yield a comparable or lower error than the $\ZDS$ scheme, with the $\ZD$ scheme presenting an error larger by about three orders of magnitude.}

\paragraph{Preservation of the Hamiltonian.}

To address the Hamiltonian preservation issue, we take $T=10\,000$ and assess the position error together with the Hamiltonian deviation, checking the evolution of the maximum error over time. A first note concerns the accuracy of the position: the errors increase with time, but the $\ZD$ and $\ZDS$ schemes provide lower deviations than the classical schemes. Indeed, for the same order, we observe a gain of two orders of magnitude when using the $\ZDS$ scheme rather than the classical schemes, while the difference between the $\ZD$ scheme and the standard schemes was reduced to a factor 10. In conclusion, the $\ZD$ and $\ZDS$ schemes offer a better accuracy when dealing with long-time simulations (i.e., numerous time steps).

\rall{Concerning the Hamiltonian conservation, the $\ZD$ and $\ZDS$ schemes reach machine precision (as expected from \cref{sec:conservation_properties}), but the standard schemes do not reach such a target.} Moreover, the deviation from the initial Hamiltonian depends on the grid size for the classical schemes, while errors with the $\ZD$ and $\ZDS$ schemes are independent of~$\Delta t$.

\begin{table}[!ht]\notapolice\centering
    \caption{\rempolice Two springs, two masses test case from \cref{sec:2S2M}: error on the position and the Hamiltonian at $T=10\,000$s. For $\ex$ and $\eH$, both rows correspond to $N=T\times 24$ and $N=T\times 96$ respectively.}
    \begin{subtable}{0.49\textwidth}\centering
        \caption{\rempolice Errors obtained with the \ZD\ scheme.}
        \label{tab:two_spring-mass_Ham_ZD}

        \begin{tabular}{@{}l@{}  r@{}l@{}  r@{}l@{}  r@{}l@{}  r@{}l@{}}
            \toprule
                             & \phantom{aaa} & {R=2}
                             & \phantom{aaa} & {R=4}
                             & \phantom{aaa} & {R=6}
                             & \phantom{aaa} & {R=8}                                              \\
            \midrule
            \ex ($\times24$) &               & 4.07e-01 &  & 3.08e-03 &  & 2.74e-05 &  & 2.60e-07 \\
            \ex ($\times96$) &               & 1.60e-03 &  & 7.59e-07 &  & 4.26e-10 &  & 2.62e-13 \\\midrule
            \eH ($\times24$) &               & 3.70e-26 &  & 9.13e-27 &  & 8.00e-27 &  & 2.57e-28 \\
            \eH ($\times96$) &               & 2.54e-25 &  & 6.04e-26 &  & 2.19e-26 &  & 1.99e-26 \\
            \bottomrule
        \end{tabular}

    \end{subtable}
    \begin{subtable}{0.49\textwidth}\centering
        \caption{\rempolice Errors obtained with the \ZDS\ scheme.}
        \label{tab:two_spring-mass_Ham_ZDS}

        \begin{tabular}{@{}l@{}  r@{}l@{}  r@{}l@{}  r@{}l@{}  r@{}l@{}}
            \toprule
                             & \phantom{aaa} & {R=1}
                             & \phantom{aaa} & {R=2}
                             & \phantom{aaa} & {R=3}
                             & \phantom{aaa} & {R=4}                                              \\
            \midrule
            \ex ($\times24$) &               & 2.55e-02 &  & 2.58e-05 &  & 3.10e-08 &  & 4.35e-11 \\
            \ex ($\times96$) &               & 9.97e-05 &  & 6.33e-09 &  & 4.76e-13 &  & 4.22e-17 \\\midrule
            \eH ($\times24$) &               & 1.01e-25 &  & 4.26e-26 &  & 7.99e-26 &  & 3.28e-26 \\
            \eH ($\times96$) &               & 6.61e-25 &  & 2.17e-25 &  & 4.49e-25 &  & 1.84e-26 \\
            \bottomrule
        \end{tabular}

    \end{subtable}

    \medskip

    \begin{subtable}{\textwidth}\centering
        \caption{\rempolice Errors obtained with the classical schemes.}
        \label{tab:two_spring-mass_Ham_classical}

        \begin{tabular}{@{}l@{}  r@{}l@{}  r@{}l@{}  r@{}l@{}  r@{}l@{}}
            \toprule
                             & \phantom{aaa} & {MA2}
                             & \phantom{aaa} & {CS4}
                             & \phantom{aaa} & {KL6}
                             & \phantom{aaa} & {KL8}                                              \\
            \midrule
            \ex ($\times24$) &               & 7.51e+00 &  & 1.36e-03 &  & 9.98e-05 &  & 4.06e-09 \\
            \ex ($\times96$) &               & 1.73e+00 &  & 5.29e-06 &  & 2.44e-08 &  & 1.06e-11 \\ \midrule
            \eH ($\times24$) &               & 6.57e-05 &  & 6.06e-08 &  & 5.56e-10 &  & 1.98e-14 \\
            \eH ($\times96$) &               & 1.02e-06 &  & 2.33e-10 &  & 1.35e-13 &  & 3.03e-19 \\
            \bottomrule
        \end{tabular}

    \end{subtable}
\end{table}
%
%
%
%
\subsection{One-dimensional pendulum problem}
\label{sec:pendulum}

We turn to a nonlinear \rb{one-body} problem, and consider the one-dimensional pendulum system governed by the Hamiltonian
\begin{equation}
    \mathcal{H}(x, p) = \frac {p^2} {2 m \ell^2} + m g \ell  (1 - \cos(x))
\end{equation}
where $m$ is the mass, $\ell$ the length of the pendulum and $g$ the gravity.

\subsubsection{Physical equations}
The differential system \ra{derived} from Hamilton's equations reads
\begin{equation}
    \label{eq:hamilton_equations_pendulum}
    \dot x=\frac{p}{m},\qquad \dot p=-mg\ell \sin(x)
\end{equation}
The solution is a periodic solution \cite{BPM2007} with a period given by the elliptic integral
$$
    T_p=4\sqrt{\frac{\ell}{mg}}\int_0^\frac{\pi}{2}\frac{du}{\sqrt{1-\omega^2\sin^2(u)}},
$$
with $\omega=\sin(\pi/8)$.
Then, differentiating the first physical equation \eqref{eq:hamilton_equations_pendulum} with respect to time provides relations involving the second-order derivative, namely
$$
    \ddot x=\frac{\dot p}{m},\qquad \ddot p=-mg\ell \cos(x)\dot x.
$$
We then deduce the first set of physical equations $\PE[1]$, given by
$$
    \Dx=\Zp/m,\qquad \Dp=-mg\ell \sin(\Zx),
$$
whereas the second set of physical equations $\PE[2]$ reads
$$
    \Sx=\Dp/m,\qquad \Sp=-mg\ell \cos(\Zx) \Dx.
$$

\subsubsection{Numerical tests}

The numerical applications have been carried out with $m = 1$, $g = 1$, $\ell = 1$, and the initial conditions $x(0) = \pi / 4$ and $p(0) = 0$. A numerical approximation of the period is $T_p=6.53$. In \cite{BPM2007}, an analytical solution is derived, whose expression, not reported here, involves a Jacobi elliptic function. The exact position and momentum at final time $T=100$ are given by
$$
    x(t=100) = -0.2633498226088722
    \text{\qquad and \qquad}
    p(t=100) = -0.7189111241830892
$$
in {\tt double} precision (since the Jacobi functions are not implemented with arbitrary precision). Hence, schemes with very high accuracy will reach machine error with a comparatively low number of points.

\begin{figure}[!ht]
    \includegraphics{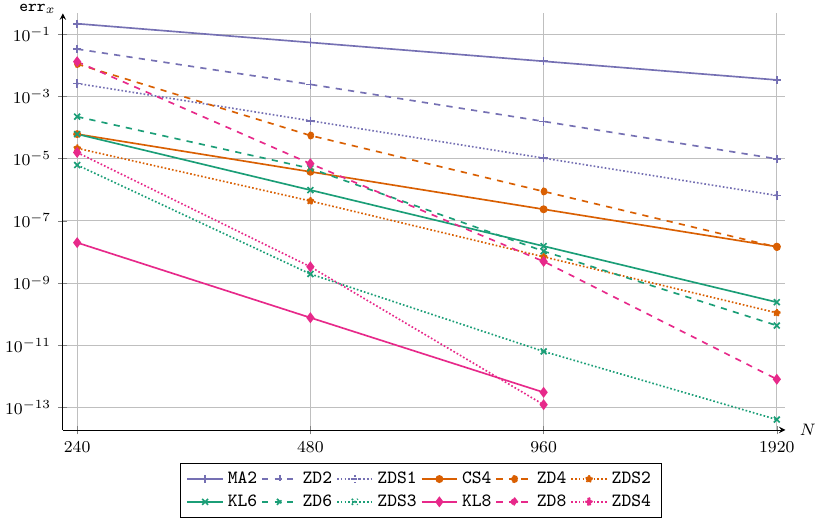}
    \caption{\rempolice Pendulum system from \cref{sec:pendulum}:  position error at $T=100$s.}
    \label{fig:pendulum_Xaccuracy}
\end{figure}

\begin{figure}[!ht]
    \includegraphics{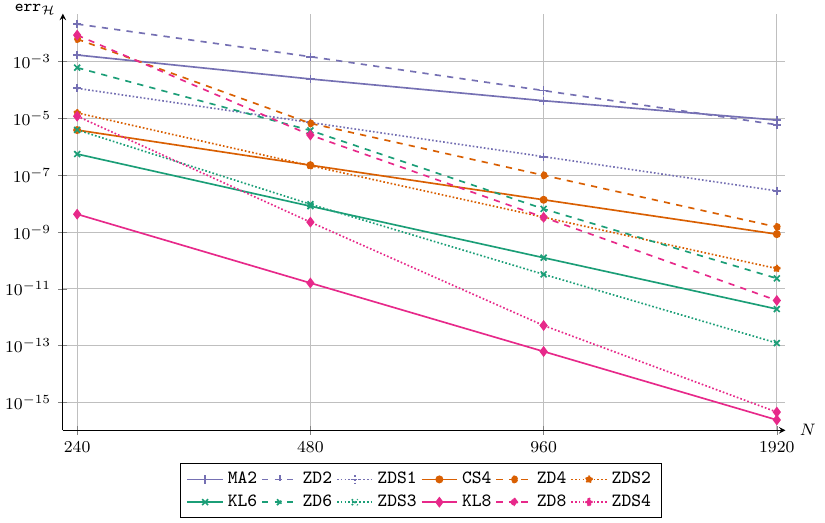}
    \caption{\rempolice Pendulum system  from \cref{sec:pendulum}:  Hamiltonian error  at $T=100s$.}
    \label{fig:pendulum_Haccuracy}
\end{figure}

Position errors, together with the deviation of the Hamiltonian are reported in \cref{fig:pendulum_Xaccuracy,fig:pendulum_Haccuracy} for $T=100$s. The deviation of $\mathcal{H}$ strongly depends on the order of the method, as well as on the time step. The $\ZD$ scheme produces larger errors, about two orders of magnitudes larger, than the $\ZDS$ scheme for the same rate of convergence.
\rbbis{The \texttt{CS4} outperforms the $\ZDS$ scheme with $R=1$, yielding a better accuracy for the same number of points.
    Going up in the order of accuracy, the \texttt{KL6} and $\ZDS$ scheme with $R=2$} present the same amount of error, whereas the \texttt{KL8} scheme provides the best accuracy compared to its competitor $\ZDS$ with $R=3$. We recover the same accuracy by using the $\ZDS$ scheme with $R=4$. This is one of the main advantages of the $\ZDS$ scheme: the order is easily increased, simply by increasing the block size. In the present case, we have a substantial gain of three orders of magnitude between $R=3$ and $R=4$.

\begin{figure}[!ht]
    \centering
    \includegraphics{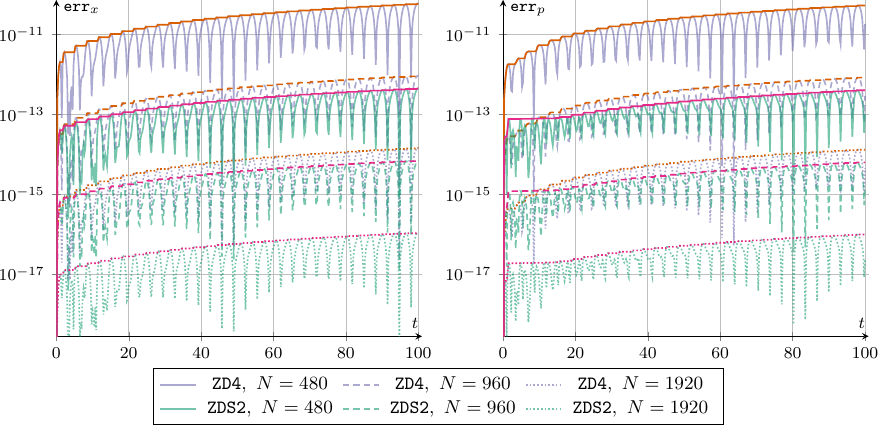}
    \caption{\rempolice \rb{Pendulum problem from \cref{sec:pendulum}: error on the position (left) and momentum (right) over time for the $\ZD4$ (blue lines) and $\ZDS2$ (green lines) schemes and several values of $N$, until $T=100$. To illustrate the growth of the error, the rolling maximal errors are displayed in orange for $\ZD4$ and in magenta for $\ZDS2$.}}
    \label{fig:HAM_pendulum_x_p}
\end{figure}

\rb{To further study the growth of the error on $x$ and $p$, we display in \cref{fig:HAM_pendulum_x_p} the errors on the position and momentum over time for the $\ZD4$ and $\ZDS2$ schemes, until the final time $T = 100$. For each value of $N$ (480, 960 and 1920), we observe a linear growth of the error in both $x$ and $p$ for both schemes. This is in line with the fact that symplectic integrators exhibit a linear error growth over time \cite{HaiLubWan2006}.}

We then carried out the simulation to $T=100\,000$ seconds to check the invariance of the Hamiltonian for a very long-time simulation. We use both a coarse grid with $N=3T$ (an intermediate between $N=2.4T$ and $N=4.8T$ we used for the case $T=100$), which corresponds to about $20$ point for a full revolution $T_p\approx6.53$, and a finer grid with $N=12T$ points, which makes it possible to compare the deviation of $\mathcal{H}$ with respect to the time step. We report in \cref{tab:pendulum_LongTime_accuracy} the errors for the $\ZD$, $\ZDS$ and classical schemes, respectively. The errors are in line with the case $T=100$ and the orders are optimal ones.

\begin{table}[!ht]\notapolice\centering
    \caption{\rempolice Pendulum system from \cref{sec:pendulum}: \rb{error on the Hamiltonian} at $T=100\,000$s. Both rows correspond to $N=T\times 3$ and $N=T\times 12$ respectively.}
    \label{tab:pendulum_LongTime_accuracy}
    \begin{subtable}{0.49\textwidth}\centering
        \caption{\rempolice Errors obtained with the $\ZD$ scheme.}
        \label{tab:pendulum_ZD_LongTime_accuracy}

        \begin{tabular}{@{}l@{}  r@{}l@{}  r@{}l@{}  r@{}l@{}  r@{}l@{}}
            \toprule
                             & \phantom{aaa} & {R=2}
                             & \phantom{aaa} & {R=4}
                             & \phantom{aaa} & {R=6}
                             & \phantom{aaa} & {R=8}                                              \\
            \midrule
            \eH ($\times3$)  &               & 9.19e-03 &  & 2.51e-04 &  & 1.62e-04 &  & 9.60e-05 \\
            \eH ($\times12$) &               & 3.92e-05 &  & 2.58e-08 &  & 1.04e-09 &  & 7.97e-11 \\
            \bottomrule
        \end{tabular}

    \end{subtable}
    \begin{subtable}{0.49\textwidth}\centering
        \caption{\rempolice Errors obtained with the $\ZDS$ scheme.}
        \label{tab:pendulum_ZDS_LongTime_accuracy}

        \begin{tabular}{@{}l@{}  r@{}l@{}  r@{}l@{}  r@{}l@{}  r@{}l@{}}
            \toprule
                             & \phantom{aaa} & {R=1}
                             & \phantom{aaa} & {R=2}
                             & \phantom{aaa} & {R=3}
                             & \phantom{aaa} & {R=4}                                              \\
            \midrule
            \eH ($\times3$)  &               & 4.72e-05 &  & 3.78e-06 &  & 8.72e-07 &  & 1.89e-07 \\
            \eH ($\times12$) &               & 1.85e-07 &  & 8.70e-10 &  & 5.39e-12 &  & 5.28e-14 \\
            \bottomrule
        \end{tabular}

    \end{subtable}

    \medskip

    \begin{subtable}{\textwidth}\centering
        \caption{\rempolice Errors obtained with the classical schemes.}
        \label{tab:pendulum_classical_LongTime_accuracy}

        \begin{tabular}{@{}l@{}  r@{}l@{}  r@{}l@{}  r@{}l@{}  r@{}l@{}}
            \toprule
                             & \phantom{aaa} & {MA2}
                             & \phantom{aaa} & {CS4}
                             & \phantom{aaa} & {KL6}
                             & \phantom{aaa} & {KL8}                                                   \\
            \midrule
            \eH ($\times3$)  &               & 9.01e-04 &  & \rb{1.58e-06} &  & 1.42e-07 &  & 7.04e-10 \\
            \eH ($\times12$) &               & 2.48e-05 &  & \rb{5.97e-09} &  & 3.27e-11 &  & 1.05e-14 \\
            \bottomrule
        \end{tabular}

    \end{subtable}
\end{table}

\begin{figure}[!ht]
    \centering
    \includegraphics[width=0.8\textwidth]{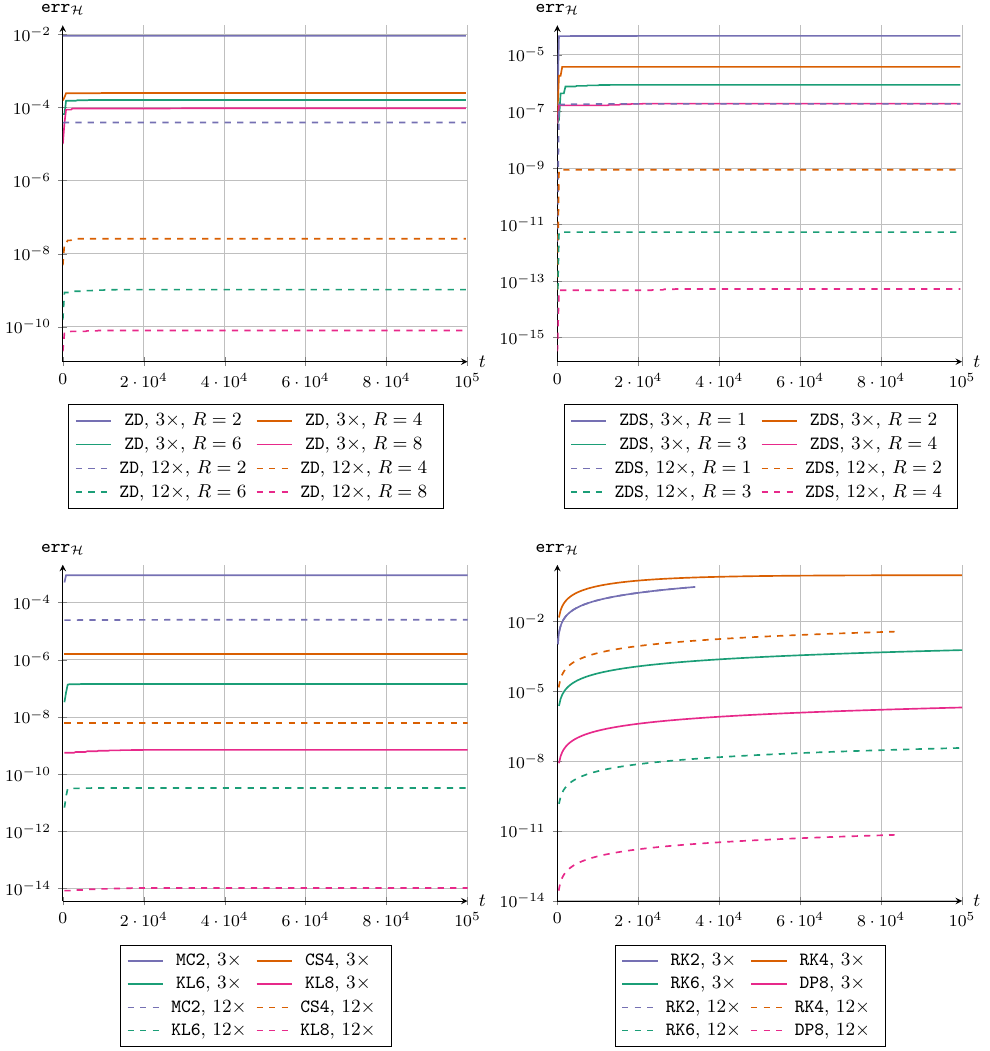}
    \caption{\rempolice Pendulum problem from \cref{sec:pendulum}: error on the Hamiltonian over time. Top panels: \ZD\ (left) and \ZDS\ (right) methods; bottom panels: classical symplectic (left) and non-symplectic (right) schemes.}
    \label{fig:HAM_pendulum_longtime}
\end{figure}

We confirm in \cref{fig:HAM_pendulum_longtime} the strict invariance in time of the errors. We also plot, in the bottom right panel, several examples of non-symplectic schemes (classical Runge-Kutta methods, see e.g. \cite{HaiLubWan2006}) to highlight their major drawback, namely a linear increase in the error over time.

%
%
%
%
\subsection{Two-dimensional Kepler problem}
\label{sec:2D_Kepler}

The two-dimensional Kepler system consists in a particle moving around a fixed point. It is characterized by the Hamiltonian
\begin{equation}\label{eq::hamiltonian_Kepler}
    \mathcal{H}(\bx, \bp) = \displaystyle \frac{\|\bp\|^2}{2} - \frac{1}{\|\bx\|},
\end{equation}
corresponding to the total mechanical energy.
The particle position $\bx(t)\in\mathbb R^2$ and the momentum $\bp(t)\in \mathbb R^2$ describe a quadratic curve (elliptic, parabolic or hyperbolic) depending on the initial total energy.

\subsubsection{Physical equations and invariants}
We first derive the dynamical system associated to the Hamiltonian \eqref{eq::hamiltonian_Kepler}
$$
    \dot\bx=\bp,\qquad  \dot\bp=-\frac{\bx}{\|\bx\|^3},
$$
and, differentiating with respect to time, we obtain the relations with the second derivatives:
$$
    \ddot \bx=\dot\bp,\qquad  \ddot\bp=-\frac{\bp}{\|\bx\|^3}+3\bx\frac{\bp\cdot \bx}{\|\bx\|^5}.
$$
We deduce the first group of physical equations \PE[1], given by
$$
    \bDx=\bZp,\qquad \bDp=-\frac{\bZx}{\|\bZx\|^3},
$$
while the second group of physical equation \PE[2] reads
$$
    \bSx=\bDp,\qquad \bSp=-\frac{\bZp}{\|\bZx\|^3}+3\bZx\frac{\bZx\cdot \bZp}{\|\bZx\|^5}.
$$

In addition to the Hamiltonian, there are two other invariant quantities for the Kepler system, namely:
\begin{enumerate}
    \item the angular momentum
          $$
              \mathcal{L}(\bx, \bp) = \bx \times \bp = p_2 x_1 - p_1 x_2;
          $$
    \item the Laplace-Runge-Lenz (LRL) vector
          $$
              \mathcal{A}(\bx, \bp) =  \mathcal{L}(\bx, \bp) \bp^\perp + \hat{\bx}
              \text{, \qquad with }
              \hat{\bx} = \frac{\bx}{\|\bx\|},
          $$
          from which we extract a one-dimensional invariant by summing the two components of $\mathcal{A}(\bx, \bp)$, yielding
          $$
              \mathcal{R}(\bx, \bp)=(p_2 x_1 - p_1 x_2)(p_1-p_2)+\frac{x_1+x_2}{\|\bx\|}.
          $$
\end{enumerate}

\subsubsection{Simulations}\label{sec:uncorrected_kepler}
Simulations are carried out until a final time $T=100$ seconds with the initial conditions $\bx(0) = (0.4, 0)$ and $\bp(0) = (0, 2)$. Coarse and fine meshes, with $N=24T$ and $N=96T$ respectively, are used to check the preservation of the invariants.
We denote by $\ea$ and $\eL$ stands the errors on the angular momentum  and LRL vector.

\begin{table}[!ht]\notapolice\centering
    \caption{\rempolice Kepler system from \cref{sec:2D_Kepler}: errors on the Hamiltonian, \ra{angular} momentum and Laplace-Runge-Lenz vector at $T=100$s. For $\eH$, $\eL$ and $\ea$, both rows correspond to $N=T\times 24$ and $N=T\times 96$ respectively.}
    \label{tab:Kepler_normal}
    \begin{subtable}{0.49\textwidth}\centering
        \caption{\rempolice Errors obtained with the $\ZD$ scheme.}
        \label{tab:Kepler_ZD_normal}

        \begin{tabular}{@{}l@{}  r@{}l@{}  r@{}l@{}  r@{}l@{}  r@{}l@{}}
            \toprule
                             & \phantom{aaa} & {R=2}
                             & \phantom{aaa} & {R=4}
                             & \phantom{aaa} & {R=6}
                             & \phantom{aaa} & {R=8}                                              \\
            \midrule
            \eH ($\times24$) &               & 4.15e-05 &  & 2.01e-05 &  & 1.52e-05 &  & 2.52e-05 \\
            \eH ($\times96$) &               & 1.81e-07 &  & 3.83e-09 &  & 1.06e-10 &  & 4.58e-12 \\
            \midrule
            \eL ($\times24$) &               & 3.25e-05 &  & 6.40e-06 &  & 1.41e-05 &  & 1.14e-05 \\
            \eL ($\times96$) &               & 1.28e-07 &  & 1.36e-09 &  & 3.89e-11 &  & 2.45e-12 \\
            \midrule
            \ea ($\times24$) &               & 4.54e-03 &  & 3.83e-04 &  & 8.63e-05 &  & 4.08e-05 \\
            \ea ($\times96$) &               & 1.82e-05 &  & 1.06e-07 &  & 1.29e-09 &  & 2.80e-11 \\
            \bottomrule
        \end{tabular}

    \end{subtable}
    \begin{subtable}{0.49\textwidth}\centering
        \caption{\rempolice Errors obtained with the $\ZDS$ scheme.}
        \label{tab:Kepler_ZDS_normal}

        \begin{tabular}{@{}l@{}  r@{}l@{}  r@{}l@{}  r@{}l@{}  r@{}l@{}}
            \toprule
                             & \phantom{aaa} & {R=1}
                             & \phantom{aaa} & {R=2}
                             & \phantom{aaa} & {R=3}
                             & \phantom{aaa} & {R=4}                                              \\
            \midrule
            \eH ($\times24$) &               & 4.83e-05 &  & 1.97e-06 &  & 3.08e-07 &  & 5.96e-08 \\
            \eH ($\times96$) &               & 1.88e-07 &  & 4.47e-10 &  & 2.47e-12 &  & 2.73e-14 \\
            \midrule
            \eL ($\times24$) &               & 1.27e-05 &  & 4.62e-07 &  & 6.02e-08 &  & 1.04e-08 \\
            \eL ($\times96$) &               & 4.93e-08 &  & 1.06e-10 &  & 5.62e-13 &  & 6.02e-15 \\
            \midrule
            \ea ($\times24$) &               & 3.26e-04 &  & 5.07e-06 &  & 3.33e-07 &  & 5.65e-08 \\
            \ea ($\times96$) &               & 1.28e-06 &  & 1.25e-09 &  & 3.52e-12 &  & 2.67e-14 \\
            \bottomrule
        \end{tabular}

    \end{subtable}

    \medskip

    \begin{subtable}{\textwidth}\centering
        \caption{\rempolice Errors obtained with the classical schemes.}
        \label{tab:Kepler_classical_normal}

        \begin{tabular}{@{}l@{}  r@{}l@{}  r@{}l@{}  r@{}l@{}  r@{}l@{}}
            \toprule
                             & \phantom{aaa} & {MC2}
                             & \phantom{aaa} & {CS4}
                             & \phantom{aaa} & {KL6}
                             & \phantom{aaa} & {KL8}                                                   \\
            \midrule
            \eH ($\times24$) &               & 3.56e-03 &  & \rb{4.82e-07} &  & 7.27e-09 &  & 3.67e-12 \\
            \eH ($\times96$) &               & 2.22e-04 &  & \rb{1.92e-09} &  & 1.73e-12 &  & 5.69e-17 \\
            \midrule
            \eL ($\times24$) &               & 7.95e-76 &  & 3.11e-76      &  & 7.00e-76 &  & 1.91e-75 \\
            \eL ($\times96$) &               & 7.17e-76 &  & 5.96e-76      &  & 1.00e-75 &  & 1.11e-75 \\
            \midrule
            \ea ($\times24$) &               & 5.23e-02 &  & \rb{5.81e-05} &  & 5.13e-07 &  & 1.49e-10 \\
            \ea ($\times96$) &               & 3.16e-03 &  & \rb{2.28e-07} &  & 1.26e-10 &  & 2.29e-15 \\
            \bottomrule
        \end{tabular}

    \end{subtable}
\end{table}

We report in \cref{tab:Kepler_normal} the deviations of the three invariants using the $\ZD$, $\ZDS$ and classical methods, respectively. First, we note that the $\mathcal H$ errors are now mesh-dependent, and decrease when the time step decreases. Second, the $\ZDS$ schemes have smaller error than their equivalent-order $\ZD$ versions (about 2 orders of magnitude). \rbbis{Comparisons with the classical symplectic schemes show that both structural methods provide larger errors than the \texttt{CS4} scheme.} On the contrary, the \texttt{KL6} (6th-order) and \texttt{KL8} (8th-order) schemes clearly provide the lowest $\mathcal H$ errors in comparison with the $\ZD\,R=6$ scheme (6th-order) and the $\ZDS\,R=3$ scheme (8th-order). We report a similar behavior for the LRL $\mathcal{R}$ invariant with a particular mention for the \texttt{KL8} method with a noticeable gain of three orders of magnitude regarding the equivalent $\ZDS\,R=3$ structural scheme.
Another comment concerns the conservation of the angular momentum $\mathcal L$. The structural schemes provide errors with the same magnitude as the Hamiltonian invariant, but the classical schemes definitively deliver a perfect invariance where the errors are only the consequence of the machine errors due to the {\tt quadruple} precision.

\subsubsection{Projection on the third invariant manifold}
\label{sec:Kepler_proj}

Since the last invariant is not fully preserved, we revisit the scheme by adding a projection onto the third invariant manifold characterized by the equation
$$
    \mathcal{R}(\bx,\bp)=\mathcal{R}(\bx_0, \bp_0)=\mathcal{R}_0.
$$
In short, we provide new vectors $\widetilde{\bZx}_{n}$ and $\widetilde{\bZp}_{n}$ close to $\bZx_{n}$ and $\bZp_{n}$ such that $\mathcal{R}(\widetilde{\bZx}_{n},\widetilde{\bZp}_{n})=\mathcal{R}_0$.

To this end, given $\bx$ and $\bp$ assumed to be close to the manifold, we consider the optimization problem: find $\widetilde \bx$ and $\widetilde \bp$ that minimize
$$
    \mathcal J= \frac{1}{2}\Big \Vert [\widetilde \bx,\widetilde \bp]^\trans-[\bx,\bp]^\trans\Big \Vert^2, \text{ \quad under the constraint }
    \mathcal{R}(\widetilde \bx, \widetilde \bp)=\mathcal{R}_0.
$$
We obtain an update $\widetilde \bx,\widetilde \bp$ given by
\begin{equation*}
    \begin{bmatrix}
        \widetilde x_1 \\\widetilde x_2\\\widetilde p_1\\\widetilde p_2
    \end{bmatrix} =
    \begin{bmatrix}
        x_1 \\x_2\\p_1\\p_2
    \end{bmatrix}
    -\lambda \displaystyle \nabla_{x,p} \mathcal {R}(\widetilde \bx,\widetilde \bp), \textrm{ with } \mathcal {R}(\widetilde \bx,\widetilde \bp)=\mathcal{R}_0,
\end{equation*}
where
$$
    \nabla_{x,p} \mathcal {R}(\bx,\bp)=\begin{bmatrix}
        +p_2(p_2-p_1)-\dfrac{(x_2)^2}{||x||^3}      \\
        -p_1(p_2-p_1)-\dfrac{(x_1)^2}{||x||^3}      \\
        2x_2p_1-x_1p_2-x_2p_2 \vphantom{\dfrac 1 2} \\
        2x_1p_2-x_1p_1-x_2p_1 \vphantom{\dfrac 1 2}
    \end{bmatrix}.
$$
We then use the following approximation
$$
    0=\mathcal R\Big (\bx-\lambda \nabla_x \mathcal{R}(\widetilde \bx,\widetilde \bp)\,,\, \bp-\lambda \nabla_p \mathcal{R}(\widetilde \bx,\widetilde \bp)\Big)-\mathcal R_0
    \approx \mathcal R(\bx,\bp)-\mathcal R_0-\lambda \Big <\nabla_{x,p}\mathcal{R}(\widetilde \bx,\widetilde \bp),\nabla_{x,p}\mathcal{R}(\bx,\bp)\Big >.
$$
Substituting $\nabla_{x,p}\mathcal{R}(\widetilde \bx,\widetilde \bp)$ with $\nabla_{x,p}\mathcal{R}(\bx,\bp)$ in the inner product, and we get the approximation
$$
    \lambda=\frac{\mathcal{R}(\bx,\bp)-\mathcal{R}_0}{\Vert \nabla_{x,p}\mathcal {R}(\bx,\bp) \Vert^2}.
$$
Thus, the approximation of the projection onto the third manifold is given by
$$
    \begin{bmatrix} \widetilde \bx\\\widetilde \bp \end{bmatrix}=
    \begin{bmatrix} \bx\\\bp \end{bmatrix}-\frac{\mathcal{R}(\bx,\bp)-\mathcal{R}_0}{\Vert \nabla_{x,p}\mathcal {R}(\bx,\bp) \Vert^2}
    \displaystyle \begin{bmatrix} \nabla_x \mathcal{R}(\bx,\bp)\\[0.3em]\nabla_p \mathcal{R}(\bx,\bp) \end{bmatrix}.
$$

We report in \cref{tab:Kepler_update} the deviation of the invariants at the final time $T=100$ using the projection. We obtain very similar errors to the non-projection case given in \cref{tab:Kepler_ZD_normal},\cref{tab:Kepler_ZDS_normal}. The correction does not bring significant change in the error magnitude for the $\ZDS$ method. We just observe a slight reduction of the LRL error and larger errors for the two other invariants for the $\ZD$ method.

\begin{table}[!ht]\notapolice\centering
    \caption{\rempolice Kepler system from \cref{sec:2D_Kepler}, with projection on the invariant LRL manifold: errors on the Hamiltonian, \ra{angular} momentum and Laplace-Runge-Lenz vector at $T=100$s. For $\eH$, $\eL$ and $\ea$, both rows correspond to $N=T\times 24$ and $N=T\times 96$ respectively.}
    \label{tab:Kepler_update}
    \begin{subtable}{0.49\textwidth}\centering
        \caption{\rempolice Errors obtained with the $\ZD$ scheme.}
        \label{tab:Kepler_ZD_update}

        \begin{tabular}{@{}l@{}  r@{}l@{}  r@{}l@{}  r@{}l@{}  r@{}l@{}}
            \toprule
                             & \phantom{aaa} & {R=2}
                             & \phantom{aaa} & {R=4}
                             & \phantom{aaa} & {R=6}
                             & \phantom{aaa} & {R=8}                                              \\
            \midrule
            \eH ($\times24$) &               & 3.75e-03 &  & 3.88e-04 &  & 1.19e-04 &  & 2.16e-05 \\
            \eH ($\times96$) &               & 8.78e-05 &  & 1.89e-07 &  & 4.21e-09 &  & 1.13e-10 \\
            \midrule
            \eL ($\times24$) &               & 1.40e-03 &  & 5.68e-05 &  & 3.38e-05 &  & 1.48e-05 \\
            \eL ($\times96$) &               & 5.16e-05 &  & 8.48e-08 &  & 1.41e-09 &  & 4.19e-11 \\
            \midrule
            \ea ($\times24$) &               & 1.39e-04 &  & 3.70e-05 &  & 2.15e-05 &  & 2.66e-05 \\
            \ea ($\times96$) &               & 2.97e-06 &  & 9.69e-09 &  & 3.95e-11 &  & 3.70e-12 \\
            \bottomrule
        \end{tabular}

    \end{subtable}
    \begin{subtable}{0.49\textwidth}\centering
        \caption{\rempolice Errors obtained with the $\ZDS$ scheme.}
        \label{tab:Kepler_ZDS_update}

        \begin{tabular}{@{}l@{}  r@{}l@{}  r@{}l@{}  r@{}l@{}  r@{}l@{}}
            \toprule
                             & \phantom{aaa} & {R=1}
                             & \phantom{aaa} & {R=2}
                             & \phantom{aaa} & {R=3}
                             & \phantom{aaa} & {R=4}                                              \\
            \midrule
            \eH ($\times24$) &               & 3.17e-03 &  & 1.66e-05 &  & 2.07e-07 &  & 9.13e-08 \\
            \eH ($\times96$) &               & 1.25e-03 &  & 5.01e-08 &  & 1.57e-10 &  & 6.59e-13 \\
            \midrule
            \eL ($\times24$) &               & 1.28e-03 &  & 2.70e-06 &  & 2.27e-07 &  & 2.25e-08 \\
            \eL ($\times96$) &               & 7.22e-04 &  & 2.78e-08 &  & 7.71e-11 &  & 2.88e-13 \\
            \midrule
            \ea ($\times24$) &               & 9.87e-05 &  & 6.85e-07 &  & 9.61e-08 &  & 4.61e-08 \\
            \ea ($\times96$) &               & 2.45e-05 &  & 1.37e-09 &  & 4.71e-12 &  & 1.93e-14 \\
            \bottomrule
        \end{tabular}

    \end{subtable}
\end{table}

\rbbis{
    \subsection{Two-dimensional three-body problem: the figure-eight solution}
    \label{sec:2D_3body}
    The two-dimensional three-body system is characterized by the position $X[k]=\bx^k\in\mathbb R^2$ and the momentum $P[k]=\bp^k=m_k\bv^k\in\mathbb R^2$ of each body, where matrices $X,P \in \Kbody$ are given by $X_{i,k}=\bx^k_i$ and $P_{i,k}=\bp^k_i$, for $k \in \{1,2,3\}$ and $i \in \{1,2\}$. The Hamiltonian corresponds to the time-invariant total mechanical energy
    \begin{equation*}\label{eq::HAM-3-Bodies}
        \mathcal{H}(X, P) = \displaystyle \sum_{k=1}^K \frac {1}{2 m_k} \|\bp^k\|^2 - \sum_{k\neq \ell} \frac {G m_k m_\ell} {2\|\bx^k - \bx^\ell\|}.
    \end{equation*}
    We choose the figure-eight orbit
    for the three-body problem ($K=3$) given by \cite{Moore1993}. To simplify the problem, we assume that all masses are equal to $1$ and take the gravity constant $G=1$. We prescribe the initial conditions with
    \begin{equation*}
        x^1=\begin{pmatrix}0.97000436 \\ -0.24308753\end{pmatrix},\quad
        x^2=\begin{pmatrix}- 0.97000436 \\ 0.24308753\end{pmatrix},\quad
        x^3=\begin{pmatrix}0 \\ 0\end{pmatrix},
    \end{equation*}
    for the position while the initial velocities are given by
    \begin{equation*}
        v^1=\begin{pmatrix}0.466203685 \\ 0.43236573\end{pmatrix},\quad
        v^2=\begin{pmatrix}0.466203685 \\ 0.43236573\end{pmatrix},\quad
        v^3=\begin{pmatrix}-0.93240737 \\ -0.86473146\end{pmatrix}.
    \end{equation*}
}%
\rbbis{
    To check the accuracy for both the position and the velocity, we build a 32 digit-accurate approximation at time $T=6$, reported in \cref{tab::32digitT6}, with a very fine mesh $N=10000$, $16$th-order method and octuple precision.
    \begin{table}[ht]\notapolice\centering
        \caption{\rempolice 32 digit-accurate approximations of the final position and velocity at time T=6 for the two-dimensional three-body problem?}
        \begin{tabular}{l|rr}
                  & \multicolumn{1}{c}{component 1}         & \multicolumn{1}{c}{component 2}         \\
            \hline
            $x^1$ & 7.46327477619643171788882241806705e-01  & -3.50642825804179469824366576996967e-01 \\
            $v^1$ & 9.15234708896112323200536200718399e-01  & 1.58591215936613597397887314351545e-01  \\
            \hline
            $v^2$ & -1.06634604893325246853521993232162e+00 & 9.40324145269491260952974444713522e-02  \\
            $v^2$ & 1.48379313616329367267908182363434e-01  & 4.67056998997661829135979976473009e-01  \\
            \hline
            $x^3$ & 3.20018571313609296746337690514915e-01  & 2.56610411277230343729069132525615e-01  \\
            $v^3$ & -1.06361402251244169046844438308183e+00 & -6.25648214934275426533867290824555e-01 \\
            \hline
        \end{tabular}
        \label{tab::32digitT6}
    \end{table}
}

\rbbis{
    We plot in \cref{fig:L1-error-Zx-Zp}, left panel, the sum of the position errors for the three bodies while the right panel provides the error summation for the velocity. Methods SV (full line) and ZD (dashed line) provides similar accuracies for the same order but from our experience, the ZD method runs much faster. When dealing with the ZDS method, we reach the optimal order but observe a gain in accuracy of up to four orders of magnitude due to the compact nature of the scheme. The introduction of the second derivative as an additional unknown together with the new physical equation strongly reduce the errors with a modest extra computational cost.
    \begin{figure}
        \centering
        \includegraphics[width=0.48\linewidth]{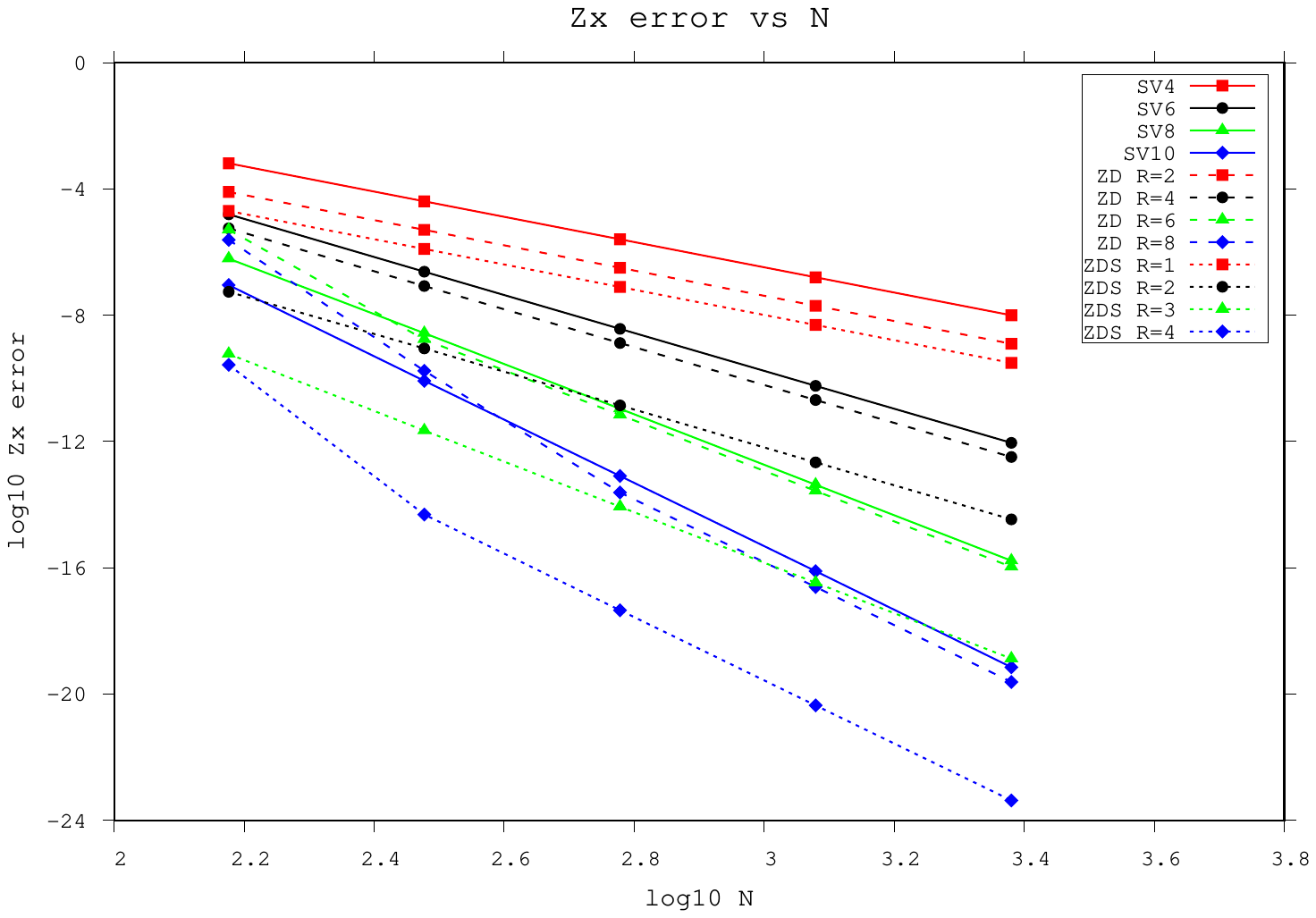}
        \hfill
        \includegraphics[width=0.48\linewidth]{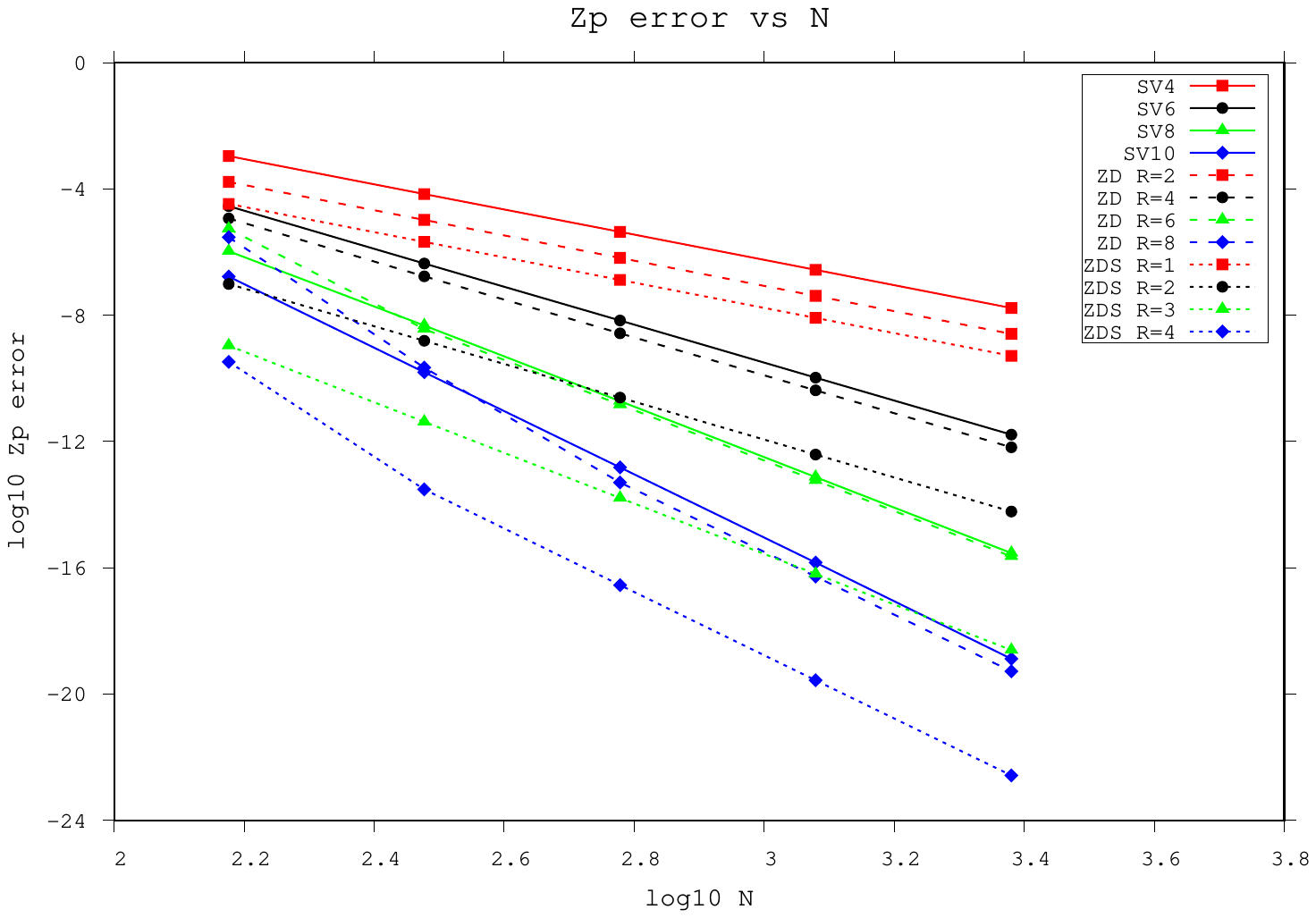}
        \caption{\rempolice $L^1$ errors and convergence orders for the 2-dimension 3-body problem.  Position errors (left panel) and the velocity errors (right panel) are given in function of the grid size for several schemes}
        \label{fig:L1-error-Zx-Zp}
    \end{figure}
}

%
%
%
%
\subsection{The outer Solar system}
\label{sec:solar_system}

{\rbparagraphbis
    The three-dimensional $K$-body system is characterized by the position $X[k]=\bx^k\in\mathbb R^3$ and the momentum $P[k]=\bp^k=m_k\bv^k\in\mathbb R^3$ of each body, where matrices $X,P \in \Kbody$ are given by $X_{i,k}=\bx^k_i$ and $P_{i,k}=\bp^k_i$, for all $k \in \{1,\dots,K\}$ and $i \in \{1,2,3\}$. The Hamiltonian corresponds to the time-invariant total mechanical energy
    \begin{equation*}\label{eq::HAM-K-Bodies}
        \mathcal{H}(X, P) = \displaystyle \sum_{k=1}^K \frac {1}{2 m_k} \|\bp^k\|^2 - \sum_{k\neq \ell} \frac {G m_k m_\ell} {2\|\bx^k - \bx^\ell\|}.
    \end{equation*}

    \subsubsection{Physical equations and invariants}

    Hamilton's equations $\dot X=\nabla_P \mathcal{H}(X,P)$ and $\dot P=-\nabla_X \mathcal{H}(X,P)$ give the first set of physical equations:
    $$
        \text{for all } k \in \{1, \dots, K\},
        \qquad
        \dot \bx^k=\frac{\bp^k}{m_k},
        \qquad
        \dot  \bp^k= -\sum_{\ell\neq k} Gm_km_\ell\frac{\bx^k-\bx^\ell}{\|\bx^k-\bx^\ell\|^3}.
    $$
    Hence, the first set of physical equations $\PE[1]$ given by relations \eqref{body_PE1_x} -- \eqref{body_PE1_p} derives from Hamilton's principle:
    \begin{eqnarray}
        &&\displaystyle \DX[k] =\frac{\ZP[k]}{m_k},\label{body_PE1_x}\\
        &&\displaystyle \DP[k]=-\sum_{\ell\neq k} G m_km_\ell\frac{\ZX[\ell]-\ZX[k]}{\big \|\ZX[\ell]-\ZX[k]\big\|^3}.\label{body_PE1_p}
    \end{eqnarray}
    Differentiating the dynamical system once more provides the second set of physical equations $\PE[2]$, given by relations \eqref{body_PE2_x} -- \eqref{body_PE2_p}:
    \begin{align}
        \label{body_PE2_x}
        \SP[k] & =\frac{\DP[k]}{m_k},                                                                              \\
        \label{body_PE2_p}
        \SP[k] & =-\sum_{\ell\neq k} G m_k m_\ell\Bigg ( \frac{\DX[\ell]-\DX[k]}{\big \|\ZX[\ell]-\ZX[k]\big\|^3}-
        3\Big\langle\ZX[k]-\ZX[\ell], \DX[k]-\DX[\ell]\Big\rangle\,\frac{\ZX[k]-\ZX[\ell]}{\big \|\ZX[k]-\ZX[\ell]\big\|^5}  \Bigg ),
    \end{align}
    with $\langle \cdot, \cdot\rangle$ the usual inner product, while $\ZX[k]$, $\DX[k]$ and $\SX[k]$ represent the position, the velocity and the acceleration of body $k$.
    A second invariant of the $K$-body system is the angular momentum $\mathcal{L}(X,P)$, given by
    $$
        \mathcal{L}(X,P) = \displaystyle \sum_{k=1}^K \mathcal{L}_k(X[k],P[k]) = \displaystyle \sum_{k=1}^K \bx^k \times \bp^k.
    $$

    \subsubsection{Numerical tests}

    As an application, we consider the $6$-body system corresponding to the so-called outer solar system: the Sun and Jupiter, Saturn, Uranus, Neptune and Pluto \cite{HaiLubWan2006}.}
The masses and initial conditions are given \footnote{\tt \url{https://dspace.mit.edu/bitstream/handle/1721.1/6442/AIM-877.pdf}}$^{,}$\footnote{\tt \url{https://rebound.readthedocs.io/en/latest/c\_examples/outer\_solar\_system/}} in \cref{tab::data_outer_system_initialization}. The initial positions are given in Astronomical Units (1 au=149 597 870 km) and the initial velocities are given in $au$ per earth day. The gravity constant is set to $G = 2.95912208286 \cdot 10^{-4}$.
Together with the Hamiltonian, we analyse the angular momentum, another quantity that is time-invariant, given by the vector
$$
    \ba=\mathcal{L}(X,P) = \displaystyle \sum_{k=1}^K \mathcal{L}_k(X[k],P[k]) = \displaystyle \sum_{k=1}^K \bx^k \times \bp^k.
$$


\begin{table}
    \centering
    \caption{\rempolice Data and initial conditions for the outer solar system problem. Masses are normalized regarding the solar mass, and the initial positions and velocities are given, respectively, in $au$ and $au$ per earth day.}
    \label{tab::data_outer_system_initialization}
    \medskip
    {\tablepolice 
        \begin{tabular}{crrrr}
            \toprule
            Celestial body & Mass               & Initial Position                                                                    & Initial Velocity                                       & period \\
            \midrule
            Sun            & 1.00000597682e+00  & \makecell[r]{0                                       \\ 0 \\ 0}                     & \makecell[r]{0 \\ 0 \\ 0}                              & 0      \\
            \midrule
            Jupiter        & 9.547861040430e-04 & \makecell[r]{-3.5023653                              \\ -3.8169847 \\ -1.5507963}   & \makecell[r]{0.00565429 \\ -0.00412490 \\ -0.00190589} & 4333   \\
            \midrule
            Saturn         & 2.855837331510e-04 & \makecell[r]{ 9.0755314                              \\ -3.0458353 \\ -1.6483708}   & \makecell[r]{0.00168318 \\ 0.00483525 \\ 0.00192462}   & 10759  \\
            \midrule
            Uranus         & 4.37273164546e-05  & \makecell[r]{ 8.3101420                              \\ -16.2901086 \\ -7.2521278}  & \makecell[r]{0.00354178 \\ 0.00137102 \\ 0.00055029}   & 30687  \\
            \midrule
            Neptune        & 5.17759138449e-05  & \makecell[r]{11.4707666                              \\ -25.7294829 \\ -10.8169456} & \makecell[r]{0.00288930 \\ 0.00114527 \\ 0.00039677}   & 60190  \\
            \midrule
            Pluto          & (10/13) e-08       & \makecell[r]{-15.5387357                             \\ -25.2225594 \\ -3.1902382}  & \makecell[r]{0.00276725 \\ -0.00170702 \\ -0.00136504} & 90560  \\
            \bottomrule
        \end{tabular}
    }
\end{table}

We provide in \cref{tab:OuterSolar_longTime_accuracy} the maximum errors for the Hamiltonian and \ra{angular} momentum using the $\ZD$, $\ZDS$ and classical schemes, with coarse and fine grids, and $T=100\,000$ years. In each situation, we obtain the optimal order for the $\ZD$, $\ZDS$ and classical methods, noting that we reach the {\tt quadruple} precision for the angular momentum with the symplectic scheme. We also observe the same gain of two or three orders of magnitude between the $\ZD$ and $\ZDS$ schemes. \rbbis{The \texttt{CS4} scheme (4th-order accurate) outperforms both of its structural equivalents.}
The \texttt{KL6} scheme (6th-order accurate) provides very similar deviations compared to its equivalent $\ZDS$ with $R=2$, while the \texttt{KL8} scheme (8th-order accurate) is comparable to the $\ZDS$ with $R=3$. The main advantage of the structural scheme is its ability to choose or adapt the order without modifying its implementation, only by changing the value of the block size $R$.

\begin{table}[!ht]\notapolice\centering
    \caption{\rempolice Three-dimensional $n$-body problem (i.e., the outer Solar system) from \cref{sec:solar_system}: errors on the Hamiltonian and angular momentum at $T=100\,000$ years. For $\eH$ and $\eL$, both rows correspond to $N=480$ and $N=1920$ respectively.}
    \label{tab:OuterSolar_longTime_accuracy}
    \begin{subtable}{0.49\textwidth}\centering
        \caption{\rempolice Errors obtained with the $\ZD$ scheme.}
        \label{tab:OuterSolar_ZD_longTime_accuracy}

        \begin{tabular}{@{}l@{}  r@{}l@{}  r@{}l@{}  r@{}l@{}  r@{}l@{}}
            \toprule
                         & \phantom{aaa} & {R=2}
                         & \phantom{aaa} & {R=4}
                         & \phantom{aaa} & {R=6}
                         & \phantom{aaa} & {R=8}                                              \\
            \midrule
            \eH ($480$)  &               & 3.06e-04 &  & 1.07e-04 &  & 4.02e-05 &  & 3.79e-05 \\
            \eH ($1920$) &               & 1.26e-06 &  & 2.72e-08 &  & 5.11e-10 &  & 2.29e-11 \\
            \midrule
            \eL ($480$)  &               & 9.12e-09 &  & 2.42e-09 &  & 9.37e-10 &  & 8.87e-10 \\
            \eL ($1920$) &               & 3.70e-11 &  & 6.15e-13 &  & 1.08e-14 &  & 5.46e-16 \\
            \bottomrule
        \end{tabular}

    \end{subtable}
    \begin{subtable}{0.49\textwidth}\centering
        \caption{\rempolice Errors obtained with the $\ZDS$ scheme.}
        \label{tab:OuterSolar_ZDS_longTime_accuracy}

        \begin{tabular}{@{}l@{}  r@{}l@{}  r@{}l@{}  r@{}l@{}  r@{}l@{}}
            \toprule
                         & \phantom{aaa} & {R=1}
                         & \phantom{aaa} & {R=2}
                         & \phantom{aaa} & {R=3}
                         & \phantom{aaa} & {R=4}                                              \\
            \midrule
            \eH ($480$)  &               & 7.92e-05 &  & 2.33e-06 &  & 9.99e-08 &  & 3.48e-08 \\
            \eH ($1920$) &               & 3.10e-07 &  & 5.69e-10 &  & 1.52e-12 &  & 1.03e-14 \\
            \midrule
            \eL ($480$)  &               & 1.76e-09 &  & 5.03e-11 &  & 2.63e-12 &  & 6.56e-13 \\
            \eL ($1920$) &               & 6.88e-12 &  & 1.23e-14 &  & 3.12e-17 &  & 2.36e-19 \\
            \bottomrule
        \end{tabular}

    \end{subtable}

    \medskip

    \begin{subtable}{\textwidth}\centering
        \caption{\rempolice Errors obtained with the classical schemes.}
        \label{tab:OuterSolar_classical_longTime_accuracy}

        \begin{tabular}{@{}l@{}  r@{}l@{}  r@{}l@{}  r@{}l@{}  r@{}l@{}}
            \toprule
                         & \phantom{aaa} & {MA2}
                         & \phantom{aaa} & {CS4}
                         & \phantom{aaa} & {KL6}
                         & \phantom{aaa} & {KL8}                                                   \\
            \midrule
            \eH ($480$)  &               & 1.24e-01 &  & \rb{6.19e-06} &  & 2.27e-05 &  & 3.29e-07 \\
            \eH ($1920$) &               & 7.43e-03 &  & \rb{2.46e-08} &  & 3.55e-09 &  & 2.71e-12 \\
            \midrule
            \eL ($480$)  &               & 7.91e-81 &  & 1.90e-80      &  & 2.06e-80 &  & 2.11e-80 \\
            \eL ($1920$) &               & 3.37e-80 &  & 2.21e-80      &  & 4.53e-80 &  & 1.17e-79 \\
            \bottomrule
        \end{tabular}

    \end{subtable}
\end{table}
%
%
%
%
\subsection{Motion of a particle in a 3D electromagnetic field}
\label{sec:3D_EM_field}

We consider a charged particle subjected to electromagnetic forces, characterized by the electric potential $\varphi=\varphi(\bx)\in\mathbb R$ and the magnetic potential $A=A(\bx)\in\mathbb R^3$.
This situation is governed by the following Hamiltonian:
\begin{equation*}
    \mathcal{H}(\bx, \bp) = \frac 1 {2 m} \| \bp - e A(\bx) \|^2 + e \varphi(\bx),
\end{equation*}
with $m$ and $e$ the mass and electric charge of the particle.
We denote by $A(\bx)=[A_1,A_2,A_3]^\trans$ the components of the magnetic potential.

This is an example of a non-separable Hamiltonian, i.e., a problem where $\nabla_\bp \mathcal{H}$ depends on both $\bx$ and $\bp$.
The classical schemes already used in the previous section and implemented in {\tt julia} failed to solve such a problem.
Hence, we have implemented the second-order symplectic \rb{St\"ormer}-Verlet schemes \cite{HaiLubWan2003}. The composition method \cite{HaiLubWan2006} make it possible to reach the fourth, sixth and eighth orders of accuracy, to provide comparisons with the structural method.

\subsubsection{Physical equations}
From the Hamiltonian, we derive the equations of motion given by
\begin{align}
    \dot\bx & = \partial_p \mathcal{H}(\bx,\bp) = \dfrac 1 m (\bp - e A), \label{EL::PE1a}                                     \\
    \dot\bp & =-\partial_x \mathcal{H}(\bx,\bp) = \frac{e}{m}\big [\partial_x A\big ]^\trans(\bp - e A) - e \partial_x \varphi 
    =e \left ( \big [\partial_x A\big ]^\trans\,\dfrac{dx}{dt} - \partial_x \varphi \right ),  \label{EL::PE1b}
\end{align}
with the matrix $\displaystyle \partial_x A(\bx)=\Big [ \partial_{j} A_{i}(\bx)\Big ]_{i,j} \in\mathbb R^{3\times 3}$ and the vector $\partial_x \varphi(\bx)=[\partial_1 \varphi(\bx),\partial_2 \varphi(\bx),\partial_3 \varphi(\bx)]^\trans \in\mathbb R^3$.

From relations \eqref{EL::PE1a} -- \eqref{EL::PE1b}, we deduce the first group of physical equations $(\PE[1])$  connecting the first-order derivatives
\begin{align*}
    m\bDx & =\bZp-eA(\bZx),  \quad
    \bDp   =e\left ( \big [\partial_x A(\bZx)\big ]^t\bDx-\partial_x \varphi(\bZx)  \right ), \\
\end{align*}
Differentiating relations \eqref{EL::PE1a} -- \eqref{EL::PE1b} with respect to time, we get the second group of physical equations
\begin{align*}
    m\ddot \bx & =\dot \bp-e\partial_x A \dot \bx-e\partial_t A \vphantom{\dfrac 1 2}
    = e \left ( \big [\partial_x A\big ]^\trans\,\dot \bx - \partial_x \varphi -\partial_x A \dot \bx-\partial_t A\right ),\vphantom{\dfrac 1 2} \\
    \ddot \bp  & =e\frac{d}{dt} \left (\big [\partial_x A\big ]^\trans\,\dot \bx \right )-e  \partial^2_x \varphi\,\dot \bx
    \vphantom{\dfrac 1 2}
    =e\Big ( \big [\partial_x A\big ]^\trans\, \ddot \bx + w- \partial^2_x \varphi\,\dot \bx \Big )\vphantom{\dfrac 1 2},
\end{align*}
with
$$
    \partial_x^2 \varphi=\big [\partial_i \partial_j \varphi\big ]_{i,j}
    \text{\qquad and, for all } \ell \in \{1, 2, 3\},
    w_\ell=\sum_{i=1}^3  \sum_{j=1}^3   \partial_\ell \partial_i \Big(A_j\Big ) \dot x_i \dot x_j.
$$
The second group of physical equations $(\PE[2])$ then reads
\begin{align*}
    \frac{m}{e}\bSx & =\Big (\big [\partial_x A(\bZx)\big ]^t-[\partial_x A(\bZx)\big ]\Big )\bDx-\partial_tA(\bZx)-\partial_x \phi(\bZx), \\
    \bSp            & =e \Big (\big [\partial_x A(\bZx)\big ]\bSx +w -\big[\partial^2_x \varphi(\bZx)\big] \bDx  \Big ).                   \\
\end{align*}

\subsubsection{A Sanity Check Benchmark (SCB)}
\label{sec:SCB}

The SCB consists in considering a pseudo-2D problem with very smooth magnetic and electric potentials.
We take
\begin{equation*}
    \varphi(\bx) = -\frac{1}{0.1+\|\bx\|},\quad
    A(\bx) =\begin{pmatrix}
        0 \\1000x_1 \\0\end{pmatrix}.
\end{equation*}
The initial conditions are $x(0) = (1, 0, 0)^\trans$ and $p(0) = (0, 1+100, 0)^\trans$, while the physical parameters are $m = 1$ and $q = 1$, with final time $T=100$.

We have carried out the simulation up to the time $T=100$ using a coarse grid ($N=12T$) and finer one ($N=48T$) to check the order of accuracy and the preservation of the Hamiltonian. We report in  \cref{tab:ORTF_accuracy} the Hamiltonian deviation for the $\ZD$, $\ZDS$ and  \rb{St\"ormer}-Verlet (SV) methods. Clearly, the structural schemes provide a better accuracy in comparison to the SV version. \rbbis{Namely, the $\ZD$ methods provide a gain of up to one order of magnitude compared to the classical SV methods, while we obtain a gain of two to three orders of magnitude with the $\ZDS$ approach.} We highlight the quality of the compact second-derivative method, able to reach such very low errors.

\begin{table}[!ht]\notapolice\centering
    \caption{\rempolice Particle in an electromagnetic field, Sanity Check Benchmark from \cref{sec:SCB}: error on the Hamiltonian at $T=100$s. Both rows correspond to $N=1200$ and $N=4800$, respectively.}
    \label{tab:ORTF_accuracy}
    \begin{subtable}{0.49\textwidth}\centering
        \caption{\rempolice Errors obtained with the $\ZD$ scheme.}
        \label{tab:ORTF_ZD_accuracy}

        \begin{tabular}{@{}l@{}  r@{}l@{}  r@{}l@{}  r@{}l@{}  r@{}l@{}}
            \toprule
                             & \phantom{aaa} & {R=2}
                             & \phantom{aaa} & {R=4}
                             & \phantom{aaa} & {R=6}
                             & \phantom{aaa} & {R=8}                                              \\
            \midrule
            \eH ($\times12$) &               & 8.69e-07 &  & 2.27e-08 &  & 9.53e-10 &  & 8.56e-11 \\
            \eH ($\times48$) &               & 3.41e-09 &  & 5.57e-12 &  & 1.43e-14 &  & 6.96e-17 \\
            \bottomrule
        \end{tabular}

    \end{subtable}
    \begin{subtable}{0.49\textwidth}\centering
        \caption{\rempolice Errors obtained with the $\ZDS$ scheme.}
        \label{tab:ORTF_ZDS_accuracy}

        \begin{tabular}{@{}l@{}  r@{}l@{}  r@{}l@{}  r@{}l@{}  r@{}l@{}}
            \toprule
                             & \phantom{aaa} & {R=1}
                             & \phantom{aaa} & {R=2}
                             & \phantom{aaa} & {R=3}
                             & \phantom{aaa} & {R=4}                                              \\
            \midrule
            \eH ($\times12$) &               & 2.89e-07 &  & 9.42e-10 &  & 7.37e-12 &  & 1.18e-13 \\
            \eH ($\times48$) &               & 1.13e-09 &  & 2.29e-13 &  & 8.42e-17 &  & 1.04e-19 \\
            \bottomrule
        \end{tabular}

    \end{subtable}

    \medskip

    \begin{subtable}{\textwidth}\centering
        \caption{\rempolice Errors obtained with the classical schemes.}
        \label{tab:ORTF_SV_accuracy}

        \begin{tabular}{@{}l@{}  r@{}l@{}  r@{}l@{}  r@{}l@{}  r@{}l@{}}
            \toprule
                             & \phantom{aaa} & {SV2}
                             & \phantom{aaa} & {SV4}
                             & \phantom{aaa} & {SV6}
                             & \phantom{aaa} & {SV8}                                              \\
            \midrule
            \eH ($\times12$) &               & 6.10e-04 &  & 4.66e-06 &  & 1.98e-08 &  & 2.07e-09 \\
            \eH ($\times48$) &               & 3.80e-05 &  & 1.82e-08 &  & 4.43e-12 &  & 3.23e-14 \\
            \bottomrule
        \end{tabular}

    \end{subtable}
\end{table}

\begin{figure}[!ht]
    \centering
    \includegraphics[width=\textwidth]{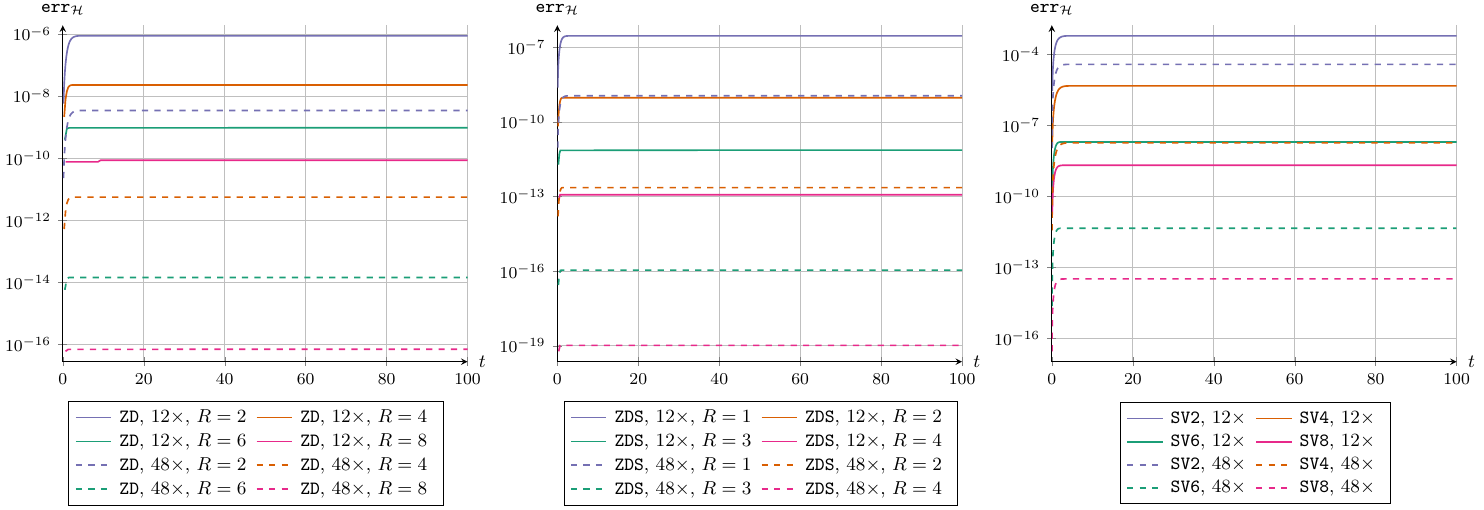}
    \caption{\rempolice Particle in an electromagnetic field, Sanity Check Benchmark from \cref{sec:SCB}: error on the Hamiltonian over time. From left to right: $\ZD$, $\ZDS$ and SV methods.}
    \label{fig:ORTF_HAM}
\end{figure}

\cref{fig:ORTF_HAM} displays the growth of the Hamiltonian error in time for the $\ZD$ (left panel), $\ZDS$ (middle panel) and SV (right panel) schemes. In each case, the Hamiltonian error becomes constant after a very short growth.

\subsubsection{Non-separable case: particle motion for challenging potentials}
\label{sec:particle_non_separable}

We propose a more challenging benchmark, where the magnetic potential presents a singularity at $x_1=0$, leading to additional numerical difficulties to preserve the Hamiltonian. The electrostatic potential is given by $\varphi(\bx) = 2 \cos(x_1)^2 + \sin(x_1)^2 (\sin(x_2) \cos(x_2) + \sin(x_3) \cos(x_3))$, while the magnetic potential reads
$$
    A(\bx) =\begin{pmatrix}
        r^2 \\r^2 \frac{x_2}{x_1} \\ -2 \log(1 + r^2)
    \end{pmatrix},
$$
with $r^2 = x_1^2 + x_2^2 + x_3^2$.
The mass and charge are set to unity, i.e., $m = 1$ and $q = 1$,
and the final simulation time is $20\,000\,s$.
At last, we use the initial conditions $\bx(0) = (0.5, - 0.25, -0.25)^\trans$ and $\bp(0) = (0, 0, -1)^\trans$.
Not that this potential leads to a truly non-separable problem.

\begin{figure}[!ht]
    \centering
    \includegraphics[width=\textwidth]{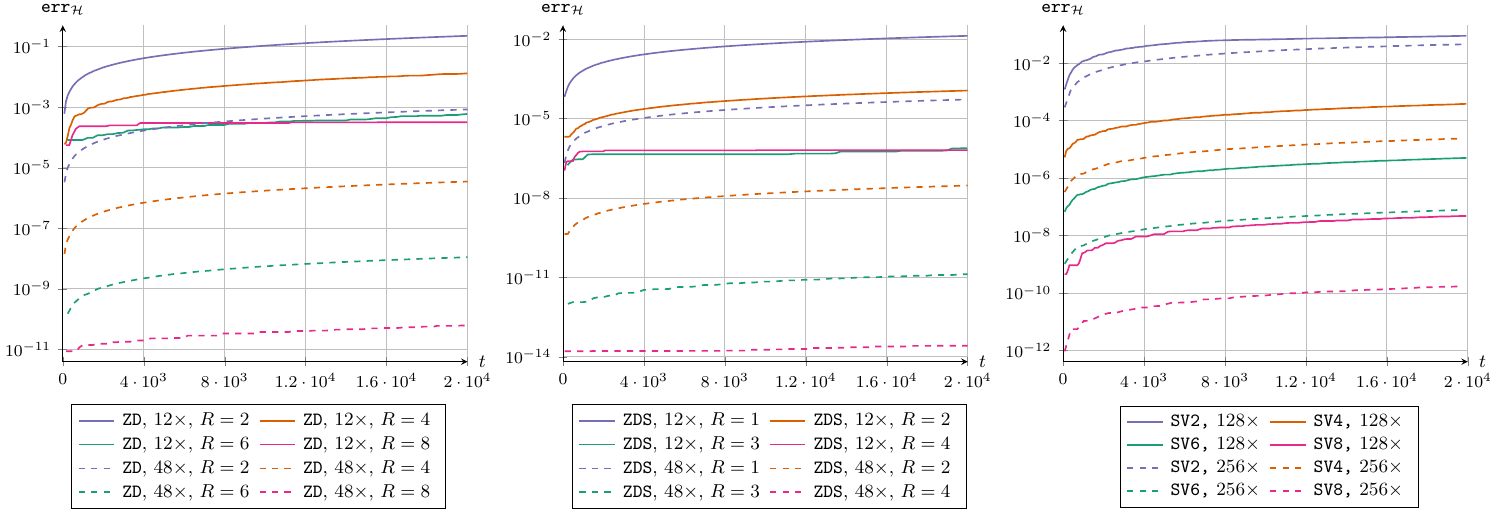}
    \caption{\rempolice Particle in an electromagnetic field, non-separable case from \cref{sec:particle_non_separable}: error on the Hamiltonian over time. From left to right: \ZD, \ZDS\ and classical methods.}
    \label{fig:Particle_HAM}
    \label{fig::HAM_electromagnetic_particule}
\end{figure}

We run the $\ZD$ and $\ZDS$ simulation for $N=12T$ and $N=48T$ whereas, for stability reason, the SV method required much finer grids with $N=128T$ and $N=256T$.
The results are reported in \cref{tab:Particle_ZD_accuracy,tab:Particle_ZDS_accuracy,tab:Particle_classical_accuracy}.
Once again, we observe the superior advantage of the structural scheme to achieve very high accuracy, even with coarser grids. As usual, we get a gain of about two orders of magnitude between the $\ZD$ and $\ZDS$ scheme for a given order of convergence. Moreover, we display in \cref{fig::HAM_electromagnetic_particule} the evolution of the deviation of the theoretically invariant Hamiltonian. On the one hand, the SV method presents a linear growth of the error as a linear function of the time, independently of the method order (almost the same slope in $\log$ scale). We observe a similar behavior for the fourth- and sixth-order structural methods. At last, we notice a dramatic cut of the error growth when dealing with very high orders of accuracy, in particular the $\ZDS$ scheme with $R=4$ shows a horizontal line that indicates the boundedness of the Hamiltonian deviation.

\begin{table}[!ht]\notapolice\centering
    \caption{\rempolice Particle in an electromagnetic field, non-separable case from \cref{sec:particle_non_separable}: error on the Hamiltonian at $T=20\,000$ s. For the $\ZD$ and $\ZDS$ schemes, both rows correspond to $N=12\times T$ and $N=48\times T$, respectively; for the classical scheme, they correspond to $N=128\times T$ and $N=256\times T$.}
    \label{tab:Particle_accuracy}
    \begin{subtable}{0.49\textwidth}\centering
        \caption{\rempolice Errors obtained with the $\ZD$ scheme.}
        \label{tab:Particle_ZD_accuracy}

        \begin{tabular}{@{}l@{}  r@{}l@{}  r@{}l@{}  r@{}l@{}  r@{}l@{}}
            \toprule
                             & \phantom{aaa} & {R=2}
                             & \phantom{aaa} & {R=4}
                             & \phantom{aaa} & {R=6}
                             & \phantom{aaa} & {R=8}                                              \\
            \midrule
            \eH ($\times12$) &               & 2.28e-01 &  & 1.31e-02 &  & 6.00e-04 &  & 3.19e-04 \\
            \eH ($\times48$) &               & 8.45e-04 &  & 3.51e-06 &  & 1.11e-08 &  & 6.18e-11 \\
            \bottomrule
        \end{tabular}

    \end{subtable}
    \begin{subtable}{0.49\textwidth}\centering
        \caption{\rempolice Errors obtained with the $\ZDS$ scheme.}
        \label{tab:Particle_ZDS_accuracy}

        \begin{tabular}{@{}l@{}  r@{}l@{}  r@{}l@{}  r@{}l@{}  r@{}l@{}}
            \toprule
                             & \phantom{aaa} & {R=1}
                             & \phantom{aaa} & {R=2}
                             & \phantom{aaa} & {R=3}
                             & \phantom{aaa} & {R=4}                                              \\
            \midrule
            \eH ($\times12$) &               & 1.34e-02 &  & 1.15e-04 &  & 7.59e-07 &  & 6.44e-07 \\
            \eH ($\times48$) &               & 5.29e-05 &  & 2.95e-08 &  & 1.24e-11 &  & 2.64e-14 \\
            \bottomrule
        \end{tabular}

    \end{subtable}

    \medskip

    \begin{subtable}{0.49\textwidth}\centering
        \caption{\rempolice Errors obtained with the classical schemes.}
        \label{tab:Particle_classical_accuracy}

        \begin{tabular}{@{}l@{}  r@{}l@{}  r@{}l@{}  r@{}l@{}  r@{}l@{}}
            \toprule
                              & \phantom{aaa} & {SV2}
                              & \phantom{aaa} & {SV4}
                              & \phantom{aaa} & {SV6}
                              & \phantom{aaa} & {SV8}                                              \\
            \midrule
            \eH ($\times128$) &               & 8.97e-02 &  & 3.82e-04 &  & 5.01e-06 &  & 4.90e-08 \\
            \eH ($\times256$) &               & 4.52e-02 &  & 2.41e-05 &  & 7.88e-08 &  & 1.76e-10 \\
            \bottomrule
        \end{tabular}

    \end{subtable}
\end{table}

\subsection{Complexity and function calls}
\label{sec:complexity}

\newcommand{\neqn}{{\notapolice n\textsubscript{eq}}}
\newcommand{\nit}{{\notapolice n\textsubscript{it}}}
\newcommand{\nitsv}{\smash{\notapolice n$^\text{SV}_\text{it}$}}
\newcommand{\nitzd}{\smash{\notapolice n$^\text{ZD}_\text{it}$}}
\newcommand{\nitzds}{\smash{\notapolice n$^\text{ZDS}_\text{it}$}}
\newcommand{\ts}{\textsuperscript}

In this section, we analyse the complexity and function calls of the structural methods (\ZD\ and \ZDS) compared to classical symplectic schemes. We first provide a detailed complexity analysis for both separable and non-separable cases in \cref{sec:complexity_analysis}. We then present, in \cref{sec:function_calls}, the average number of iterations and function calls required by the fixed-point method used in the structural schemes. This complexity and function call analysis helps to understand the computational efficiency and performance of the proposed methods.

\subsubsection{Complexity analysis for a single time step}
\label{sec:complexity_analysis}

We start with comparing the methods introduced in this article to classical ones, in terms of algorithmic complexity.
We assume that the Hamiltonian system governs the behavior of $K$ bodies in dimension $I$: therefore, the system has $\neqn \coloneqq 2IK$ equations. Moreover, in the structural method (as well as the classical ones in the non-separable case), non-linear equations have to be solved. We denote by $\nit$ the number of iterations of the non-linear solver (a fixed-point method in the structural method, and Newton's method for the classical schemes).
This number of iterations depends on the method, and we denote by \nitsv, \nitzd and \nitzds the number of iterations for the SV, \ZD\ and \ZDS\ methods, respectively. To compute the complexity, we split the problem into two main configurations: the separable and non-separable cases. We assume that evaluating $\Ham{X}{P}$ and its derivatives has complexity $O(1)$.
The complexity of a single time step of each algorithm is reported in \cref{tab_compl1}.

\begin{table}[h!]\notapolice
    \centering
    \renewcommand{\arraystretch}{1.75}
    \setlength{\tabcolsep}{10pt}
    \caption{\rempolice Complexity of the different schemes used and introduced in this paper.}\label{tab_compl1}
    \begin{subtable}{0.49\textwidth}\centering
        \caption{\rempolice separable case}
        \begin{tabular}{cccc}
            \toprule
            \makebox[0pt][c]{Method}
                & \makebox[0pt][c]{order}
                & \makebox[0pt][c]{sub-steps}
                & \makebox[0pt][c]{complexity}
            \\
            \midrule
            MA2 & 2                            & 3                & O(3\neqn)                             \\
            CS4 & 4                            & \rbbis{4}        & O(\rbbis{4}\neqn)                     \\
            KL6 & 6                            & 7 \texttimes\ 3  & O(21\neqn)                            \\
            KL8 & 8                            & 15 \texttimes\ 3 & O(45\neqn)                            \\
            ZD  & R+2                          & N/A              & O\big(2R\neqn\big(1+\nitzd\big)\big)  \\
            ZDS & 2(R+1)                       & N/A              & O\big(4R\neqn\big(1+\nitzds\big)\big) \\
            \bottomrule
        \end{tabular}
    \end{subtable}
    \begin{subtable}{0.49\textwidth}\centering
        \caption{\rempolice non-separable case}
        \begin{tabular}{cccc}
            \toprule
            \makebox[0pt][c]{Method}    &
            \makebox[0pt][c]{order}     &
            \makebox[0pt][c]{sub-steps} &
            \makebox[0pt][c]{complexity}
            \\
            \midrule
            SV2                         & 2      & 2                & O\big(2\neqn\big(1+\nitsv\big)\big)   \\
            SV4                         & 4      & 3 \texttimes\ 2  & O\big(6\neqn\big(1+\nitsv\big)\big)   \\
            SV6                         & 6      & 9 \texttimes\ 2  & O\big(18\neqn\big(1+\nitsv\big)\big)  \\
            SV8                         & 8      & 27 \texttimes\ 2 & O\big(54\neqn\big(1+\nitsv\big)\big)  \\
            ZD                          & R+2    & N/A              & O\big(2R\neqn\big(1+\nitzd\big)\big)  \\
            ZDS                         & 2(R+1) & N/A              & O\big(6R\neqn\big(1+\nitzds\big)\big) \\
            \bottomrule
        \end{tabular}
    \end{subtable}
\end{table}

In the \textbf{separable case}, the complexity of the classical methods
only depends on the order of accuracy and system size.
On the contrary, for the structural method, the complexity also depends on the number of iterations $\nit$ in the fixed-point method, with one additional iteration corresponding to the initialization stage. The \ZDS\ method has twice the complexity of the \ZD\ method, as it requires evaluating the second derivatives of the Hamiltonian in addition to its first derivatives.
Therefore, in terms of complexity, the structural methods are slightly less favorable for low orders, but become comparable or even slightly better for very high orders, depending on the number of iterations. These conclusions are summarized in \cref{tab:complexity_ratio}, where we report the ratio of the complexity of the structural methods with respect to the classical ones.

In the \textbf{non-separable case}, the complexity of the structural and classical methods depend on the order, the system size, and the number of iterations. Indeed, in the SV2 method, one has to solve two non-linear equations and two linear ones per unknown in the system, leading to a complexity in $2\neqn(1+\nit)$. This complexity is then multiplied by a constant depending on the number of compositions, itself depending on the order of the method. For the structural methods, it turns out that the complexity is almost the same as in the separable case. The only difference is that the cross derivatives are now needed in the \ZDS\ method, leading to an extra two evaluations of the Hamiltonian and its derivatives compared to the separable case. \cref{tab:complexity_ratio} once again reports the ratio between the complexity of the structural methods and the classical ones.

\begin{table}[h!]\notapolice
    \centering
    \renewcommand{\arraystretch}{1.75}
    \caption{\rempolice Ratio of the complexity of the \ZD\ and \ZDS\ schemes with the classical ones.}
    \label{tab:complexity_ratio}
    \medskip
    \begin{tabular}{ccccc}
        \toprule
        \multirow{2}{*}{Order}
          & \multicolumn{2}{c}{Separable}
          & \multicolumn{2}{c}{Non-separable}                      \\
        \cmidrule(lr){2-3}
        \cmidrule(lr){4-5}
          & ZD
          & ZDS
          & ZD
          & ZDS                                                    \\
        \cmidrule(lr){1-5}
        4 & \big(1+\nitzd\big) \texttimes\ 4/7
          & \big(1+\nitzds\big) \texttimes\ 4/7
          & \big(1+\nitzd\big)/\big(1+\nitsv\big) \texttimes\ 2/3
          & \big(1+\nitzds\big)/\big(1+\nitsv\big) \texttimes\ 1   \\
        6 & \big(1+\nitzd\big) \texttimes\ 8/21
          & \big(1+\nitzds\big) \texttimes\ 8/21
          & \big(1+\nitzd\big)/\big(1+\nitsv\big) \texttimes\ 4/9
          & \big(1+\nitzds\big)/\big(1+\nitsv\big) \texttimes\ 2/3 \\
        8 & \big(1+\nitzd\big) \texttimes\ 4/15
          & \big(1+\nitzds\big) \texttimes\ 4/15
          & \big(1+\nitzd\big)/\big(1+\nitsv\big) \texttimes\ 2/9
          & \big(1+\nitzds\big)/\big(1+\nitsv\big) \texttimes\ 1/3 \\
        \bottomrule
    \end{tabular}
\end{table}

In summary, in all cases, the gain in complexity by using the structural method depends on the number of iterations of the non-linear solver. In addition, the number of sub-steps increases linearly with the order of accuracy for the structural method, while this increase is super-linear, or even quadratic, for the classical methods. Therefore, to conclude on the potential gains in computation time secured by our schemes, one has to compare the number of iterations in the non-linear solver, as well as the total number of time steps. This is easily done through the number of function calls, whose analysis is presented in the next section.

\subsubsection{Function calls}
\label{sec:function_calls}

On the one hand, the $\ZD$ and $\ZDS$ schemes solve the system using blocks of size $R$, while the total number of steps is denoted by $N$. We simply call \avrOne the ratio between the total number of iterations, including the fixed-point solves, and the total number of steps $N$. Conversely, since one iteration over the block corresponds to $R$ evaluations of the Physical Equations (PE), the average number of ``PE function calls'' \avrTwo is given by the ratio between the total number of calls to PE and the number of steps $N$, that is \texttt{\avrTwo = R \texttimes \avrOne\!\!\!.}
The classical schemes, on the other hand, compose several semi-implicit steps, with one linear and potentially one nonlinear solve for every semi-implicit step. The steps alternate between solving for $\nabla_X \Ham{X}{P}$ and $\nabla_P \Ham{X}{P}$.
In the separable case, we have two linear solves; in the non-separable one, a linear and a non-linear solve are carried out.
Note that, in the specific case of the separable particle problem from \cref{sec:particle_non_separable}, $\nabla_P \Ham{X}{P}$ is linear in $P$, and the non-linear solve for $\nabla_P \Ham{X}{P}$ is not needed, which reduces the total number of non-linear solves.

\begin{table}[!ht]\notapolice\centering
\rbbis{
    \caption{\rempolice Number of iterations of the fixed point method for the outer solar system case, $\ZD$ scheme.  \label{tab_iter_outer_ZD}}
    \begin{subtable}{0.49\textwidth}\centering
        \caption{\rempolice  $N=3000$, \texttt{tol = 10\ts{-30}}.}
        \begin{tabular}{@{}l@{}  r@{}l@{}  r@{}l@{}  r@{}l@{}  r@{}l@{}}
            \toprule
                       & \phantom{aaa} & {R=2}
                       & \phantom{aaa} & {R=4}
                       & \phantom{aaa} & {R=6}
                       & \phantom{aaa} & {R=8}                               \\
            \midrule
            total iter &               & 33405 &  & 18454 &  & 12609 &  & 10103 \\
            \avrOne    &               & 11     &  & 6 &  & 4.2 &  & 3.4  \\
            \avrTwo    &               & 22    &  & 24   &  & 25   &  & 27   \\
            \bottomrule
        \end{tabular}
    \end{subtable}
    \begin{subtable}{0.49\textwidth}\centering
        \caption{\rempolice  $N=12000$, \texttt{tol = 10\ts{-30}}.}
        \begin{tabular}{@{}l@{}  r@{}l@{}  r@{}l@{}  r@{}l@{}  r@{}l@{}}
            \toprule
                       & \phantom{aaa} & {R=2}
                       & \phantom{aaa} & {R=4}
                       & \phantom{aaa} & {R=6}
                       & \phantom{aaa} & {R=8}                                  \\
            \midrule
            total iter &               & 91439 &  & 51000 &  & 34548 &  & 26357 \\
            \avrOne    &               & 7.6   &  & 4.25   &  & 2.88  &  & 2.2  \\
            \avrTwo    &               & 15.2   &  & 17  &  & 17.3  &  & 17.6  \\
            \bottomrule
        \end{tabular}
    \end{subtable}
    \medskip
    \begin{subtable}{0.49\textwidth}\centering
        \caption{\rempolice $N=12000$, \texttt{tol = 10\ts{-15}}.}
        \begin{tabular}{@{}l@{}  r@{}l@{}  r@{}l@{}  r@{}l@{}  r@{}l@{}}
            \toprule
                       & \phantom{aaa} & {R=2}
                       & \phantom{aaa} & {R=4}
                       & \phantom{aaa} & {R=6}
                       & \phantom{aaa} & {R=8}                                 \\
            \midrule
            total iter &               & 47258 &  & 27000 &  & 18617 &  & 15749 \\
            \avrOne    &               & 3.9  &  & 2.2   &  & 1.6   &  & 1.3  \\
            \avrTwo    &               & 7.8  &  & 9  &  & 9.6  &  & 10.4 \\
            \bottomrule
        \end{tabular}
    \end{subtable}
}   
\end{table}

\begin{table}[!ht]\notapolice\centering
\rbbis{
    \caption{\rempolice Number of iterations of the fixed point method for the outer solar system case, $\ZDS$ scheme. \label{tab_iter_outer_ZDS}}
    \begin{subtable}{0.49\textwidth}\centering
        \caption{\rempolice  $N=3000$, \texttt{tol = 10\ts{-30}}.}
        \begin{tabular}{@{}l@{}  r@{}l@{}  r@{}l@{}  r@{}l@{}  r@{}l@{}}
            \toprule
                       & \phantom{aaa} & {R=1}
                       & \phantom{aaa} & {R=2}
                       & \phantom{aaa} & {R=3}
                       & \phantom{aaa} & {R=4}                                  \\
            \midrule
            total iter &               & 60000 &  & 33591 &  & 25434 &  & 21719 \\
            \avrOne    &               & 20.0  &  & 11.2  &  & 8.5 &  & 8     \\
            \avrTwo    &               & 20.0  &  & 22.4  &  & 25.4  &  & 29  \\
            \bottomrule
        \end{tabular}
    \end{subtable}
    \begin{subtable}{0.49\textwidth}\centering
        \caption{\rempolice  $N=12000$, \texttt{tol = 10\ts{-30}}.}
        \begin{tabular}{@{}l@{}  r@{}l@{}  r@{}l@{}  r@{}l@{}  r@{}l@{}}
            \toprule
                       & \phantom{aaa} & {R=1}
                       & \phantom{aaa} & {R=2}
                       & \phantom{aaa} & {R=3}
                       & \phantom{aaa} & {R=4}                                  \\
            \midrule
            total iter &               & 168000 &  & 91169 &  & 64000 &  & 53994 \\
            \avrOne    &               &  14  &  & 7.6   &  & 5.3  &  & 4.5   \\
            \avrTwo    &               &  14  &  & 15.2   &  & 16   & & 18   \\
            \bottomrule
        \end{tabular}
    \end{subtable}
    \medskip
    \begin{subtable}{0.49\textwidth}\centering
        \caption{\rempolice $N=12000$, \texttt{tol = 10\ts{-15}}.}
        \begin{tabular}{@{}l@{}  r@{}l@{}  r@{}l@{}  r@{}l@{}  r@{}l@{}}
            \toprule
                       & \phantom{aaa} & {R=1}
                       & \phantom{aaa} & {R=2}
                       & \phantom{aaa} & {R=3}
                       & \phantom{aaa} & {R=4}                                  \\
            \midrule
            total iter &               & 84000 &  & 44891 &  & 32000 &  & 24000   \\
            \avrOne    &               & 7  &  & 3.7  &  & 2.67 &  & 2 \\
            \avrTwo    &               & 7  &  & 7.5  &  & 8    &  & 8  \\
            \bottomrule
        \end{tabular}
    \end{subtable}
}   
\end{table}

\begin{table}[!ht]\notapolice\centering
    \caption{\rempolice Number of iterations of the fixed point method for the particle in electromagnetic field case, $\ZD$ scheme.  \label{tab_iter_particule_ZD}}
    \begin{subtable}{0.49\textwidth}\centering
        \caption{\rempolice  $N=12000$, \texttt{tol = 10\ts{-30}}.}

        \begin{tabular}{@{}l@{}  r@{}l@{}  r@{}l@{}  r@{}l@{}  r@{}l@{}}
            \toprule
                       & \phantom{aaa} & {R=2}
                       & \phantom{aaa} & {R=4}
                       & \phantom{aaa} & {R=6}
                       & \phantom{aaa} & {R=8}                                    \\
            \midrule
            total iter &               & 188656 &  & 107441 &  & 77665 &  & 62836 \\
            \avrOne    &               & 15.7   &  & 8.9    &  & 6.5   &  & 5.2   \\
            \avrTwo    &               & 31.4   &  & 35.9   &  & 39.0  &  & 42.0  \\
            \bottomrule
        \end{tabular}

    \end{subtable}
    \begin{subtable}{0.49\textwidth}\centering
        \caption{\rempolice  $N=48000$, \texttt{tol = 10\ts{-30}}.}

        \begin{tabular}{@{}l@{}  r@{}l@{}  r@{}l@{}  r@{}l@{}  r@{}l@{}}
            \toprule
                       & \phantom{aaa} & {R=2}
                       & \phantom{aaa} & {R=4}
                       & \phantom{aaa} & {R=6}
                       & \phantom{aaa} & {R=8}                                      \\
            \midrule
            total iter &               & 473032 &  & 256917 &  & 182448 &  & 142823 \\
            \avrOne    &               & 9.6    &  & 5.4    &  & 3.8    &  & 3.0    \\
            \avrTwo    &               & 19.7   &  & 21.4   &  & 22.8   &  & 23.8   \\
            \bottomrule
        \end{tabular}

    \end{subtable}

    \medskip

    \begin{subtable}{0.49\textwidth}\centering
        \caption{\rempolice $N=48000$, \texttt{tol = 10\ts{-15}}.}

        \begin{tabular}{@{}l@{}  r@{}l@{}  r@{}l@{}  r@{}l@{}  r@{}l@{}}
            \toprule
                       & \phantom{aaa} & {R=2}
                       & \phantom{aaa} & {R=4}
                       & \phantom{aaa} & {R=6}
                       & \phantom{aaa} & {R=8}                                  \\
            \midrule
            total iter &               & 98541 &  & 57103 &  & 41795 &  & 34031 \\
            \avrOne    &               & 2.1   &  & 1.2   &  & 0.87  &  & 0.71  \\
            \avrTwo    &               & 4.1   &  & 4.8   &  & 5.2   &  & 5.7   \\
            \bottomrule
        \end{tabular}

    \end{subtable}
\end{table}

\begin{table}[!ht]\notapolice\centering
    \caption{\rempolice Number of iterations of the fixed point method for the particle in electromagnetic field case, $\ZDS$ scheme. \label{tab_iter_particule_ZDS}}
    \begin{subtable}{0.49\textwidth}\centering
        \caption{\rempolice  $N=12000$, \texttt{tol = 10\ts{-30}}.}

        \begin{tabular}{@{}l@{}  r@{}l@{}  r@{}l@{}  r@{}l@{}  r@{}l@{}}
            \toprule
                       & \phantom{aaa} & {R=1}
                       & \phantom{aaa} & {R=2}
                       & \phantom{aaa} & {R=3}
                       & \phantom{aaa} & {R=4}                                      \\
            \midrule
            total iter &               & 340494 &  & 202183 &  & 161918 &  & 152791 \\
            \avrOne    &               & 28.4   &  & 16.8   &  & 13.5   &  & 12.7   \\
            \avrTwo    &               & 28.4   &  & 33.7   &  & 40.5   &  & 51.0   \\
            \bottomrule
        \end{tabular}

    \end{subtable}
    \begin{subtable}{0.49\textwidth}\centering
        \caption{\rempolice  $N=48000$, \texttt{tol = 10\ts{-30}}.}

        \begin{tabular}{@{}l@{}  r@{}l@{}  r@{}l@{}  r@{}l@{}  r@{}l@{}}
            \toprule
                       & \phantom{aaa} & {R=1}
                       & \phantom{aaa} & {R=2}
                       & \phantom{aaa} & {R=3}
                       & \phantom{aaa} & {R=4}                                      \\
            \midrule
            total iter &               & 856330 &  & 474821 &  & 349537 &  & 292265 \\
            \avrOne    &               & 17.8   &  & 9.9    &  & 7.3    &  & 6.1    \\
            \avrTwo    &               & 17.8   &  & 19.8   &  & 21.8   &  & 24.4   \\
            \bottomrule
        \end{tabular}

    \end{subtable}

    \medskip

    \begin{subtable}{0.49\textwidth}\centering
        \caption{\rempolice $N=48000$, \texttt{tol = 10\ts{-15}}.}

        \begin{tabular}{@{}l@{}  r@{}l@{}  r@{}l@{}  r@{}l@{}  r@{}l@{}}
            \toprule
                       & \phantom{aaa} & {R=1}
                       & \phantom{aaa} & {R=2}
                       & \phantom{aaa} & {R=3}
                       & \phantom{aaa} & {R=4}                                      \\
            \midrule
            total iter &               & 422185 &  & 238998 &  & 172695 &  & 143948 \\
            \avrOne    &               & 8.8    &  & 5.0    &  & 3.4    &  & 3.0    \\
            \avrTwo    &               & 8.8    &  & 10     &  & 10.8   &  & 12.0   \\
            \bottomrule
        \end{tabular}

    \end{subtable}
\end{table}

We report in \cref{tab_iter_outer_ZD,tab_iter_outer_ZDS}, the total number of iterations of the fixed point method (including the initialization), the average number of iterations and call to PE for the outer solar system (a separable example from \cref{sec:solar_system}) and in \cref{tab_iter_particule_ZD,tab_iter_particule_ZDS},
for the particle in an electromagnetic field (a non-separable example from \cref{sec:3D_EM_field}), respectively. The number of iterations decreases as we increase the value of $R$, since we handle larger blocks. We then observe that the \avrOne is also decreasing in both simulations. Moreover, the ratio also decreases for larger~$N$ since we get a better predictor for the fixed point method; hence we reduce the number of iterations to reach the tolerance value. If one releases the constraint on the tolerance to \texttt{10\ts{-15}} in the last case, we almost divide the number of iterations by two. We run faster but less \ra{accurately}; therefore the user has to determine the best trade-off between computational effort and solution quality.

Additional remarks concern the outer solar system when we take $N=12000$ instead of $N=3000$. We observe that we have a substantial gain since the iteration ratio strongly diminishes with a larger number of steps. This comes from a better initialization of the fixed point method, making a smaller number of iterations with a more efficient predictor. We highlight that the proposed initialization based on a simple Taylor expansion could be strongly improved by taking advantage of the former values of the function and derivative \ra{computed} in the previous $R$ block.

\begin{table}[!ht]\notapolice\centering
    \caption{\rempolice Number of time steps, function calls and Newton iterations for the classical schemes applied to the particle in electromagnetic field case. On average, $3$ Newton iterations are used at every time step, leading to $6$ calls to the gradient of the Hamiltonian and $3$ calls to its second derivatives per time step.}
    \label{tab_iter_particule_CL}
    \begin{subtable}{\textwidth}\centering
        \caption{\rempolice  $N=2\,560\,000$, \texttt{tol = 10\ts{-15}}}

        \begin{tabular}{@{}l@{}  r@{}l@{}  r@{}l@{}  r@{}l@{}  r@{}l@{}}
            \toprule
                                 & \phantom{aaa} & {SV2}
                                 & \phantom{aaa} & {SV4}
                                 & \phantom{aaa} & {SV6}
                                 & \phantom{aaa} & {SV8}                                                                \\
            \midrule
            \# time steps        &               & 2\,560\,000  &  & 2\,560\,000  &  & 2\,560\,000   &  & 2\,560\,000   \\
            \# function calls    &               & 23\,039\,068 &  & 69\,119\,694 &  & 207\,359\,500 &  & 622\,079\,462 \\
            \# Newton iterations &               & 7\,679\,534  &  & 23\,039\,847 &  & 69\,119\,750  &  & 207\,359\,731 \\
            \bottomrule
        \end{tabular}

    \end{subtable}

    \medskip

    \begin{subtable}{\textwidth}\centering
        \caption{\rempolice  $N=5\,120\,000$, \texttt{tol = 10\ts{-15}}}

        \begin{tabular}{@{}l@{}  r@{}l@{}  r@{}l@{}  r@{}l@{}  r@{}l@{}}
            \toprule
                                 & \phantom{aaa} & {SV2}
                                 & \phantom{aaa} & {SV4}
                                 & \phantom{aaa} & {SV6}
                                 & \phantom{aaa} & {SV8}                                                                    \\
            \midrule
            \# time steps        &               & 5\,120\,000  &  & 5\,120\,000   &  & 5\,120\,000   &  & 5\,120\,000      \\
            \# function calls    &               & 46\,034\,006 &  & 138\,218\,150 &  & 414\,695\,588 &  & 1\,244\,123\,488 \\
            \# Newton iterations &               & 15\,337\,003 &  & 46\,069\,075  &  & 138\,227\,794 &  & 414\,701\,744    \\
            \bottomrule
        \end{tabular}

    \end{subtable}
\end{table}

To compare the structural method to the classical ones, we provide in \cref{tab_iter_particule_CL} the number of time steps, function calls and Newton iterations for the classical schemes, in the case of the particle in an electromagnetic field, and with double-precision accuracy (i.e., \texttt{tol = 10\ts{-15}}). For the second-order SV2 schemes, in this case, three explicit sub-steps and one implicit sub-step are carried out for each time step. With, on average, three Newton iterations per implicit sub-step, we obtain, on average, nine function calls per time step: $3$ for the explicit steps, $3$ for the gradient in the Newton iterations, and $3$ for the Hessian in the Newton iteration, assuming that computing the Hessian is as expensive as computing the gradient, which is a very conservative assumption. This number of function calls per time step is on par with the \ZD\ and \ZDS\ schemes for \texttt{tol = 10\ts{-15}}. However, since the classical schemes require much more time steps in order to ensure stability, it turns out that the total number of function calls is much higher for the classical schemes than for the structural ones.

Based on the preceding analysis, we can summarize the key conclusions as follows:
\begin{itemize}
    \item One of the key benefits of the structural method is its unconditional stability.
          Indeed, one can take very large time steps, which leads to a significant reduction in the number of function calls without loss in accuracy.
    \item In the separable case, both approaches allow for large time steps.
          For high orders of accuracy and large number~$N$ of time steps, the structural schemes are faster because \nitzd\ and \nitzds\ remain small. At lower orders and small values of ~$N$, the structural method is slower, but its precision remains superior. As a result, we maintain a better speed-to-accuracy ratio overall.
    \item In the non-separable case, regardless of $N$, our method has a smaller complexity if a similar amount of iterations in the non-linear solver is required. Even when the number of iterations is larger, the structural method still has a significant advantage, namely its unconditional stability. In the particle case (where the number of iterations was comparable), the classical method required time steps about $20$ times smaller than the structural method to ensure stability.
\end{itemize}
\rb {In conclusion, the structural method should, in general, be slightly faster than classical methods in the separable case, and significantly faster in the non-separable case. To make sure this is truly the case, in the next section, we compare the runtime of both methods in the two cases.}

\subsection{Computation time}
\label{sec:computation_time}

\rb{This last section is dedicated to the comparison of the computation time required by the structural schemes and the classical ones. We tackle a separable setting (the outer Solar system problem) and a non-separable one (the particle in a 3D electromagnetic field). To provide a fair comparison, they were coded in the same framework, using the same representation of floating-point numbers. Since the ZD and ZDS schemes are more stable, we additionally choose the lowest time step imposed by the classical symplectic schemes. The tests were run on a personal computer with an AMD Ryzen 7 5800H @ 4.46 GHz and 16 GB of RAM with only one thread.}

\begin{figure}[!ht]
    \centering
    \includegraphics{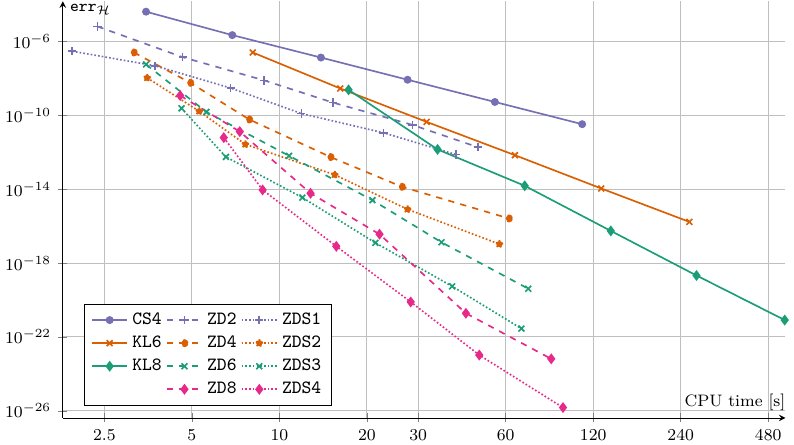}
    \caption{\rempolice \rbbis{Hamiltonian error} with respect to computation time for the outer solar system problem from \cref{sec:solar_system}. The structural schemes ($\ZD$ and $\ZDS$) are compared to the classical symplectic ones.}
    \label{fig:cpu_time_outer_planets}
\end{figure}

\rb{In \cref{fig:cpu_time_outer_planets}, we display the efficiency curves (error with respect to computation time) for the (separable) outer planets benchmark, until final time $T=10^6$. In all cases, we observe that the structural scheme have an enormous gain in CPU time compared to the classical symplectic ones. For instance, if one wishes to pay a budget of around $30$ seconds of computation time with a scheme of order $8$, one can reach an error of around $10^{-12}$ with the classical schemes, while the structural schemes provide an error of around $10^{-16}$ (for $\ZD$) and $10^{-19}$ (for $\ZDS$), which is a gain of $7$ orders of magnitude in accuracy for the same computational effort. To this gain, one could add the larger time steps allowed by the structural schemes, which would further increase the advantage of our method.}

\begin{figure}[!ht]
    \centering
    \includegraphics{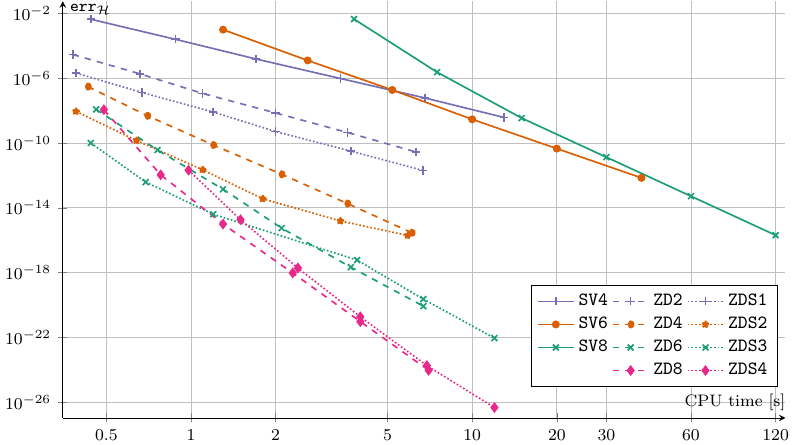}
    \caption{\rempolice \rbbis{Hamiltonian error} with respect to computation time for the particle in a 3D electromagnetic field problem from \cref{sec:3D_EM_field}. The structural schemes ($\ZD$ and $\ZDS$) are compared to the classical symplectic ones.}
    \label{fig:cpu_time_particle}
\end{figure}

\rb{In \cref{fig:cpu_time_particle}, we once again display efficiency curves, this time for the (non-separable) particle in a 3D electromagnetic field problem, using the benchmark from \cref{sec:particle_non_separable}, which we ran until final time $T=100$. We again observe dramatic gains in efficiency with the structural schemes. For instance, if one's computational budget is $5$ seconds, one can reach an error of around $10^{-7}$ with the classical schemes, while both structural schemes provide an error of around $10^{-19}$. This is a gain of~$12$ orders of magnitude in accuracy for the same computational budget. This verifies the observation from \cref{sec:complexity} that the structural schemes are more efficient in the non-separable case.}

\section{Conclusions}
In this paper, evidence was brought forward of the ability of the structural method to provide stable and accurate solutions in the specific context of Hamiltonian systems. We have detailed the design of the method, and provided numerical experiments with different classical Hamiltonian systems, systematically comparing the results to well-known symplectic solvers. These experiments showed that our method is quite efficient and fast, and, in the $\ZDS$ case, provides a solution of excellent quality. We also tackled a non-separable system to prove the superiority of the structural method with \ra{respect} to other schemes. In short, we have provided a general framework for the numerical approximation of Hamiltonian systems, that can easily be adapted to large classes of problems.

We would like to highlight that our method is quite versatile (e.g. by adapting the block size to fix the convergence rate), very easy to implement, and with a low memory consumption. Unconditional stability is also a desirable property, since it makes it possible to reduce the computational effort involved in dealing with long-time simulations by taking large time steps. Additionally, we also mention the simplicity of the fixed point method to solve non-linear problems where all the matrices and inverse matrices associated to the structural equations are computed once in a pre-processing stage, that dramatically reduces runtime. At last, we achieve very high accuracy that requires quad-precision (even octuple precision) when we use the $\ZDS$ scheme with $R>2$.

A final remark is the ability to produce more physical equations by differentiating the physical equations twice (or more), \ra{and by including} the third-order (or higher-order) derivatives in the structural equations. That way, extreme accuracy (more than one hundred correct digits) would be obtained, which could be useful for applications requiring very long simulation times and very accurate approximations (planet or satellite positioning, for instance).

\subsection*{Funding}
S. Clain  was financially supported by the Fundação para a Ciência e a Tecnologia (Portuguese Foundation for Science and Technology) under the scope of the projects UID/00324/2025 \\ {\tt https://doi.org/10.54499/UID/00324/2025} (Centre for Mathematics of the University of Coimbra).\\[0.4em]
S. Clain acknowledges the financial support of the Portuguese Foundation for Science and Technology (FCT) via a national funding for projects IC\&DT with the reference 2023.16854.ICDT.\\[0.4em]
The authors extend their thanks to ANR-24-CE46-7505 SMEAGOL.\\[0.4em]
V. Michel-Dansac thanks ANR-22-CE25-0017 OptiTrust. \\[0.4em]
The Shark-FV conference has greatly contributed to this work.

\bibliographystyle{plain}
\bibliography{references}

\appendix

\section{Analytic expressions for the structural schemes}

For $R \in \{1, 2, 3, 4\}$, we provide in \cref{table::SE_equation_ZD_analytic} the analytic expressions of the matrices $A_z$ and $A_d$, and vectors $\ba_z$ and $\ba_d$ involved in the \ZD\ scheme from \eqref{eq:ZD_scheme_matrix_form}. Similarly, in \cref{table::SE_equation_ZDS_analytic}, we provide the analytic expressions of the matrices $A_z$, $A_d$ and $A_s$, and vectors $\ba_z$, $\ba_d$ and $\ba_s$ involved in the \ZDS\ scheme~\eqref{eq:ZDS_scheme_matrix_form}.

\begin{table}[!ht]
    \tablepolice
    \begin{tabular}{cl}
        \toprule
        $R=1$: \quad & \scalebox{0.9}{
                           $A_z = 2$, $A_d = -1$, $\ba_z = -2$, $\ba_d = -1$
                       }                                                        \\
        \addlinespace[0.75em]
        $R=2$: \quad & \scalebox{0.9}{
                           $A_z = \begin{pmatrix}
                                   0 & 3 \\ -4 & 2
                               \end{pmatrix}$,
                           $A_d = \begin{pmatrix}
                                   -4 & -1 \\ 0 & -1
                               \end{pmatrix}$,
                           $\ba_z = \begin{pmatrix}
                                   3 \\ -2
                               \end{pmatrix}$,
                           $\ba_d = \begin{pmatrix}
                                   1 \\ -1
                               \end{pmatrix}$
                       }                                                                             \\
        \addlinespace[0.75em]
        $R=3$: \quad & \scalebox{0.9}{
                           $A_z = \begin{pmatrix}
                                   27 & -27 & -11 \\ 0 & -27 & 8 \\ 27 & -27 & 13
                               \end{pmatrix}$,
                           $A_d = \begin{pmatrix}
                                   27 & 27 & 3 \\ 27 & 0 & -3 \\ 0 & 0 & -6
                               \end{pmatrix}$,
                           $\ba_z = \begin{pmatrix}
                                   -11 \\ -19 \\ 13
                               \end{pmatrix}$,
                           $\ba_d = \begin{pmatrix}
                                   -3 \\ -6 \\ 6
                               \end{pmatrix}$
                       }                                                   \\
        \addlinespace[0.75em]
        $R=4$: \quad & \scalebox{0.9}{\makecell[l]{
                               $A_z = \begin{pmatrix}
                                       160 & 0 & -160 & -25 \\ -224 & 216 & 96 & -19 \\ 16 & -108 & 48 & -16 \\ 32 & -36 & 32 & -14
                                   \end{pmatrix}$,
                               $\ba_z = \begin{pmatrix}
                                       25 \\ -69 \\ 60 \\ -14
                                   \end{pmatrix}$, \\
                               $A_d = \begin{pmatrix}
                                       96 & 216 & 96 & 6 \\ -192 & -216 & 0 & 6\\ 96 & 0 & 0 & 6 \\ 0 & 0 & 0 & 6
                                   \end{pmatrix}$,
                               $\ba_d = \begin{pmatrix}
                                       6 \\ -18 \\ 18 \\ -6
                                   \end{pmatrix}$,
                           }}
        \\ \bottomrule
    \end{tabular}
    \caption{\rempolice\ZD\ scheme from \eqref{eq:ZD_scheme_matrix_form}: Expressions of the matrices $A_z$ and $A_d$, and vectors $\ba_z$ and $\ba_d$, for $R \in \{1, 2, 3, 4\}$.}
    \label{table::SE_equation_ZD_analytic}
\end{table}

\begin{table}[!ht]
    \tablepolice
    \begin{tabular}{cl}
        \toprule
        $R=1$: \quad & \scalebox{0.9}{
                           $A_z = 12$, $A_d = -6$, $A_s = 1$,
                           $\ba_z = -12$, $\ba_d = -6$, $\ba_s = -1$
                       }                            \\
        \addlinespace[1.5em]
        $R=2$: \quad & \scalebox{0.9}{\makecell[l]{
                               $A_z = \begin{pmatrix}
                                       -48 & 24 \\ 0 & 15
                                   \end{pmatrix}$,
                               $A_d = \begin{pmatrix}
                                       0 & -9 \\ -16 & -7
                                   \end{pmatrix}$,
                               $A_s = \begin{pmatrix}
                                       -8 & 1 \\ 0 & 1
                                   \end{pmatrix}$, \\
                               $\ba_z = \begin{pmatrix}
                                       -24 \\ 15
                                   \end{pmatrix}$,
                               $\ba_d = \begin{pmatrix}
                                       -9 \\ 7
                                   \end{pmatrix}$,
                               $\ba_s = \begin{pmatrix}
                                       -1 \\ 1
                                   \end{pmatrix}$
                           }}                                     \\
        \addlinespace[1.5em]
        $R=3$: \quad & \scalebox{0.9}{\makecell[l]{
                               $A_z = \begin{pmatrix}
                                       729 & -729 & 103 \\ -486 & 243 & 70 \\ -81 & 81 & 113
                                   \end{pmatrix}$,
                               $A_d = \begin{pmatrix}
                                       243 & 243 & -33 \\ -81 & -162 & -27 \\ -162 & -162 & -48
                                   \end{pmatrix}$,
                               $A_s = \begin{pmatrix}
                                       81 & -81 & 3 \\ -81 & 0 & 3 \\ 0 & 0 & 6
                                   \end{pmatrix}$, \\
                               $\ba_z = \begin{pmatrix}
                                       103 \\ -173 \\ 113
                                   \end{pmatrix}$,
                               $\ba_d = \begin{pmatrix}
                                       33 \\ -60 \\ 48
                                   \end{pmatrix}$,
                               $\ba_s = \begin{pmatrix}
                                       3 \\ -6 \\ 6
                                   \end{pmatrix}$
                           }} \\
        \addlinespace[1.5em]
        $R=4$: \quad & \scalebox{0.9}{\makecell[l]{
                               $A_z = \begin{pmatrix}
                                       -10240 & 19440  & -10240 & 520 \\
                                       16640  & -19440 & 3840   & 370 \\
                                       -5248  & 2592   & 1152   & 307 \\
                                       -512   & 0      & 512    & 269
                                   \end{pmatrix}$,
                               $\ba_z = \begin{pmatrix}
                                       520 \\-1410\\1197\\-269
                                   \end{pmatrix}$, \\
                               $A_d = \begin{pmatrix}
                                       -3840 & 0     & 3840  & -150 \\
                                       6144  & 5184  & -1536 & -126 \\
                                       -1536 & -2592 & -768  & -114 \\
                                       -512  & -864  & -512  & -106
                                   \end{pmatrix}$,
                               $\ba_d = \begin{pmatrix}
                                       150 \\ -426\\ 390\\ -106
                                   \end{pmatrix}$,\\
                               $A_s = \begin{pmatrix}
                                       -768 & 2592  & -768 & 12 \\
                                       1536 & -2592 & 0    & 12 \\
                                       -768 & 0     & 0    & 12 \\
                                       0    & 0     & 0    & 12
                                   \end{pmatrix}$,
                               $\ba_s = \begin{pmatrix}
                                       12 \\-36\\36\\-12
                                   \end{pmatrix}$
                           }}                          \\  \bottomrule
    \end{tabular}
    \caption{\rempolice\ZDS\ scheme from \eqref{eq:ZDS_scheme_matrix_form}: Expressions of the matrices $A_z$, $A_d$ and $A_s$, and vectors $\ba_z$, $\ba_d$ and $\ba_s$, for $R \in \{1, 2, 3, 4\}$.}
    \label{table::SE_equation_ZDS_analytic}
\end{table}

\section{Exact conservation of non-separable quadratic Hamiltonians}
\label{sec:non_separable_quadratic_Hamiltonians}

We consider a non-separable quadratic Hamiltonian of the form
\eqref{eq:general_quadratic_Hamiltonian}, i.e.,
\begin{equation*}
    \Ham{x}{p} = \frac 1 2 p^2 + \gamma x p + \frac 1 2 \omega^2 x^2.
\end{equation*}
In this section, we exhibit conditions on the $R$-block \ZD\ scheme so that it
exactly preserves this Hamiltonian at step $1 \leq r \leq R$, i.e.,
\begin{equation*}
    \Ham{\Zx_{n+r}}{\Zp_{n+r}} = \Ham{\Zx_n}{\Zp_n}.
\end{equation*}
The proof for the \ZDS\ scheme is carried out in the same manner,
and it is left to the reader.
To exhibit the conditions, we follow the same steps as in
the proof of \propref{prop:energy_conservation_ZD}.
Indeed, plugging the physical equations \eqref{PE1x}--\eqref{PE1p}
into the structural equations \eqref{ZD_SEx}--\eqref{ZD_SEp} yields
\begin{align}
    \label{eq:ZD_SEx_appendix}
    \blockZx_n
    + \gamma \Delta t \, B_d\, \blockZx_n
    + \Delta t \, B_d\,\blockZp_n
     & =
    - \Zx_n\, b_z
    - \gamma \Delta t \, \Zx_n\, b_d
    - \Delta t \, \Zp_n\, b_d, \\
    \label{eq:ZD_SEp_appendix}
    \blockZp_n
    - \gamma \Delta t \, B_d\,\blockZp_n
    - \omega^2 \Delta t \, B_d\,\blockZx_n
     & =
    - \Zp_n\, b_z
    + \gamma \Delta t \, \Zp_n\, b_d
    + \omega^2 \Delta t \, \Zx_n\, b_d.
\end{align}
Isolating $\blockZx_n$ in \eqref{eq:ZD_SEx_appendix}
and $\blockZp_n$ in \eqref{eq:ZD_SEp_appendix}, we obtain
\begin{align}
    \label{eq:ZD_SEx_appendix_2}
    (I_R + \gamma \Delta t \, B_d\, \blockZx_n) \, \blockZx_n
     & =
    - \Delta t \, B_d\,\blockZp_n
    - \Zx_n\, (b_z + \gamma \Delta t \, b_d)
    - \Delta t \, \Zp_n\, b_d, \\
    \label{eq:ZD_SEp_appendix_2}
    (I_R - \gamma \Delta t \, B_d\, \blockZx_n) \, \blockZp_n
     & =
    \omega^2 \Delta t \, B_d\,\blockZx_n
    - \Zp_n\, (b_z - \gamma \Delta t \, b_d)
    + \omega^2 \Delta t \, \Zx_n\, b_d.
\end{align}
Plugging \eqref{eq:ZD_SEx_appendix_2} into \eqref{eq:ZD_SEp_appendix}
and \eqref{eq:ZD_SEp_appendix_2} into \eqref{eq:ZD_SEx_appendix},
rearranging terms, and performing some straightforward
(but tedious) algebraic manipulations,
we obtain that the $R$-blocks satisfy
\begin{equation}
    \label{eq:ZD_SE_blocks_final}
    A_+ \, \blockZx_n = \Zx_n \, a_z^+ + \Delta t\, \Zp_n \, a_d^+
    \text{\qquad and \qquad}
    A_- \, \blockZp_n = \Zp_n \, a_z^- - \omega^2 \Delta t\, \Zx_n \, a_d^-,
\end{equation}
where we have set \begin{align*}
    A_\pm   & =
    M_\pm + \omega^2 \Delta t^2 \, B_d M_\pm^{-1} B_d,                            \\
    a_z^\pm & =
    - b_z \mp \gamma \Delta t \, b_d - \omega^2 \Delta t^2 \, B_d M_\pm^{-1} b_d, \\
    a_d^\pm & = -b_d + B_d M_\pm^{-1} (b_z \pm \gamma \Delta t \, b_d),           \\
    M_\pm   & = I_R \pm \gamma \Delta t \, B_d.
\end{align*}
This gives eaasily checkable conditions for the preservation of the Hamiltonian,
by simply computing the $r$\textsuperscript{th} entry of the vector
\begin{equation*}
    \frac 1 2 \blockZp_n \odot \blockZp_n
    + \gamma \blockZx_n \odot \blockZp_n
    + \frac 1 2 \omega^2 \blockZx_n \odot \blockZp_n
    - \left(
    \frac 1 2 \Zp_n^2
    + \gamma \Zx_n \Zp_n
    + \frac 1 2 \omega^2 \Zx_n^2
    \right),
\end{equation*}
where $\blockZx_n$ and $\blockZp_n$ are given by \eqref{eq:ZD_SE_blocks_final},
and $\odot$ denotes the Hadamard (componentwise) product of two vectors,
and checking whether it vanishes.

\end{document}